\documentclass[final]{siamltex}

\usepackage{url}
\usepackage{graphicx}
\usepackage{epstopdf}
\usepackage{amsmath,amsfonts,amssymb}
\usepackage[usenames]{color}
\usepackage[mathscr]{eucal}

\usepackage{stmaryrd}

\usepackage{amsbsy}

\usepackage{bm}

\usepackage{xspace}

\usepackage[margin=5pt,justification=centering]{subfig}

\usepackage{paralist}
\usepackage{braket}
\usepackage{tabulary}

\usepackage{algorithm}
\usepackage{algpseudocode}

\newcommand{\MI}[1]{\mathbf{#1}}

\newcommand{\V}[1]{{\bm{\mathbf{\MakeLowercase{#1}}}}} % vector
\newcommand{\Vhat}[1]{{\bm{\hat \mathbf{\MakeLowercase{#1}}}}} % vector
\newcommand{\M}[1]{{\bm{\mathbf{\MakeUppercase{#1}}}}} % matrix
\newcommand{\Mtilde}[1]{{\bm{\tilde \mathbf{\MakeUppercase{#1}}}}} % matrix
\newcommand{\Mn}[2]{\M{#1}^{(#2)}} % n-th matrix
\newcommand{\T}[1]{\boldsymbol{\mathscr{\MakeUppercase{#1}}}} %tensor
\newcommand{\TE}[2]{\MakeLowercase{#1}_{\MI{#2}}} % tensor element with multi-index

\newcommand{\Mz}[2]{\M{#1}_{(#2)}} % n-th mode matricize
\newcommand{\KT}[1]{\left\llbracket #1 \right\rrbracket} % kruskal operator
\newcommand{\KTsmall}[1]{\llbracket #1 \rrbracket} % small kruskal operator
\def\TX{\T{X}}

\def\TM{\T{M}} %%%%%%%%%%%%%

\def\Vl{\V{\lambda}}

\def\rowx{\hat{\V{x}}}
\def\kktviol{\textrm{kkt}_{\textrm{viol}}}

\newcommand{\ib}{\textbf{i}}

\newcommand{\eb}{\textbf{e}}

\def\RefSecIntro{Section~1.1 }
\def\RefSecExp{Section~4 }

\def\RefSecExpFullprobPdnrPqnrMu{Section~4.2.2 }
\def\RefSecAppendix{Appendix~A }
\def\RefTblSparsity{Table~1 }
\def\RefTblFullprob{Table~6 }
\def\RefBibScoreK{[1]}
\def\RefBibPhanZdunek2010{Phan et al.~[33]}
\def\RefBibZdunekCichocki{Zdunek and Cichocki~[42] }
\def\RefBibZdunekCichockiNumber{[42] }

\title{Newton-Based Optimization for Kullback-Leibler
       Nonnegative Tensor Factorizations}
\author{Samantha Hansen, Todd Plantenga, Tamara G. Kolda}

\begin{document}
\maketitle

\begin{abstract}
Tensor factorizations with nonnegative constraints have found application
in analyzing data from cyber traffic, social networks, and other areas.
We consider application data best described
as being generated by a Poisson process (e.g., count data),
which leads to sparse tensors that can be modeled
by sparse factor matrices.
In this paper we investigate efficient techniques for computing an appropriate
canonical polyadic tensor factorization based on the Kullback-Leibler divergence function.
We propose novel subproblem solvers within the
standard alternating block variable approach.
Our new methods exploit structure and reformulate the optimization problem
as small independent subproblems.  We employ bound-constrained Newton
and quasi-Newton methods.
We compare our algorithms against other codes, demonstrating superior speed
for high accuracy results and the ability to quickly find sparse solutions.
\end{abstract}

\section{Introduction}

Multilinear models have proved useful in analyzing data in a variety of fields.
We focus on data that derives from a Poisson process,
such as the number of packets sent from one IP address to another
on a specific port~\cite{Sun:DynamicTensorCyber}, the number of papers
published by an author at a given conference~\cite{Kolda:DBLP}, or
the count of emails between users in a given time period~\cite{Bader:Enron2007}.
Data in these applications is nonnegative and often quite sparse, i.e.,
most tensor elements have a count of zero.
The tensor factorization model corresponding to such sparse count data
is computed from a nonlinear optimization problem that minimizes the
Kullback-Leibler (K-L) divergence function and contains nonnegativity
constraints on all variables.

In this paper we 
show how to make second-order optimization
methods suitable for Poisson-based tensor models of large sparse count data.
Multiplicative update
is one of the most widely implemented methods for this model,
but it suffers from slow convergence and inaccuracy in discovering the underlying
sparsity.  In large sparse tensors, the application of nonlinear optimization
techniques requires consideration of sparsity and problem structure
to get better performance.
We show that, by exploiting the partial separability of the subproblems,
we can successfully apply second-order methods.
We develop algorithms that scale to large sparse tensor applications
and are quick in identifying sparsity in the factors of the model.

There is a need for
second-order methods because computing factor matrices to
high accuracy, as measured by satisfaction
of the first-order KKT conditions, is effective in revealing sparsity.
By contrast, multiplicative update methods can make elements small
but are slow to reach the variable bound at zero, forcing the
user to guess when ``small'' means zero.   We demonstrate that guessing
a threshold is inherently difficult, making the high accuracy obtained with
second-order methods desirable. 
  
We start from a standard Gauss-Seidel alternating block
framework and show that each block subproblem is further separable
into a set of independent functions, 
each of which depends on only a subset of variables.
We optimize each subset of variables independently, an obvious idea which has
nevertheless not previously appeared in the setting of sparse tensors.
We call this a \emph{row subproblem formulation} because the subset of
variables corresponds to one row of a factor matrix.
Each row subproblem amounts to minimizing a strictly convex function with
nonnegativity constraints, which we solve
using two-metric gradient projection techniques and exact or approximate
second-order information.

The importance of the row subproblem formulation is demonstrated in Section
\ref{subsec:exp-subprob}, where we show that applying a second-order method
directly to the block subproblem is highly inefficient.  We provide evidence
that a more effective way to apply second-order methods
is through the use of the row subproblem formulation.

Our contributions in this paper are as follows:
\begin{enumerate}
  \item A new formulation for nonnegative tensor factorization based on the
        Kullback-Leibler divergence objective that allows for the effective
        use of second-order optimization methods.  The optimization problem
        is separated into row subproblems containing $R$ variables,
        where $R$ is the number of factors in the model.
        The formulation makes row subproblems independent, suggesting
        a parallel method, although we do not explore parallelism in
        this paper.
  \item Two Matlab algorithms for computing factorizations of sparse
        nonnegative tensors:  one using second derivatives and the other using
        limited-memory quasi-Newton approximations.  The algorithms are
        made robust with an Armijo line search, damping modifications when the
        Hessian is ill conditioned, and projections to the bound of zero based
        on two-metric gradient projection ideas in
        \cite{bertsekas1982projected}.
        The two algorithms have different computational costs:
        the second derivative method is preferred when $R$ is small,
        and the quasi-Newton when $R$ is large.
  \item Test results that compare the performance of our two new
        algorithms with the best available multiplicative update method and
        a related quasi-Newton algorithm that does not
        formulate using row subproblems.
        The most significant test results are reported in this paper;
        detailed results of all
        experiments are available in the supplementary material (Appendix B).
  \item Test results showing the ability of our methods to quickly and
        accurately determine which elements of the factorization model
        are zero without using problem-specific thresholds.
\end{enumerate}

The paper is outlined as follows:  the remainder of Section 1 surveys related
work and provides a review of basic tensor properties.
Section 2 formalizes the Poisson nonnegative tensor factorization 
optimization problem, shows how the Gauss-Seidel alternating block framework can
be applied, and converts the block subproblem into independent row subproblems.
Section 3 outlines two algorithms for solving the row subproblem, one based
on the damped Hessian (PDN-R for projected damped Newton), and one based
on a limited-memory approximation (PQN-R for projected quasi-Newton).
Section 4 details numerical results on synthetic and real data sets and
quantifies the accuracy of finding a truly sparse factorization.
Additional test results are available in the supplementary material.
Section 5 contains a summary of the paper and concluding remarks.

%%%%%%%%%%%%%%%
\subsection{Related Work}  \label{subsec-relatedwork}
%%%%%%%%%%%%%%%

In this paper, we specifically consider nonnegative tensor factorization (NTF) 
in the case of the canonical polyadic (also known as CANDECOMP/PARAFAC) tensor decomposition. Our focus is on the K-L divergence objective function, but we also mention related work for the least squares (LS) case.
Additionally, we consider related work for nonnegative matrix factorization (NMF) for both K-L and LS. Note that there is much more work in the LS case, but the K-L objective function is different enough that it deserves its own attention.
We do not discuss other decompositions such as Tucker.

NMF in the LS case was first proposed by
Paatero and Tapper \cite{PaateroTapper1994} and also studied by Bro~\cite[p.~169]{BroPhD}.
Lee and Seung later consider the problem for both LS and K-L formulations and introduce multiplicative updates based on the convex subproblems
\cite{seungLeeNature1999, seung2001algorithms}.
Their work is extended to tensors by
Welling and Weber~\cite{WellingWeber2001}.
Many other works have been published on the LS versions of NMF
\cite{lin2007projected, KimSraDhillon2008, KimPark2008,
      PaateroTapper1994, ZhengZhang}
and NTF
\cite{BroDeJong1997, cichocki2009fast, FriedlanderHatz2008,
      kim2012fast, xu2012block}.

Lee and Seung's multiplicative update method \cite{seungLeeNature1999, seung2001algorithms,WellingWeber2001}
is the basis for most NTF algorithms that minimize the K-L divergence function.
Chi and Kolda provide an improved multiplicative update scheme for
K-L that addresses
performance and convergence issues as elements approach
zero \cite{ChiKolda:cpapr}; we compare to their method in
Section~\ref{sec:exp}.  By interpreting the  K-L divergence as an
alternative Csiszar-Tusnady procedure,
Zafeiriou and Petrou~\cite{zafeiriou2011nonnegative}
provide a probabilistic interpretation of NTF along with a new multiplicative
update scheme.
The multiplicative update is equivalent to a scaled steepest-descent step
\cite{seung2001algorithms}, so it is a first-order optimization method.
Since our method uses second-order information, it
allows for convergence to
higher accuracy and a better determination of sparsity in the factorization.

Second-order information has been used before in connection with the
K-L objective.
Zdunek and Cichocki~\cite{zdunek2006non, zdunek2007nonnegative}
propose a hybrid method for blind source separation applications via NMF that
uses a damped Hessian method similar to ours.
They recognize that the Hessian of the K-L objective
has a block diagonal structure but do not reformulate the optimization problem
further as we do. Consequently, their Hessian matrix is large, and
they switch to the LS objective function for the larger
mode of the matrix because their Newton method cannot scale up.
Mixing objective functions in this manner is undesirable because it combines
two different underlying models.
As a point of comparison, a problem in~\cite{zdunek2007nonnegative}
of size $200 \times 1000$ is considered too large for their Newton method,
but our algorithms can factor a data set of this size with $R=50$ components
to high accuracy in less than ten minutes
(see the supplementary material).
The Hessian-based method in~\cite{zdunek2007nonnegative} has most of the
advanced optimization features that we use (though details differ), including
an Armijo line search, active set identification, and an adjustable
Hessian damping factor.
We also note that
Zheng and Zhang~\cite{ZhengZhang} compute a damped Hessian search direction
and find an iterate with a backtracking line search, though this work is for
the LS objective in NMF.

Recently, Hsiel and Dhillon \cite{DhillonOneVar2011}
reported algorithms for NMF with both LS and K-L objectives. Their method updates one variable at a time, solving
a nonlinear scalar function using Newton's method with a constant step size.
They achieve good performance for the LS objective by taking the
variables in a particular order based on gradient information; however, for the
more complex K-L objective, they must cycle through all the variables one by one.
Our algorithms solve convex row subproblems with $R$ variables using second-order
information; solving these subproblems one variable at a time by coordinate
descent will likely have a much slower rate of
convergence~\cite[pp.~230-231]{nocedal1999numerical}.

A row subproblem reformulation similar to ours is noted in earlier papers
exploring the LS objective, but it never led to Hessian-based
methods that exploit sparsity as ours do.
Gonzales and Zhang use the reformulation with a multiplicative update method
for NMF~\cite{gonzalez2005accelerating} but do not generalize to tensors or
the K-L objective.
Phan et al.\@~\cite{PhanZdunek2010} note the reformulation is suitable
for parallelizing a Hessian-based method for NTF using LS.
Kim and Park use the reformulation for NTF with LS~\cite{kim2012fast}, deriving
small bound-constrained LS subproblems.  Their method solves the
LS subproblems by exact matrix factorization, without exploiting
sparsity, and features a block principal pivoting method for choosing
the active set.
Other works solve the LS objective by taking advantage of
row-by-row or column-by-column subproblem decomposition
\cite{cichocki2009fast, liu2012sparse, PhanCichocki2011}.

Our algorithms are similar in spirit to the work of Kim,
Sra and Dhillon~\cite{KimSuvritDhillon:BoxConstraints},
which applies a projected quasi-Newton
algorithm (called PQN in this paper) to solving NMF with a K-L objective.
Like PQN, our algorithms
identify active variables, compute a Newton-like direction in the space of free
variables, and find a new iterate using a projected backtracking line search.
We differ from PQN in reformulating the subproblem and in computing
a damped Newton direction; both improvements make a huge difference
in performance for large-scale tensor problems.
We compare to PQN in Section~\ref{sec:exp}.

All-at-once optimization methods, including Hessian-based algorithms, have been
applied to NTF with the LS objective function.
As an example, Paatero replaces the nonnegativity constraints with
a barrier function~\cite{Paatero1997} to yield an unconstrained optimization
problem,
and Phan, Tichavsky and Cichocki~\cite{PhanCichocki2011} apply
a fast damped Gauss-Newton algorithm for minimizing a similar penalized
objective.
We are not aware of any work on all-at-once methods for
the K-L objective in NTF.

Finally, we note that all methods, including ours, find only a locally optimal
solution to the NTF problem.  Finding the global solution is generally much
harder; for instance, Vavasis~\cite{vavasis2009complexity} proves it is
NP-hard for an NMF model that fits the data exactly.

%%%%%%%%%%%%%%%
\subsection{Tensor Review}  \label{subsec:tensorreview}
%%%%%%%%%%%%%%%

For a thorough introduction to tensors, see~\cite{KoldaBader:SiamReview} and
references therein; we only review concepts that are necessary for
understanding this paper.
A tensor is a multidimensional array.  An $N$-way tensor $\TX$ has size
$I_1 \times I_2 \times \ldots \times I_N$.
To differentiate between tensors, matrices, vectors, and scalars, we use
the following notational convention:
$\TX$ is a tensor (bold, capitalized, calligraphic),
$\M{X}$ is a matrix (bold, capitalized),
$\V{x}$ is a vector (bold, lowercase),
and $x$ is a scalar (lowercase).
Additionally, given a matrix $\M{X}$, $\V{x}_j$ denotes its $j$th
column and $\rowx_i$ denotes its $i$th row.

\begin{figure}[bht!]
 \centering
 \includegraphics[scale=.20]{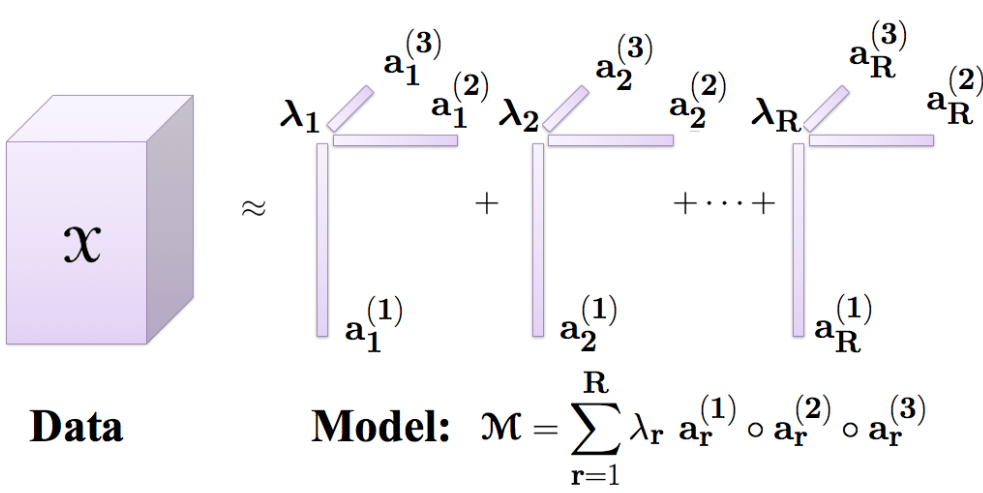}
 \caption{CANDECOMP/PARAFAC decomposition of a three-way tensor
          into $R$ components.}
 \label{fig:CP}
\end{figure}

Just as a matrix can be decomposed into a sum of
outer products between two vectors, an $N$-way tensor can be decomposed
into a sum of outer products between $N$ vectors. Each of these outer products
(called components) yields an $N$-way tensor of rank one.
The CP (CANDECOMP/PARAFAC) decomposition~\cite{CarrollChang:1970, Harshman:1970}
represents a tensor as a sum of rank-one tensors
(see Figure~\ref{fig:CP}):
\begin{equation}  \label{kruskal}
  \TX  \approx \KT{\Vl; \Mn{A}{1},\dots,\Mn{A}{N}}
       = \sum_{r=1}^R \lambda_r \, \V{a}^{(1)}_r \circ \ldots \circ \V{a}^{(N)}_r
\end{equation}
where $\Vl$ is a vector and each $\Mn{A}{n}$ is an $I_n \times R$
\emph{factor matrix} containing the $R$ vectors contributed to the outer
products by mode $n$, i.e.,
\begin{equation}  \label{factors}
  \Mn{A}{n} = [ \V{a}^{(n)}_1 \cdots \V{a}^{(n)}_R ] .
\end{equation}
Equality holds in~(\ref{kruskal}) when $R$ equals the rank of $\TX$,
but often a tensor is approximated by a smaller number of terms.
We let $\V{i}$ denote the multi-index $(i_1, i_2, \ldots, i_N)$ of an element
$\TE{x}{i}$ of $\TX$.

%%%%%%%%%%%%%%%%%%%%
%%%%%%%%%%%%%%%%%%%%%%

We use of the idea of matricization, or unfolding a tensor
into a matrix.  Specifically,
unfolding along mode $n$ yields a matrix of size $I_n \cdot J_n$, where
\begin{displaymath}
  J_n = I_1 \cdot I_2 \; \cdots \; I_{n-1} \cdot I_{n+1}
            \; \cdots \; I_{N-1} \cdot I_N .
\end{displaymath}
We use the notation $\Mz{X}{n}$ to represent a tensor $\TX$ that has been
unfolded so that its $n$th mode forms the rows of the matrix,
and $x^{(n)}_{ij}$ for its $(i,j)$ element.
If a tensor $\TX$ is written in
Kruskal form (\ref{kruskal}), then the mode-$n$ matricization
is given by
\begin{equation}  \label{matricization}
  \Mz{X}{n} = \Mn{A}{n} \M{\Lambda} ( \Mn{A}{N} \odot \ldots \odot \Mn{A}{n+1}
                                \odot \Mn{A}{n-1} \odot \ldots \odot\Mn{A}{1})^T
\end{equation}
where $\M{\Lambda} = \text{diag} (\Vl)$
and $\odot$ denotes the Khatri-Rao product~\cite{KoldaBader:SiamReview}.

Tensor results are generally easier to interpret when the factors (\ref{factors})
are sparse.
Moreover, many sparse count applications can reasonably expect sparsity in the
factors.  For example, the 3-way data considered in~\cite{Kolda:DBLP}
counts publications by authors at various conferences over a ten year
period. The tensor representation has a
sparsity of 0.14\% (only 0.14\% of the data elements are nonzero), and the
factors computed by our algorithm with $R=60$
(see Section~\ref{subsec:exp-fullprob}) have sparsity 9.3\%, 2.7\%,
and 77.5\% over the three modes.
One meaning of sparsity in the factors is to say
that a typical outer product term connects about 9\% of the authors with
3\% of the conferences in 8 of the 10 years.  Linking particular authors and
conferences
is an important outcome of the tensor analysis, requiring clear distinction
between zero and nonzero elements in the factors.

\section{Poisson Nonnegative Tensor Factorization}

In this section we state the optimization problem, examine its structure,
and show how to separate it into simpler subproblems.

\subsection{Gauss-Seidel Alternating Block Formulation}

We seek a tensor model in CP form to approximate data $\TX$:
\begin{equation*}
  \TX \approx \TM
      = \KT {\Vl; \Mn{A}{1},\dots,\Mn{A}{N}}
      = \sum_{r=1}^R \lambda_r \, \V{a}^{(1)}_r \circ \ldots \circ \V{a}^{(N)}_r .
\end{equation*}
The value of $R$ is chosen empirically,
and the scaling vector $\Vl$ and factor matrices
$\Mn{A}{n}$ are the model parameters that we compute.

In \cite{ChiKolda:cpapr}, it is shown that a K-L divergence objective function
results when data elements follow Poisson distributions
with multilinear parameters.
The best-fitting tensor model under this assumption satisfies:
\begin{align}
  \min_{\Vl, \Mn{A}{1}, \ldots, \Mn{A}{N}}  \;
  & f(\TM) = \sum_{\V{i}} \TE{m}{i} - \TE{x}{i} \log \TE{m}{i}
    \nonumber
  \\
  \text{s.t.}  \;\;\;\;
  & \displaystyle
    \TM = \sum_{r=1}^R \lambda_r \, \V{a}^{(1)}_r \circ \ldots \circ \V{a}^{(N)}_r,
    \label{fullprob}
  \\
  & \lambda_r \geq 0, \V{a}^{(n)}_r \geq 0, \| \V{a}^{(n)}_r \|_1 = 1,
      \;\; \forall r \in \{1, \ldots R \},
      \;   \forall n \in \{1, \ldots N \} ,
    \nonumber
\end{align}
where $\TE{x}{i}$ denotes element $(i_1, \ldots, i_n)$ of tensor $\T{X}$
and $\TE{m}{i}$ denotes element $(i_1, \ldots, i_n)$ of the model $\T{M}$.
The model may have terms where $\TE{m}{i} = 0$ and $\TE{x}{i} = 0$; for this
case we define $0 \log 0 = 0$.
Note that for matrix factorization,
(\ref{fullprob}) reduces
to the K-L divergence used by
Lee and Seung \cite{seungLeeNature1999, seung2001algorithms}.
The constraint that normalizes the column sum of the factor matrices serves
to remove an inherent scaling ambiguity in the CP factor model.

As in \cite{ChiKolda:cpapr}, we unfold $\TX$ and $\TM$ into their $n$th
matricized mode, and use (\ref{matricization}) to express the objective as
\begin{equation*}
  f(\TM) = \V{e}^T [ \Mn{A}{n} \M{\Lambda} \Mn{\Pi}{n}
                     - \M{X}_{(n)} \ast \log (\Mn{A}{n} \M{\Lambda} \Mn{\Pi}{n}) ]
                   \V{e} \; ,
\end{equation*}
where $\V{e}$ is a vector of all ones,
the operator $\ast$ denotes elementwise multiplication,
$\log(\cdot)$ is taken elementwise,
\begin{eqnarray}
  \label{pimatrix}  
  \M{\Pi}^{(n)} & = & ( \Mn{A}{N} \odot \ldots \odot \Mn{A}{n+1} \odot \Mn{A}{n-1}
                      \odot \ldots \odot \Mn{A}{1} )^T
                    \in {\bf R}^{R \times J_n}  \; , \; \mbox{and}
      \\ 
      \nonumber
  \M{\Lambda}   & = & \text{diag} (\Vl) \in {\bf R}^{R \times R} .  
\end{eqnarray}
Note that by expanding the Khatri-Rao products in~(\ref{pimatrix}) and remembering
that column vectors $\V{a}^{(n)}_r$ are normalized,
each row of $\M{\Pi}^{(n)}$ conveniently sums to one.
This is a consequence of using the $\ell_1$ norm in (\ref{fullprob}).

The above representation of the objective motivates the use of an
alternating block optimization method where only one factor matrix is optimized
at a time.
Holding the other factor matrices fixed, the optimization problem for
$\Mn{A}{n}$ and $\M{\Lambda}$ is
\begin{equation}
  \begin{split}
    \min_{\M{\Lambda}, \Mn{A}{n}} \;
      & f(\M{\Lambda}, \; \Mn{A}{n}) =
             \V{e}^T [ \Mn{A}{n} \M{\Lambda} \Mn{\Pi}{n}
                       - \M{X}_{(n)} \ast
                                     \log (\Mn{A}{n} \M{\Lambda} \Mn{\Pi}{n})
                     ] \V{e}  \\
    \text{s.t.} \quad
      & \M{\Lambda} \geq 0, \; \Mn{A}{n} \geq 0, \; \V{e}^T \Mn{A}{n} = \V{1} .
  \end{split}
  \label{modeNprobNonconvex}
\end{equation}

Problem (\ref{modeNprobNonconvex}) is not convex.
However, ignoring the equality constraint and letting
$\Mn{B}{n} = \Mn{A}{n} \M{\Lambda}$, we have
\begin{equation}
  \min_{\Mn{B}{n} \geq 0} \;\;
      f(\Mn{B}{n}) = \V{e}^T [ \Mn{B}{n} \Mn{\Pi}{n}
                               - \Mz{X}{n} \ast \log( \Mn{B}{n} \Mn{\Pi}{n}) ]
                     \V{e}  \label{modeNprob}
\end{equation}
which is convex with respect to $\Mn{B}{n}$.
The two formulations are equivalent in that a KKT point of (\ref{modeNprob})
can be used to find a KKT point of (\ref{modeNprobNonconvex}).
Chi and Kolda show in \cite{ChiKolda:cpapr} that (\ref{modeNprob})
is \emph{strictly} convex
given certain assumptions on the sparsity pattern of $\Mz{X}{n}$.

We pause to think about (\ref{modeNprob}) when the tensor is two-way.
In this case, we solve for two factor matrices by alternating over two block
subproblems; for instance, with $n=1$ the subproblem~(\ref{modeNprob})
finds $\Mn{B}{1}$ with $\Mn{\Pi}{1} = (\Mn{A}{2})^T$.
For an $N$-way problem, the only change to (\ref{modeNprob}) is $\Mn{\Pi}{n}$,
which grows in size exponentially with each additional factor matrix.
To efficiently solve the subproblems (\ref{modeNprob}) for large sparse
tensors we can exploit sparsity to reduce the computational costs.
As discussed in the next section, columns of $\Mn{\Pi}{n}$ need to be computed
only when the corresponding column in the unfolded tensor $\Mz{X}{n}$
has a nonzero element.
Our row subproblem formulation carries this idea further because each
row subproblem
generally uses only a fraction of the nonzero elements in $\Mz{X}{n}$.

At this point we define Algorithm~\ref{alg:outline}, a Gauss-Seidel
alternating block method.  The algorithm iterates over each mode of the tensor,
solving the convex optimization block subproblem.  Steps~\ref{algoutline:rescale1}
and \ref{algoutline:rescale2} rescale the factor matrix columns, redistributing
the weight into $\Vl$.  For the moment, we leave the subproblem solution method
in Step~\ref{algoutline:subproblem} unspecified.
A proof that Algorithm~\ref{alg:outline} convergences to a local minimum
of (\ref{fullprob}) is given in \cite{ChiKolda:cpapr}.

\begin{algorithm}
  \caption{Alternating Block Framework}
  \label{alg:outline}
  Given data tensor $\TX$ of size $I_1 \times I_2 \times \cdots \times I_N$,
  and the number of components $R$
  \\
  Return a model $\TM = [\V{\lambda}; \Mn{A}{1} \ldots \Mn{A}{N}]$

  \begin{algorithmic}[1]
    \State
    Initialize $\Mn{A}{n} \in {\bf R}^{I_n \times R}$ for $n = 1, \ldots, N$
    \Repeat
      \For{$n=1, \dots, N$} 
        \label{algoutline:for}
        \State
	  Let $\M{\Pi}^{(n)}= ( \Mn{A}{N} \odot \ldots \odot \Mn{A}{n+1}
                               \odot \Mn{A}{n-1} \odot \ldots \odot\Mn{A}{1})^T$
          \label{algoutline:formPi}
        \State
          Use Algorithm~\ref{alg:subproblem-sparse} to
          compute $\M{B}^*$ that minimizes $f(\Mn{B}{n})$
                                      s.t. $\Mn{B}{n} \geq 0$
          \label{algoutline:subproblem}
        \State
          $\V{\lambda} \leftarrow \eb^T\M{B}^*$
          \label{algoutline:rescale1}
        \State
          $\Mn{A}{n} \leftarrow \M{B}^* \M{\Lambda}^{-1}, \;$
          where $\M{\Lambda} = \text{diag}(\Vl)$
          \label{algoutline:rescale2}
      \EndFor
          \label{algoutline:endfor}
    \Until all mode subproblems have converged
  \end{algorithmic}
\end{algorithm}

This outline of Algorithm~\ref{alg:outline} corresponds exactly with the
method proposed in \cite{ChiKolda:cpapr}; where we differ is in how to
solve subproblem~(\ref{modeNprob}) in Step~\ref{algoutline:subproblem}.
Note also that this algorithm is the same as for the least squares objective (references were given in
Section~\ref{subsec-relatedwork}); there the subproblem
in Step~\ref{algoutline:subproblem} is replaced by a linear least
squares subproblem.
We now proceed to describe our method for solving~(\ref{modeNprob}).

%%%%%%%%%%%%%%%%%%%%%%%%%%%%%%%%
\subsection{Row Subproblem Reformulation}  \label{subsec:ptf-reformulation}
%%%%%%%%%%%%%%%%%%%%%%%%%%%%%%%%

We examine the objective function $f(\Mn{B}{n})$ in (\ref{modeNprob})
and show that it can be reformulated into independent functions.
As mentioned in the previous section, rows of $\M{\Pi}^{(n)}$ sum to one if
the columns of factor matrices are nonnegative and sum to one.
When $\M{\Pi}^{(n)}$ is formed at Step~\ref{algoutline:formPi}
of Algorithm~\ref{alg:outline}, the factor matrices satisfy these conditions
by virtue of Steps~\ref{algoutline:rescale1} and \ref{algoutline:rescale2};
hence, the first term of $f(\Mn{B}{n})$ is
\begin{displaymath}
  \V{e}^T \Mn{B}{n} \Mn{\Pi}{n} \V{e}
    \; = \; \V{e}^T \Mn{B}{n} \V{e}
    \; = \; \sum_{i=1}^{I_n} \sum_{r=1}^R b^{(n)}_{ir} .
\end{displaymath}

The second term of $f(\Mn{B}{n})$ is a sum of elements from the $I_n \times J_n$
matrix $\Mz{X}{n} \ast \log( \Mn{B}{n} \Mn{\Pi}{n})$.  Recall that the
operations in this expression are elementwise, so the scalar matrix
element $(i,j)$ of the term can be written as
\begin{equation*}
  x^{(n)}_{ij} \log \left( \sum_{r=1}^R b^{(n)}_{ir} \pi^{(n)}_{rj} \right) .
\end{equation*}
Adding all the elements and combining with the first term gives
\begin{eqnarray*}
  f(\Mn{B}{n}) & = & \sum_{i=1}^{I_n} \sum_{r=1}^R b^{(n)}_{ir}
                     - \sum_{i=1}^{I_n} \sum_{j=1}^{J_n}
      x^{(n)}_{ij} \log \left( \sum_{r=1}^R b^{(n)}_{ir} \pi^{(n)}_{rj} \right)
      \\
               & = & \sum_{i=1}^{I_n} f_{\rm row}
                                        ( \Vhat{b}_i^{(n)},
                                          \Vhat{x}_i^{(n)},
                                          \Mn{\Pi}{n} ) .
\end{eqnarray*}
where $\Vhat{b}_i$ and $\Vhat{x}_i$ are the $i$th row vectors of their
corresponding matrices, and
\begin{equation}
  f_{\rm row} (\Vhat{b}, \Vhat{x}, \M{\Pi})
  = \sum_{r=1}^R \hat{b}_r
    - \sum_{j=1}^{J_n} \hat{x}_j
                       \log \left( \sum_{r=1}^R \hat{b}_r \pi_{rj} \right) .
  \label{eq:frow}
\end{equation}
Problem~(\ref{modeNprob}) can now be rewritten as
\begin{equation}
  \min_{\Vhat{b}_1, \ldots, \Vhat{b}_{I_n} \; \geq 0} \;\;
      \sum_{i=1}^{I_n} f_{\rm row} ( \Vhat{b}_i^{(n)},
                                     \Vhat{x}_i^{(n)},
                                     \Mn{\Pi}{n} ) .
  \label{rsp-ptf}
\end{equation}
This is a completely separable set of $I_n$ \emph{row subproblems}, each one
a convex nonlinear optimization problem containing $R$ variables.  The relatively
small number of variables makes second-order optimization methods tractable,
and that is the direction we pursue in this paper.
Algorithm~\ref{alg:subproblem-sparse} describes how the reformulation fits into
Algorithm~\ref{alg:outline}.

\begin{algorithm}
  \caption{Row Subproblem Framework for Solving (\ref{rsp-ptf})}
  \label{alg:subproblem-sparse}
  Given $\Mz{X}{n}$ of size $I_n \times J_n$,
  and $\Mn{\Pi}{n}$ of size $R \times J_n$
  \\
  Return a solution $\M{B}^*$  consisting of row vectors
  $\Vhat{b}_1^*, \ldots, \Vhat{b}_{I_N}^*$

  \begin{algorithmic}[1]
    \For{$i=1, \dots, I_n$}
      \State Select row $\Vhat{x}_i$ of $\Mz{X}{n}$
      \State Generate one column of $\Mn{\Pi}{n}$ for each nonzero in $\Vhat{x}_i$
          \label{algsubproblem:sparsePi}
      \State Use Algorithm~\ref{alg:rowsubprobPDN} or \ref{alg:rowsubprobPQN}
             to compute $\Vhat{b}_i^*$ that solves
             \vspace*{-0.10in}
             \begin{equation*}
               \min \; f_{\rm row} (\Vhat{b}_i^{(n)}, \Vhat{x}_i^{(n)}, \Mn{\Pi}{n})
               \;\; \mbox{subject to} \;\; \Vhat{b}_i \geq 0
             \end{equation*}
             \vspace*{-0.20in}
          \label{algsubproblem:rowsubproblem}
    \EndFor
  \end{algorithmic}
\end{algorithm}

The independence of row subproblems is crucial for handling large tensors.
For example, if a three-way tensor of size $1000 \times 1000 \times 1000$
is factored
into $R=100$ components, then $\Mn{\Pi}{n}$ is a $100 \times 10^6$ matrix.
However, elements of $\Mn{\Pi}{n}$ appear in the optimization objective only
where the matricized tensor $\Mz{X}{n}$ has nonzero elements, so in a sparse
tensor many columns of $\Mn{\Pi}{n}$ can be ignored;
this point was first published in \cite{ChiKolda:cpapr}.
Algorithm~\ref{alg:subproblem-sparse} exploits this fact in
Step~\ref{algsubproblem:sparsePi}.

Algorithm~\ref{alg:subproblem-sparse} also points the way to a parallel
implementation of the CP tensor factorization.
We note, as did \cite{PhanZdunek2010}, that each row subproblem can be run in
parallel and storage costs are determined by the sparsity of the data.
In a distributed computing architecture, an algorithm could
identify the nonzero elements of each row subproblem at the beginning of
execution and collect only the data needed to form appropriate columns
of $\Mn{\Pi}{n}$ at a given processing element.
We do not implement a parallel version of the algorithm in this paper.

\section{Solving the Row Subproblem}

In this section we show how to solve the row subproblem~(\ref{rsp-ptf})
using second-order information.  We describe two algorithms, one applying
second derivatives in the form of a damped Hessian matrix, and the other
using a quasi-Newton approximation of the Hessian.  Both algorithms use
projection, but the details differ.

Each row subproblem consists of minimizing a strictly convex function
of $R$ variables with nonnegativity constraints.
One of the most effective methods for solving bound-constrained problems is
second-order gradient projection; see~\cite{schmidt2012}.
We employ a form of two-metric gradient projection from
Bertsekas~\cite{bertsekas1982projected}.
Each variable is marked in one of three states based on its gradient and
current location:  fixed at its bound of zero,
allowed to move in the direction of steepest-descent, or free to move
along a Newton or quasi-Newton search direction.
Details are in Section~\ref{subsec:subprobActiveSet}.

An alternative to Bertsekas is to use methods that employ gradient projection
searches to determine the active variables (those set to zero). Examples
include the generalized Cauchy point~\cite{conn1988global} and gradient projection
along the steepest-descent direction with a line search.
We experimented with using the generalized Cauchy point to 
determine the active variables, but preliminary results indicated
that this approach sets too many variables to be at their bound,
leading to more iterations and poor overall performance.
Gradient projection steps with a line search calls for an extra
function evaluation, which is computationally expensive.
Given a more efficient method for evaluating the function, this might be a
better approach since, under mild conditions, gradient projection methods
find the active set in a finite number of
iterations~\cite{bertsekas1976goldstein}.

For notational convenience, in this section we use $\V{b}$ for the
column vector representation of row vector $\Vhat{b}_i$;
that is, $\V{b} = \Vhat{b}_i^T$.  Iterations are denoted with superscript $k$,
and $\nabla_r$ represents the derivative with respect to the $r$th
variable.
Let $P_{+}[\V{v}]$ be the projection operator that restricts each element of
vector $\V{v}$ to be nonnegative.
We make use of the first and second derivatives of
$f_{\rm row}$, given by
\begin{eqnarray}
  \nabla_r f_{\rm row}(\V{b}) & = &
  \frac{ \partial f_{\rm row} (\V{b}, \V{x}, \M{\Pi}) }
       { \partial {b}_r }
    \; = \; 1 - \sum_{j=1}^{J_n} \frac{ x_j \pi_{rj} }
                                      { \sum_{i=1}^R {b}_i \pi_{ij} } ,
  \label{eq:grad_frow}
  \\
  \nabla^2_{rs} f_{\rm row}(\V{b}) & = &
  \frac{ \partial^2 f_{\rm row} (\V{b}, \V{x}, \M{\Pi}) }
       { \partial {b}_r \; \partial {b}_s }
    \; = \; \sum_{j=1}^{J_n} \frac{ x_j \pi_{rj} \pi_{sj} }
                                  { (\sum_{i=1}^R {b}_i \pi_{ij})^2 } .
  \nonumber
\end{eqnarray}

%%%%%%%%%%%%%%%%%%%%%%%%%%%%%%%%
\subsection{Two-Metric Projection}  \label{subsec:subprobActiveSet}
%%%%%%%%%%%%%%%%%%%%%%%%%%%%%%%%

At each iteration $k$ we must choose a set of variables to update such that
progress is made in decreasing the objective.
Bertsekas demonstrated in~\cite{bertsekas1982projected} that iterative updates
of the form
\begin{equation*}
  \V{b}^{k+1} = P_{+}[ \V{b}^k - \alpha \M{M}^k \; \nabla f_{\rm row}(\V{b}^k) ]
\end{equation*}
are not guaranteed to decrease the objective function unless $\M{M}^k$
is a positive diagonal matrix.  Instead,
it is necessary to predict the variables that have the potential to make
progress in decreasing the objective and then update just those variables
using a positive definite matrix.
We present the two-metric technique of Bertsekas as it is executed in our
algorithm, which differs superficially from the presentation
in~\cite{bertsekas1982projected}.

A variable's potential effect on the objective is determined by how close it is
to zero and by its direction of steepest-descent.
If a variable is close to zero and its steepest-descent direction points
towards the negative orthant, then the next update will likely project
the variable to zero and its small displacement will have little effect on
the objective.
A closeness threshold $\epsilon_k$ is computed from a user-defined parameter
$\epsilon > 0$ as
\begin{equation}
  \epsilon_k = \min(w_k, \epsilon) , \;\;
    w_k = \left\| \V{b}^k-P_+[\V{b}^k-  \nabla f_{\rm row}(\V{b}^k) ] \right\|_2 .
  \label{ActiveSetEpsilon}
\end{equation}
We then define index sets
\begin{equation}
  \begin{aligned}
    \mathcal{A}(\V{b}^k) &=  \left\{ r \ | \ {b}^k_r  = 0, \;
                  \nabla_{r} f_{\rm row}(\V{b}^k) > 0 \right\} ,
    \\
    \mathcal{G}(\V{b}^k) &=  \left\{ r \ | \ 0 < {b}^k_r \leq \epsilon_k, \;
                  \nabla_{r} f_{\rm row}(\V{b}^k) > 0 \right\} ,
    \\
    \mathcal{F}(\V{b}^k)   &= \left( \mathcal{A}(\V{b}^k)  \ \cup \
                                     \mathcal{G}(\V{b}^k) \right)^c ,
  \end{aligned}
  \label{ActiveSetDef}
\end{equation}
where superscript $c$ denotes the set complement.
Variables in the set $\mathcal{A}$ are fixed at zero, variables in
$\mathcal{G}$ move in the direction of the negative gradient, and variables
in $\mathcal{F}$ are free to move according to second-order information.
Note that if $\epsilon = 0$ then $\epsilon_k = 0$, $\mathcal{G}$ is empty,
and the method reduces to defining an active set of variables by instantaneous
line search~\cite{andrew2007scalable}.

%%%%%%%%%%%%%%%%%%%%%%%%%%%%%%%%
\subsubsection{Damped Newton Step}  \label{subsec:DampedNewton}
%%%%%%%%%%%%%%%%%%%%%%%%%%%%%%%%

The damped Newton direction is taken with respect to only the variables
in the set $\mathcal{F}$ from (\ref{ActiveSetDef}).  Let
\begin{equation*}
  \V{g}^k_F = [\nabla f_{\rm row}(\V{b}^k)]_{\mathcal{F}}, \
  \M{H}^k_F = [\nabla^2 f_{\rm row}(\V{b}^k)]_{\mathcal{F}}, \
  \V{b}^k_F = [\V{b}^k]_{\mathcal{F}} ,
\end{equation*}
where $[\V{v}]_{\mathcal{F}}$ chooses the elements of vector $\V{v}$ corresponding
to variables in the set $\mathcal{F}$.
Since the row subproblems are strictly convex, the full Hessian and $\M{H}^k_F$
are positive definite.

The damped Hessian has its roots in trust region methods.
At every iteration we form a quadratic approximation $m_k$ of the objective
plus a quadratic penalty.
The penalty serves to ensure that the next iterate
does not move too far away from the current iterate, which is important when the
Hessian is ill conditioned.  The quadratic model plus penalty expanded
about $\V{b}^k$ for variables $\V{d}_F \in \mathbb{R}^{|\mathcal{F}|}$ is
\begin{equation}
  m_k(\V{d}_F; \mu_k) = f_{\rm row}(\V{b}^k) +  \V{d}^T_F \V{g}^k_F
                        + \frac{1}{2} \V{d}^T_F \M{H}^k_F \V{d}_F
                        + \frac{\mu_k}{2} \left\| \V{d}_F \right\|_2^2 .
    \label{quadModel}
\end{equation}
The unique minimum of $m_k(\cdot)$ is
\begin{equation*}
  \V{d}^k_F = -(\M{H}^k_F + \mu_k \M{I})^{-1} \V{g}^k_F ,
\end{equation*}
where $\M{H}^k_F + \mu_k \M{I}$ is known as the damped Hessian.
Adding a multiple of the identity to $\M{H}^k_F$ increases each of
its eigenvalues by $\mu_k$, which has the effect of diminishing the length of
$\V{D}^k_F$, similar to the action of a trust region.
The step $\V{D}^k_F$ is computed using a Cholesky factorization of the
damped Hessian, and the full space search direction $\V{d}^k \in \mathbb{R}^R$
is then given by
\begin{equation} \label{finaldirection}
  d^k_r =
    \begin{cases}
        (E \, \V{d}^k_F)_r                 & r \in \mathcal{F}  \\
        -\nabla_r f_{\rm row}(\V{b}^k)     & r \in \mathcal{G}  \\
        0                                  & r \in \mathcal{A} ,
    \end{cases}
\end{equation}
where index sets $\mathcal{A}$, $\mathcal{G}$, and $\mathcal{F}$
are defined in (\ref{ActiveSetDef}),
and $E$ is an $R \times |\mathcal{F}|$ matrix that elongates the
vector $\V{d}^k_F$ to length $R$.
Specifically, $e_{ij} = 1$ if $i$ is the row subproblem variable corresponding
to the $j$th variable in $\mathcal{F}$,
and zero otherwise.

The damping parameter $\mu_k$ is adjusted by a Levenberg-Marquardt
strategy~\cite{nocedal1999numerical}.
First define the ratio of actual reduction over predicted reduction,
\begin{equation}
\rho = \frac{f_{\rm row}(\V{b}^k +\V{d}^k ) - f_{\rm row}(\V{b}^k )}
            {m_k(\V{d}^k_F; 0) - m_k(0; 0)} ,  \label{rho}
\end{equation}
where $m_k(\cdot)$ is defined by (\ref{quadModel}).
Note the numerator of (\ref{rho}) calculates $f_{\rm row}$ using all variables,
while the denominator calculates $m_k(\cdot)$ using only the variable
in $\mathcal{F}$.
The damping parameter is updated by the following rule
\begin{equation} \label{updatedamping}
  \mu_{k+1} =
      \begin{cases} \frac{7}{2}\mu_k & \text{if } \rho < \frac{1}{4} , \\
                    \frac{2}{7}\mu_k & \text{if } \rho > \frac{3}{4} , \\
                    \mu_k & \text{otherwise.}
      \end{cases}
\end{equation}
Since $\V{d}^k_F$ is the minimum of (\ref{quadModel}), the denominator of
(\ref{rho}) is always negative. If the search direction $\V{d}^k$ increases
the objective function, then the numerator of (\ref{rho}) will be positive;
hence $\rho < 0$ and the damping parameter will be increased for the next
iteration. On the other hand, if the search direction $\V{d}^k$ decreases the
objective function, then the numerator will be negative; hence $\rho >0$ and
the relative sizes of the actual reduction and predicted reduction will
determine how the damping parameter is adjusted.

%%%%%%%%%%%%%%%%%%%%%%%%%%%%%%%%
\subsubsection{Line Search}
%%%%%%%%%%%%%%%%%%%%%%%%%%%%%%%%

After computing the search direction $\V{d}^k$, we ensure the next iterate
decreases the objective by using a projected backtracking line search that
satisfies the Armijo condition~\cite{nocedal1999numerical}.
Given scalars $0 < \beta$ and $\sigma < 1$,
we find the smallest nonnegative integer $t$ that satisfies the inequality
\begin{equation}
  f_{\rm row}(P_{+}[\V{b}^k + \beta^t \V{d}^k] )- f_{\rm row}(\V{b}^k)
  \leq
  \sigma (P_{+}[\V{b}^k + \beta^t \V{d}^k]-\V{b}^k)^T
         \nabla f_{\rm row}(\V{b}^k ) .   \label{armijo}
\end{equation}
We set $\alpha_k = \beta^t$ and the next iterate is given by
\begin{equation*}
  \V{b}^{k+1} = P_+[ \V{b}^k + \alpha_k\V{d}^{k} ] .
\end{equation*}

%%%%%%%%%%%%%%%%%%%%%%%%%%%%%%%%
\subsection{Projected Quasi-Newton Step}  \label{subsec:PQNmods}
%%%%%%%%%%%%%%%%%%%%%%%%%%%%%%%%

As an alternative to the damped Hessian step, we adapt the projected
quasi-Newton step from~\cite{KimSuvritDhillon:BoxConstraints}.  Their work
employs a limited-memory BFGS (L-BFGS) approximation~\cite{nocedalLBFGS}
in a framework suitable for any convex, bound-constrained problem.

L-BFGS estimates Hessian properties based on the most
recent $M$ update pairs $\{ \V{s}^i, \V{y}^i \}$,
$i \in [ \, \max \{ 1, k-M \}, k \, ]$, where
\begin{equation} \label{curvPairs}
  \V{s}^i = \V{b}^{i+1} - \V{b}^{i}, \qquad
  \V{y}^i = \nabla f_{\rm row}(\V{b}^{i+1}) - \nabla f_{\rm row}(\V{b}^i) .
\end{equation}
L-BFGS uses a two-loop recursion through the stored pairs to efficiently
compute a vector
 $\V{p}^k = \Mtilde{B}^k \V{g}^k$, where $\Mtilde{B}^k$ approximates
the inverse of the Hessian $[\M{H}^k]^{-1}$
using the pairs $\{ \V{s}^i, \V{y}^i \}$.
Storage is set to $M=3$ pairs in all experiments.
On the first iterate when $k=0$, we use a multiple of the identity matrix
so that $\V{p}^0$ is in the direction of the gradient.
L-BFGS updates require the quantity $1 / ((\V{s}^i)^T \V{y}^i)$ to be
positive.  We check this condition and skip the update pair if it is violated.
This can happen if all row variables are at their bound of zero, or from
numerical roundoff at a point near a minimizer.
See~\cite[Chapter 7]{nocedal1999numerical} for further detail.

The projected quasi-Newton search direction $\V{d}^k$, analogous to
(\ref{finaldirection}), is
\begin{equation}
  d^k_r = \begin{cases}
             -(E \, p^k)_r                   & r \in \mathcal{F} ,  \\
             -\nabla_r f_{\rm row}(\V{b}^k)  & r \in \mathcal{G} ,  \\
             0                               & r \in \mathcal{A} ,
          \end{cases}
      \label{lbfgs_final}
\end{equation}
where $\mathcal{F}$ and $\mathcal{A}$ are determined from (\ref{ActiveSetDef}),
and $E$ is the elongation matrix defined in (\ref{finaldirection}).

The step $p^k$ is computed from an L-BFGS approximation over all
variables in the row subproblem; in contrast, the step $\V{d}^k_F$
computed from the damped Hessian in Section~\ref{subsec:DampedNewton}
is derived from the second derivatives of only the free variables in
$\mathcal{F}$.
We could build an L-BFGS model over just the free variables as is done
in~\cite{byrd1995limited}, but the computational cost is higher.
Our L-BFGS step is therefore influenced by second-order
information from variables not in $\mathcal{F}$.
This information is irrelevant to the step, but we find that algorithm
performance is still good.
We now express the influence in terms of the reduced Hessian and
inverse of the reduced Hessian.

Let $\M{H}$ and $\M{B}$ denote the true Hessian and inverse Hessian matrices
over all variables in a row subproblem.
Suppose the variables in $\mathcal{F}$ are the first $|\mathcal{F}|$
variables, and the remaining variables are in
${\mathcal{N}} = \mathcal{A} \ \cup \ \mathcal{G}$.
Then we can write $\M{H}$ in block form as
\begin{equation*}
  \M{H} = \left[ \begin{array}{cc}
                   \M{H}_{FF}  &  \M{H}_{NF}^T  \\
                   \M{H}_{NF}  &  \M{H}_{NN}
                 \end{array}
          \right] ,
\end{equation*}
with $\M{H}_{FF} \in \mathbb{R}^{|\mathcal{F}|\times |\mathcal{F}|}$,
     $\M{H}_{NF} \in \mathbb{R}^{|\mathcal{N}|\times |\mathcal{F}|}$,
and  $\M{H}_{NN} \in \mathbb{R}^{|\mathcal{N}|\times |\mathcal{N}|}$.
The damped Hessian search direction in~(\ref{finaldirection}) is computed from
the inverse of the reduced Hessian; that is,
$\M{B}_{F} = \M{H}_{FF}^{-1}$.

Let $\Mtilde{H}$ and $\Mtilde{B}$ denote the L-BFGS approximation to the
true and inverse Hessian.
To obtain the step $p^k$ we use the inverse approximation $\Mtilde{B}$,
then extract just the free variables for use in (\ref{lbfgs_final});
hence, we compute the search direction using the approximation $\Mtilde{B}_{F}$.
Assuming the Schur complement exists, this matrix is
\begin{equation*}
  \Mtilde{B}_{F} = ( \Mtilde{H}_{FF}
                    - \Mtilde{H}_{NF}^T \Mtilde{H}_{NN}^{-1} \Mtilde{H}_{NF} )^{-1} .
\end{equation*}
Comparing with the true reduced Hessian, we see the extra term
$\Mtilde{H}_{NF}^T \Mtilde{H}_{NN}^{-1} \Mtilde{H}_{NF}$,
a matrix of rank $|\mathcal{N}|$.
This is the influence in the L-BFGS approximation of variables
not in $\mathcal{F}$; we are effectively using the L-BFGS approximation
of the reduced inverse Hessian to compute the step.
Note that a small value of the tuning parameter $\epsilon$
in~(\ref{ActiveSetEpsilon}) can help reduce the size of $|\mathcal{N}|$,
lessening the influence.

%%%%%%%%%%%%%%%%%%%%%%%%%%%%%%%%
\subsection{Stopping Criterion}
%%%%%%%%%%%%%%%%%%%%%%%%%%%%%%%%

Since the row subproblems are convex, any point satisfying the first-order KKT
conditions is the optimal solution.  Specifically, $\V{b}^*$ is a KKT point of
(\ref{rsp-ptf}) if it satisfies
\begin{equation*}
  \nabla f_{\rm row}(\V{b}^*) - \V{\upsilon}^* = \V{0}, \ \
  (\V{b}^*)^T \V{\upsilon}^* = 0, \ \
  \V{b}^* \geq \V{0}, \ \
  \V{\upsilon}^* \geq \V{0},
\end{equation*}
where $\V{\upsilon}^*$ is the vector of dual variables associated with
the nonnegativity constraints.
Knowing the algorithm keeps all iterates $\V{b}^k$ nonnegative,
we can express the KKT condition for component $r$ as
\begin{equation*}
  \left| \min \{ b^k_r , \nabla_r f_{\rm row}(\V{b}^k) \} \right| = 0 .
\end{equation*}
A suitable stopping criterion is to approximately satisfy the
KKT conditions to a tolerance $\tau > 0$.
We achieve this by requiring that all row subproblems satisfy
\begin{equation}
  \kktviol = \max_r \left\{ \left| \min \{ b^k_r, \nabla_r f_{\rm row}(\V{b}^k) \}
                            \right| \right\}
           \leq \tau .
      \label{kktviolcond}
\end{equation}

The full algorithm solves to an overall tolerance $\tau$ when the
$\kktviol$ of every row subproblem satisfies (\ref{kktviolcond}).
This condition is enforced for all the row subproblems
(Step~\ref{algsubproblem:rowsubproblem} of Algorithm~\ref{alg:subproblem-sparse})
generated from all the tensor modes
(Step~\ref{algoutline:subproblem} of Algorithm~\ref{alg:outline}).
Note that enforcement requires examination of $\kktviol$ for all row subproblems
whenever the solution of any subproblem mode is updated, because the solution
modifies the $\Mn{\Pi}{n}$ matrices of other modes.

%%%%%%%%%%%%%%%%%%%%%%%%%%%%%%%%
\subsection{Row Subproblem Algorithms}
%%%%%%%%%%%%%%%%%%%%%%%%%%%%%%%%

Having described the ingredients, we pull everything together into complete
algorithms for solving the row subproblem in
Step~\ref{algsubproblem:rowsubproblem} of Algorithm~\ref{alg:subproblem-sparse}.
We present two methods in
Algorithm~\ref{alg:rowsubprobPDN} and Algorithm~\ref{alg:rowsubprobPQN}:
PDN-R uses a damped Hessian matrix,
and PQN-R uses a quasi-Newton Hessian approximation
(the `-R' designates a row subproblem formulation).
Both algorithms employ a two-metric projection framework for handling
bound constraints and a line search satisfying the Armijo condition.

\begin{algorithm}[t]
  \caption{Projected Newton-Based Solver for the Row Subproblem (PDN-R)}
  \label{alg:rowsubprobPDN}
  Given data $\Vhat{x}$ and $\M{\Pi}$, constants $\mu_0$, $\sigma$, $\beta$,
    $K_{\text{max}}$, stop tolerance $\tau$, and initial values $\V{b}^0$
  \\
  Return a solution $\V{b}^*$ to Step~\ref{algsubproblem:rowsubproblem}
  of Algorithm~\ref{alg:subproblem-sparse}

  \begin{algorithmic}[1]
    \For{$k = 0, 1, \dots, K_{\text{max}}$}
      \State Compute the gradient, $ \V{g}^k = \nabla f_{\rm row}(\V{b}^k)$,
             using $\Vhat{x}$ and $\M{\Pi}$ in (\ref{eq:grad_frow})

      \State Compute the first-order KKT violation
        \begin{equation*}
          \kktviol = \max_r \left\{ \left| \min \{ b^k_r, g^k_r \}
                            \right| \right\}
        \end{equation*}
      \If {$\kktviol \leq \tau$}
        \State \Return $\V{b}^* = \V{b}^k$
               \Comment{Converged to tolerance.}
      \EndIf	
      \State Find the indices of free variables from (\ref{ActiveSetDef})
             with $\epsilon = 10^{-3}$ in (\ref{ActiveSetEpsilon})
             \label{algPDN:activeset}
      \State Calculate the Hessian for free variables
        \begin{equation*}
          \M{H}^k_F = [\nabla^2 f_{\rm row}(\V{b}^k)]_{\mathcal{F}}
        \end{equation*}
      \State Compute the damped Newton direction
             $\V{d}_F^k = -(\M{H}^k_F + \mu_k I)^{-1} \V{g}^k_F$
               \label{algPDN:searchdir}
      \State Construct search direction $\V{d}^k$ over all variables
             using $\V{d}_F^k$ and $\V{g}^k$ in (\ref{finaldirection})
      \State Perform the projected line search (\ref{armijo})
             using $\sigma$ and $\beta$ to find step length $\alpha_k$
               \label{algPDN:linesearch}
      \State Update the current iterate
        \begin{equation*}
          \V{b}^{k+1} = P_{+}[\V{b}^k  + \alpha_k \V{d}^k]
        \end{equation*}
      \State Update the damping parameter $\mu_{k+1}$ according
             to (\ref{rho})-(\ref{updatedamping})
    \EndFor
    \State \Return $\V{b}^* = \V{b}^k$
           \Comment{Iteration limit reached.}
  \end{algorithmic}
\end{algorithm}

\begin{algorithm}[t]
  \caption{Projected Quasi-Newton Solver for the Row Subproblem (PQN-R)}
  \label{alg:rowsubprobPQN}
  Given data $\Vhat{x}$ and $\M{\Pi}$, constants $\mu_0$, $\sigma$, $\beta$,
    $K_{\text{max}}$, stop tolerance $\tau$, and initial values $\V{b}^0$
  \\
  Return a solution $\V{b}^*$ to Step~\ref{algsubproblem:rowsubproblem}
  of Algorithm~\ref{alg:subproblem-sparse}

  \begin{algorithmic}[1]
    \For{$k = 0, 1, \dots, K_{\text{max}}$}
      \State Compute the gradient, $ \V{g}^k = \nabla f_{\rm row}(\V{b}^k)$,
             using $\Vhat{x}$ and $\M{\Pi}$ in (\ref{eq:grad_frow})

      \State Compute the first-order KKT violation
        \begin{equation*}
          \kktviol = \max_r \left\{ \left| \min \{ b^k_r, g^k_r \}
                            \right| \right\}
        \end{equation*}
      \If {$\kktviol \leq \tau$}
        \State \Return $\V{b}^* = \V{b}^k$
               \Comment{Converged to tolerance.}
      \EndIf	
      \State
             Find the indices of free variables from (\ref{ActiveSetDef})
             with $\epsilon = 10^{-8}$ in (\ref{ActiveSetEpsilon})
             \label{algPQN:activeset}
      \State Construct search direction $\V{d}^k$ using $\V{g}^k$
             in (\ref{lbfgs_final})
      \State Perform the projected line search (\ref{armijo})
             using $\sigma$ and $\beta$ to find step length $\alpha_k$
               \label{algPQN:linesearch}
      \State Update the current iterate
        \begin{equation*}
          \V{b}^{k+1} = P_{+}[\V{b}^k  + \alpha_k \V{d}^k]
        \end{equation*}
      \State Update the L-BFGS approximation with $\V{b}^{k+1}$ and $\V{g}^{k+1}$
               \label{algPQN:lbfgs}
    \EndFor
    \State \Return $\V{b}^* = \V{b}^k$
           \Comment{Iteration limit reached.}
  \end{algorithmic}
\end{algorithm}

As mentioned, PQN-R is related to~\cite{KimSuvritDhillon:BoxConstraints}.
Specifically, we note
\begin{enumerate}
  \item
    The free variables chosen in Step~\ref{algPQN:activeset} of PQN-R
    are found with $\epsilon = 10^{-8}$ in (\ref{ActiveSetEpsilon}),
    while~\cite{KimSuvritDhillon:BoxConstraints} effectively uses
    $\epsilon = 0$ for an instantaneous line search.  As noted in
    Section~\ref{subsec:subprobActiveSet}, convergence is guaranteed
    when $\epsilon > 0$.
    Step~\ref{algPDN:activeset} of PDN-R uses the default value
    $\epsilon = 10^{-3}$ because we find it generally leads to faster convergence.
  \item
    The line search in Step~\ref{algPQN:linesearch} of PQN-R and
    Step~\ref{algPDN:linesearch} of PDN-R satisfies the Armijo condition.
    This differs from~\cite{KimSuvritDhillon:BoxConstraints},
    which used $\sigma \alpha (\V{d}^k)^T \nabla f_{\rm row}(\V{b}^k)$ on the
    right-hand side of (\ref{armijo}).  We use (\ref{armijo}) because it
    correctly measures predicted progress.  In particular, it is easier to
    satisfy when $(\V{d}^k)^T \nabla f_{\rm row}(\V{b}^k)$ is large and many
    variables hit their bound for small $\alpha$.
  \item
    Updates to the L-BFGS approximation in Step~\ref{algPQN:lbfgs} of PQN-R
    are unchanged from~\cite{KimSuvritDhillon:BoxConstraints}.
    Information is included from all row subproblem variables, whether
    active or free.

\end{enumerate}

We express the computational cost of PDN-R and PQN-R in terms of the cost
per iteration of Algorithm~\ref{alg:outline}; that is, the cost of
executing Steps~\ref{algoutline:for} through \ref{algoutline:endfor}.
The matrix $\M{\Pi}^{(n)}$ is formed for the row subproblems of every mode,
with the cost for each mode
proportional to the number of nonzeros in the data tensor, nnz$(\TX)$.
This should dominate the cost of reweighting factor matrices in
Steps~\ref{algoutline:rescale1} and \ref{algoutline:rescale2}.
The $n$th mode solves $I_n$ convex row subproblems, each with $R$
unknowns, using Algorithm~\ref{alg:rowsubprobPDN} (PDN-R)
or Algorithm~\ref{alg:rowsubprobPQN} (PQN-R).
Row subproblems execute over at most $K_{\text{max}}$ inner iterations.
Near a local minimum we expect PDN-R to take fewer inner iterations than
PQN-R because the damped Newton method converges asymptotically at a
quadratic rate, while L-BFGS convergence is at best R-linear.
However, the cost estimate
will assume the worst case of $K_{\text{max}}$ iterations for all row subproblems.
The dominant cost of Algorithm~\ref{alg:rowsubprobPDN} is solution of the
damped Newton direction in Step~\ref{algPDN:searchdir}, which
costs $O(R^3)$ operations to solve the $R \times R$ dense linear system.
Hence, the cost per iteration of PDN-R is
\begin{equation}
\label{cost-pdnr}
  N \cdot O(\text{nnz}(\TX)) + K_{\text{max}} \cdot O(R^3)
                                              \cdot \sum_{n=1}^N I_n .
\end{equation}
The dominant costs of Algorithm~\ref{alg:rowsubprobPQN} are computation of
the search direction and updating the L-BFGS matrix, both $O(R)$ operations.
Hence, the cost per iteration of PQN-R is
\begin{equation}
\label{cost-pqnr}
  N \cdot O(\text{nnz}(\TX)) + K_{\text{max}} \cdot O(R)
                                              \cdot \sum_{n=1}^N I_n .
\end{equation}

\section{Experiments}  \label{sec:exp}

This section characterizes the performance of our algorithms, comparing
them with multiplicative update \cite{ChiKolda:cpapr}
and second-order methods that do not use the row subproblem formulation.
All algorithms fit in the alternating block framework of
Algorithm~\ref{alg:outline}, differing in how they solve~(\ref{modeNprob})
in Step~\ref{algoutline:subproblem}.

Our two algorithms are the projected damped Hessian method (PDN-R)
and the projected quasi-Newton method (PQN-R), from
Algorithms~\ref{alg:rowsubprobPDN} and~\ref{alg:rowsubprobPQN}, respectively.
Recall that `-R' means the row subproblem formulation is applied.
In this paper we do not tune the algorithms to each test case,
but instead chose a single set of parameter values:
$\mu_0 = 10^{-5}$, $\sigma = 10^{-4}$, and $\beta = 1/2$.
The bound constraint threshold in PDN-R from~(\ref{ActiveSetEpsilon})
was set to $\epsilon = 10^{-3}$ for PDN-R
and $\epsilon = 10^{-8}$ for PQN-R, values that are observed to give
best algorithm performance.
The L-BFGS approximations in PQN-R stored the $M=3$ most recent
update pairs (\ref{curvPairs}).

The multiplicative update (MU) algorithm that we compare with is that of Chi
and Kolda \cite{ChiKolda:cpapr}, available as function {\tt cp\_apr}
in the Matlab Tensor Toolbox \cite{BaderKolda:MatlabTTB}.
It builds on tensor generalizations of the
Lee and Seung method, specifically treating \emph{inadmissible zeros}
(their term for factor elements that are active but close to zero) to
improve the convergence rate.  Algorithm MU can be tuned by selecting the
number of inner iterations for approximately solving the subproblem at
Step~\ref{algoutline:subproblem} of Algorithm~\ref{alg:outline}.
We found that ten inner iterations worked well in all experiments.

We also compare to a projected quasi-Newton (PQN) algorithm adopted from
Kim et al.\@~\cite{KimSuvritDhillon:BoxConstraints}.
PQN is similar to PQN-R but solves~(\ref{modeNprob})
without reformulating the block subproblem into row subproblems.
PQN identifies
the active set using $\epsilon = 0$ in (\ref{ActiveSetEpsilon}) and maintains
a limited-memory BFGS approximation of the Hessian.
However, PQN uses one L-BFGS matrix for the entire subproblem, storing
the three most recent update pairs.
We used Matlab code from the authors of \cite{KimSuvritDhillon:BoxConstraints},
embedding it in the alternating framework of Algorithm~\ref{alg:outline},
with the modifications described in Section~\ref{subsec:PQNmods}.

Additionally,
we compare PDN-R to a projected damped Hessian (PDN) method that uses one matrix
for the block subproblem instead of a matrix for every row subproblem.
PDN exploits the block diagonal nature of the
Hessian to construct a search direction for the same computational cost
as PDN-R; i.e., one search direction of PDN takes the same effort as
computing one search direction for all row subproblems in PDN-R.
Similar remarks apply to computation of the objective function for the
subproblem (\ref{modeNprob}).
However, PDN applies a single damping parameter $\mu_k$ to
the block subproblem Hessian and updates all variables in the block subproblem
from a single line search along the search direction.

All algorithms were coded in Matlab using the sparse tensor objects of the
Tensor Toolbox \cite{BaderKolda:MatlabTTB}.  All experiments were performed
on a Linux workstation with 12GB memory.  Data sets were large enough
to be demanding but small enough to fit in machine memory; hence, performance
results are not biased by disk access issues.

The experiments that follow show three important results, as follows.
\begin{enumerate}
  \item The row subproblem formulation is better suited to second-order
        methods than the block subproblem formulation because it controls the
        number of iterations for each row subproblem independently, and
        because its convergence is more robust.
  \item PDN-R and PQN-R are faster than the other algorithms in terms of
        in reducing the $\kktviol$, especially when solving to high accuracy.
        This holds for any number of components.  PQN-R becomes faster than
        PDN-R as the number of components increases.
  \item PDN-R and PQN-R reach good solutions with high sparsity more quickly
        than the other algorithms, a desirable feature when the factor matrices
        are expected to be sparse.
\end{enumerate}

In Section~\ref{subsec:exp-subprob} we report only
performance in solving a single block subproblem~(\ref{modeNprob})
since the time is representative of the total time it will take to solve
the full tensor factorization problem~(\ref{fullprob}).
In Section~\ref{subsec:exp-fullprob} we report results from solving the full
problem within the alternating block framework (Algorithm 1).

%%%%%%%%%%%%%%%%%%%%%%%%%%
\subsection{Solving the Convex Block Subproblem}  \label{subsec:exp-subprob}
%%%%%%%%%%%%%%%%%%%%%%%%%%%

We begin by examining algorithm performance on the convex
subproblem~(\ref{modeNprob}) of the alternating block framework.
Here we look at a single
representative subproblem.  Our goal is to characterize the relative behavior
of algorithms on the representative block subproblem.

Appendix~\ref{app-generate} describes our method for generating synthetic
test problems with reasonable sparsity.  We investigate a three-way tensor of size
$200 \times 300 \times 400$, generating $S = 500,000$ data samples.
The number of components,
$R$, is varied over the set $\{ 20, 40, 60, 80, 100 \}$.  For each value of $R$,
the procedure generates a sparse multilinear model
$\TM = \KTsmall{\Vl; \Mn{A}{1}, \Mn{A}{2}, \Mn{A}{3}}$ and data tensor $\TX$.
Table~\ref{tbl:subproblem-sizes} lists the number of nonzero elements
found in the data tensor $\TX$ that results from $500,000$ data samples,
averaged over ten random seeds.  The number of nonzeros, a key determiner
of algorithm cost in equations~(\ref{cost-pdnr}) and (\ref{cost-pqnr}),
is approximately the same for all values of $R$.
\begin{table}[ht]
  \centering
  \small
  \caption{Subproblem sparsity for number of components $R$}
      \label{tbl:subproblem-sizes}

  \begin{tabular}{rcc}
    $R$ & \textsc{Number Nonzeros} & \textsc{Density}  \\
    \hline
     20  & 413,460  &  1.72\%  \\
     40  & 450,760  &  1.88\%  \\
     60  & 464,440  &  1.94\%  \\
     80  & 470,950  &  1.96\%  \\
    100  & 475,450  &  1.98\%
  \end{tabular}
\end{table}

We consider just the subproblem obtained by unfolding along mode 1;
hence, the test case contains 200 row subproblems of the form~(\ref{rsp-ptf}).
To solve just the mode-1 subproblem, the {\tt for} loop
at Step~\ref{algoutline:for} of
Algorithm~\ref{alg:outline} is changed to $n = 1$.

We run several trials of the subproblem solver from different initial
guesses of the unknowns, holding $\Mn{A}{2}$ and $\Mn{A}{3}$ from $\TM$ constant.
The initial guess draws each element of $\Mn{A}{1}$ from a uniform distribution
on $[0,1)$ and sets each element of $\V{\lambda}$ to one.
To satisfy constraints in (\ref{fullprob}), the columns
of $\Mn{A}{1}$ are normalized and the normalization
factor is absorbed into $\V{\lambda}$.
The mode-1 subproblem~(\ref{modeNprob}) is now defined with
$\M{\Pi}=(\Mn{A}{3} \odot \Mn{A}{2})^T$, $\M{X}=\M{X}_{(1)}$,
and $\M{B}=\Mn{A}{1}\M{\Lambda}$, with unknowns $\M{B}$ initialized
using the initial guess for $\V{\lambda}$ and $\Mn{A}{1}$.

%%%%%%%%%%%%%%%%%%%%%%%%%%
\subsubsection{PDN-R and PDN on the Convex Subproblem}  \label{subsec:exp-PDNRsub}
%%%%%%%%%%%%%%%%%%%%%%%%%%%

We first characterize the behavior of our Newton-based algorithm, PDN-R,
and compare it with PDN.
Row subproblems are solved using Algorithm~\ref{alg:rowsubprobPDN}
with stop tolerance $\tau = 10^{-8}$ and the parameter values mentioned
at the beginning of Section~\ref{sec:exp}.
The value of $K_{\text{max}}$ in Algorithm~\ref{alg:rowsubprobPDN}
is large enough that the $\kktviol$ converges to $\tau$ before $K_{\text{max}}$
is reached.

%\begin{figure}[p]
\begin{figure}
\centering
\subfloat[$R=20$]{\includegraphics[scale=0.5]
    {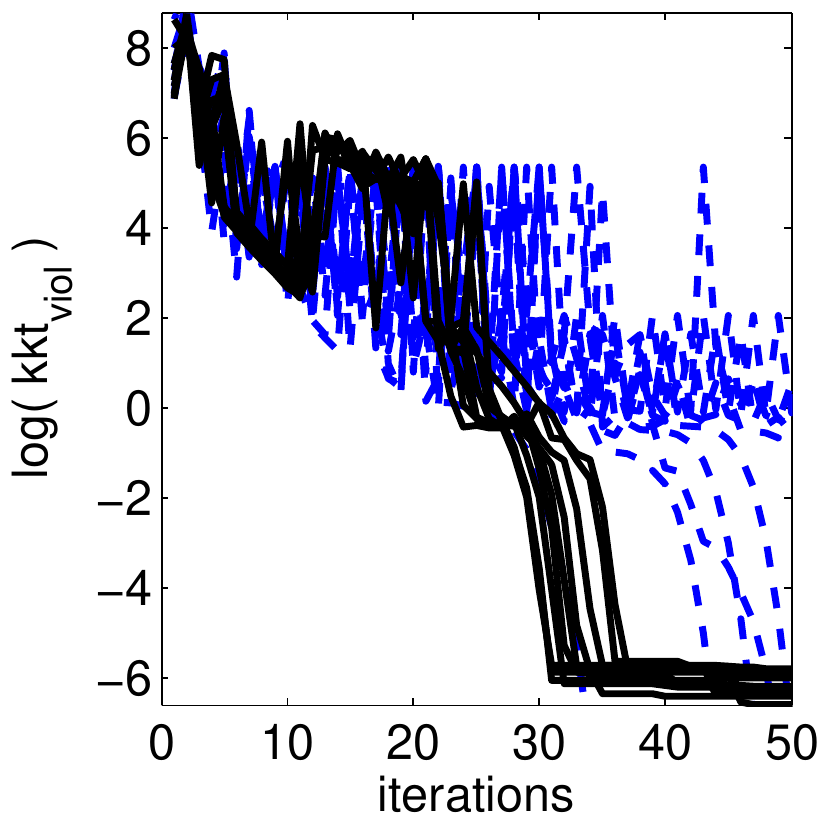}  \label{fig:subprob-Newton-kktR20}}
~~
\subfloat[$R=60$]{\includegraphics[scale=0.5,trim=0 0 0  0,clip]
    {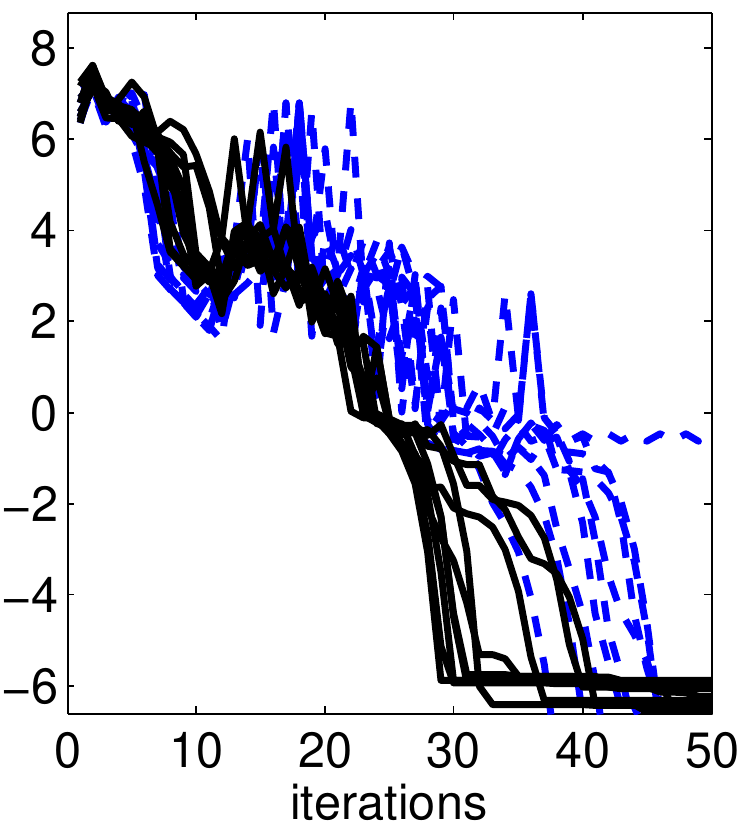}  \label{fig:subprob-Newton-kktR60}}
~~
\subfloat[$R=100$]{\includegraphics[scale=0.5,trim=0 0 0 0,clip]
    {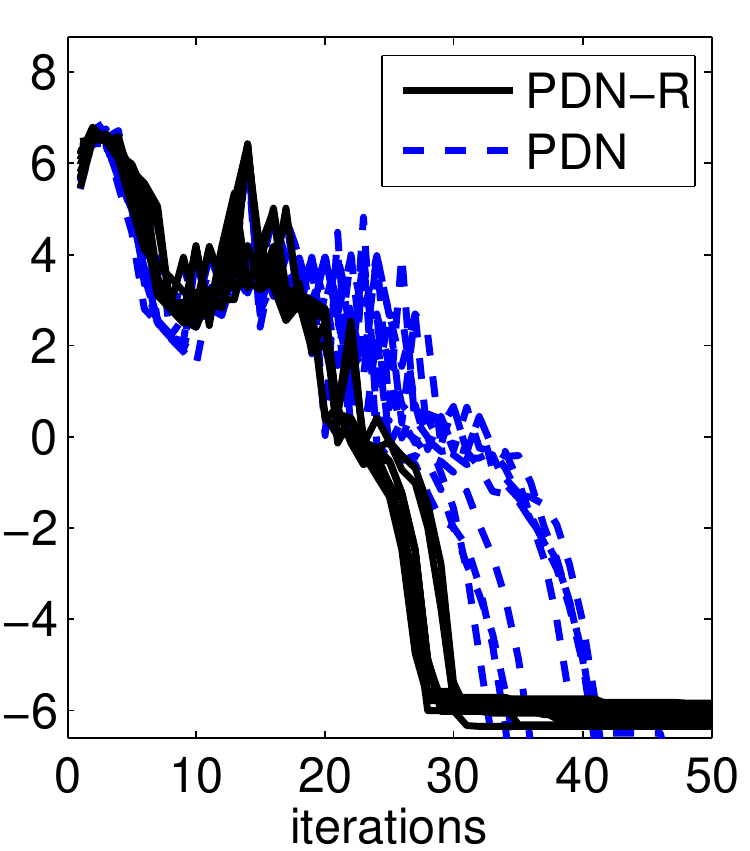}  \label{fig:subprob-Newton-kktR100}}
\\
\subfloat[$R=20$]{\includegraphics[scale=0.5]
    {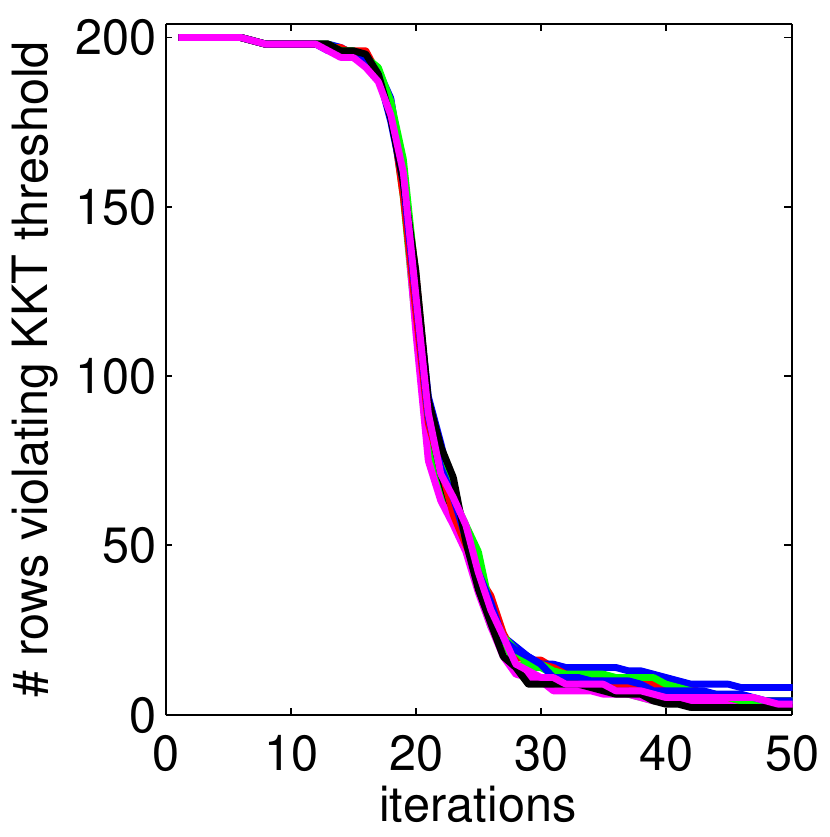}  \label{fig:subprob-Newton-rowsatR20}}
~~
\subfloat[$R=60$]{\includegraphics[scale=0.5,trim=0 0 0  0,clip]
    {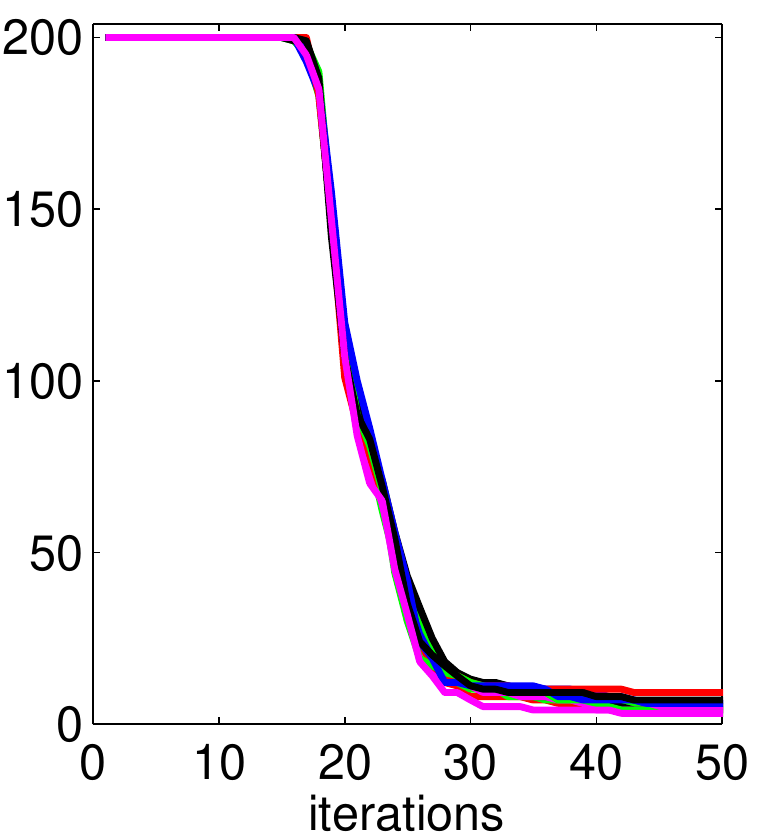}  \label{fig:subprob-Newton-rowsatR60}}
~~
\subfloat[$R=100$]{\includegraphics[scale=0.5,trim=0 0 0  0,clip]
    {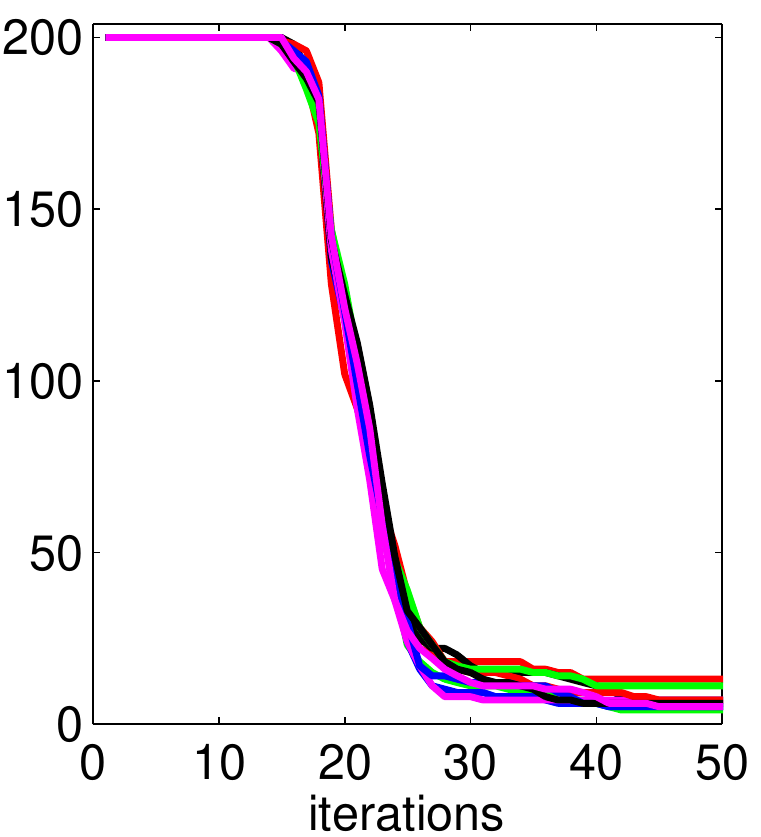}  \label{fig:subprob-Newton-rowsatR100}}
  \caption{Convergence behavior of PDN-R (Algorithm~\ref{alg:rowsubprobPDN})
           and PDN over ten runs with different start points, for three values
           of $R$.
           The upper graphs plot $\log(\kktviol)$, showing how the maximum
           violation over all row subproblems varies as the number of
           iterations increases.  Solid lines are PDN-R, dashed lines are PDN.
           The lower graphs plot the number of rows violating the KKT-based
           stop tolerance (ten runs of only PDN-R).}
\end{figure}

%\begin{figure}[p]
\begin{figure}
\centering
\subfloat[Number of zero entries in $\Mn{A}{1}$ for $R=100$]
    {\includegraphics[scale=0.5]{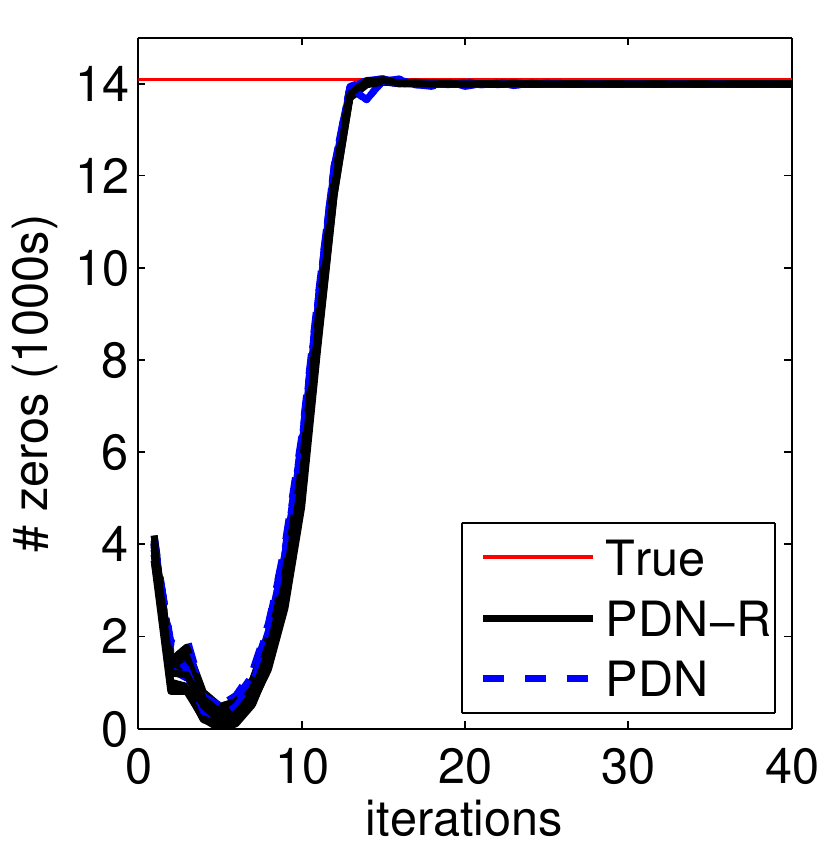}
     \label{fig:subprob-Newton-nzero}}
~~
\subfloat[Cumulative execution time for $R=100$]
    {\includegraphics[scale=0.5]{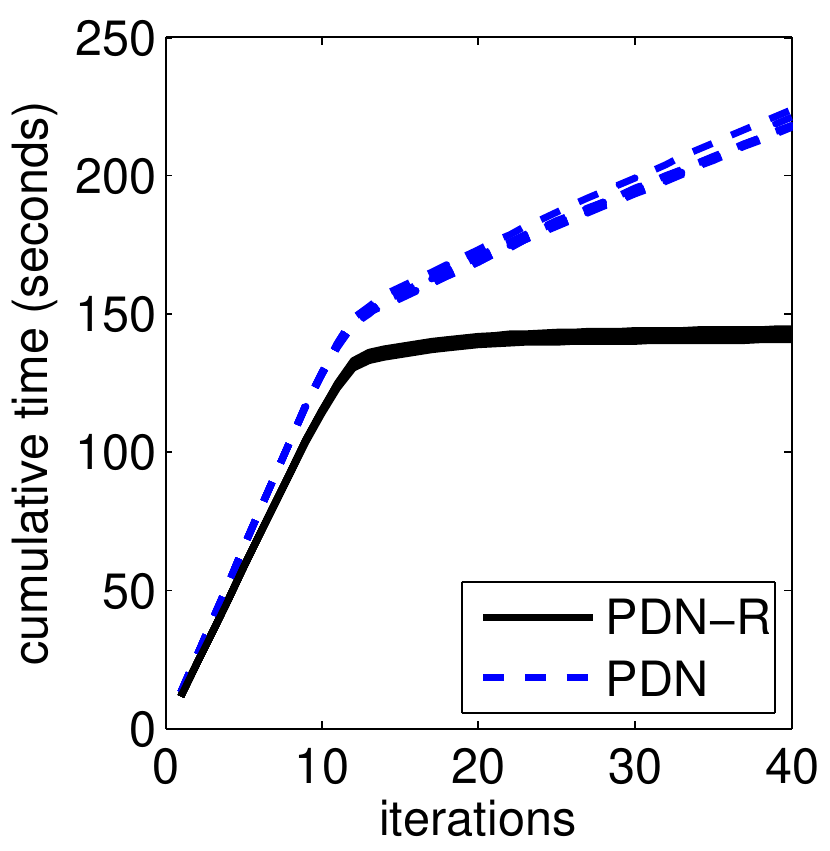}
     \label{fig:subprob-Newton-time}}
  \caption{Additional convergence behavior over ten runs with different start
           points.  Solid lines are PDN-R, dashed lines are PDN.
           The red horizontal line shows the true number of elements
           exactly equal to zero.}
      \label{fig:subprob-Newton-results-2}
\end{figure}

Figures~\ref{fig:subprob-Newton-kktR20}~-~\ref{fig:subprob-Newton-kktR100}
show how
KKT violations decrease with iteration for three different values of $R$.
The subproblem was solved ten times from different randomly chosen start points.
(Since the subproblem is strictly convex, there is a single unique minimum that
is reached from every start point.)
Each solid line plots the maximum $\kktviol$ over all 200 row subproblems for
one of the ten PDN-R runs.
Each dashed line plots the $\kktviol$ of the block subproblem for one of
the ten PDN runs.
Note the $y$-axis is the $\log_{10}$ of $\kktviol$.
The figure demonstrates that after some initial slow progress, both algorithms
exhibit the fast quadratic convergence rate typical of Newton methods.
PDN-R clearly takes fewer iterations to compute a factorization
with small $\kktviol$.

Figures~\ref{fig:subprob-Newton-rowsatR20}~-~\ref{fig:subprob-Newton-rowsatR100}
show the number of row subproblems in PDN-R
that satisfy the KKT-based stop tolerance after a given number
of iterations.  Remember that \emph{all} row subproblems must satisfy the
KKT tolerance before the algorithm declares a solution.

Figure~\ref{fig:subprob-Newton-results-2} shows additional features of the
convergence, for just the case of $R=100$ components (behavior is similar
for other values of $R$).
In Figure~\ref{fig:subprob-Newton-nzero}
we see the number of elements of $\Mn{A}{1}$ exactly equal to zero.
Data for this experiment was generated stochastically
from sparse factor matrices (see Appendix~\ref{app-generate});
hence, we expect a sparse solution.
The plot indicates that sparsity can be achieved after reducing
$\kktviol$ to a moderately small tolerance (around $10^{-2}$ in this example).
We return to sparsity of the solution in the sections below.

In Figure~\ref{fig:subprob-Newton-time}
we see that execution time per iteration decreases when variables are
closer to a solution.
PDN-R execution time becomes very small because only
a few row subproblems need to satisfy the convergence tolerance, and
only these are updated.
PDN takes more time per iteration because it computes a single
search direction for updating all variables in the block subproblem,
even when most of the variables are near an optimal value.
These experiments show that PDN-R and PDN behave similarly for the convex
subproblem and that PDN-R is a little faster; much larger differences appear
when the full factorization is computed in Section~\ref{subsec:exp-PDNRfull}.

%%%%%%%%%%%%%%%%%%%%%%%%%%
\subsubsection{PQN-R and PQN on the Convex Subproblem}  \label{subsec:exp-PQNRsub}
%%%%%%%%%%%%%%%%%%%%%%%%%%%

In this section we demonstrate the importance of the row subproblem
formulation by comparing PQN-R with PQN,
showing the huge speedup achieved with our row subproblem formulation.
We compare the algorithms on the mode-1 subproblem described above,
from the same ten random initial guesses for $\Mn{A}{1}$.
Table~\ref{tbl:subproblem-QN-methods}
lists the average CPU times over ten runs.
PQN-R was executed until the $\kktviol$ was less
than $\tau = 10^{-8}$.  PQN was unable to achieve this level of accuracy, so
execution was stopped at a tolerance of $10^{-3}$.
Results in the table show that PQN-R is much faster at decreasing the KKT
violation.  We note that a KKT violation of $10^{-8}$ is approximately the
square root of machine epsilon, the smallest practical value that can be
attained.

\begin{table}[ht]
  \centering
  \small
  \caption{Convergence comparison in solving the subproblem}
      \label{tbl:subproblem-QN-methods}
  \begin{tabular}{r|rr|rr}
        & \multicolumn{2}{c|}{Algorithm PQN}
        & \multicolumn{2}{c}{Algorithm PQN-R}   \\
    $R$ & {$\kktviol < 10^{-1}$} & {$\kktviol < 10^{-3}$ \;}
        & {$\kktviol < 10^{-1}$} & {$\kktviol < 10^{-8}$}   \\
    \hline
     20  &  625 secs &  690 secs & 12.4 secs & 17.1 secs   \\
     40  &  755 secs &  846 secs & 10.9 secs & 16.4 secs   \\
     60  &  822 secs &  920 secs & 11.3 secs & 16.8 secs   \\
     80  & 1022 secs & 1141 secs & 13.7 secs & 19.5 secs   \\
    100  &  993 secs & 1125 secs & 13.1 secs & 20.2 secs
  \end{tabular}
\end{table}

\begin{figure}[htbp]
\centering
\subfloat[$R=20$]{\includegraphics[scale=0.5]{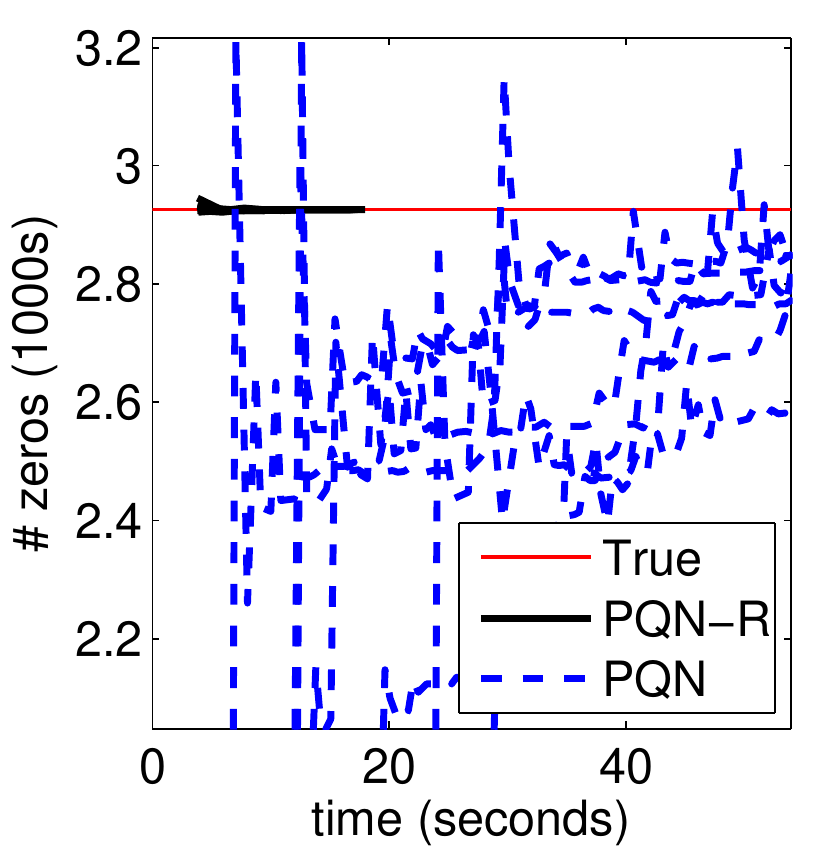}}
~~
\subfloat[$R=60$]{\includegraphics[scale=0.5]{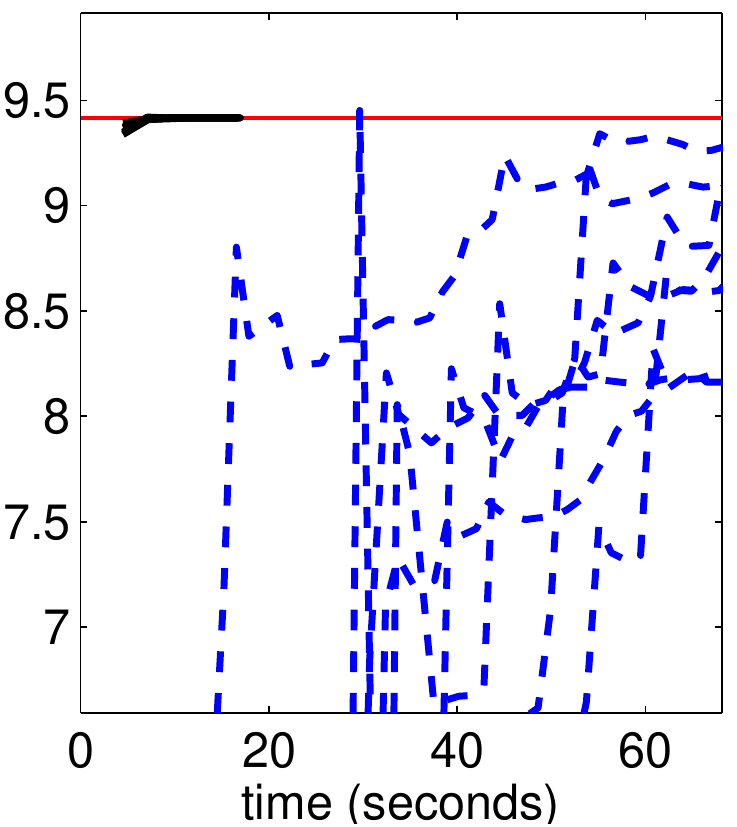}}
~~
\subfloat[$R=100$]{\includegraphics[scale=0.5]{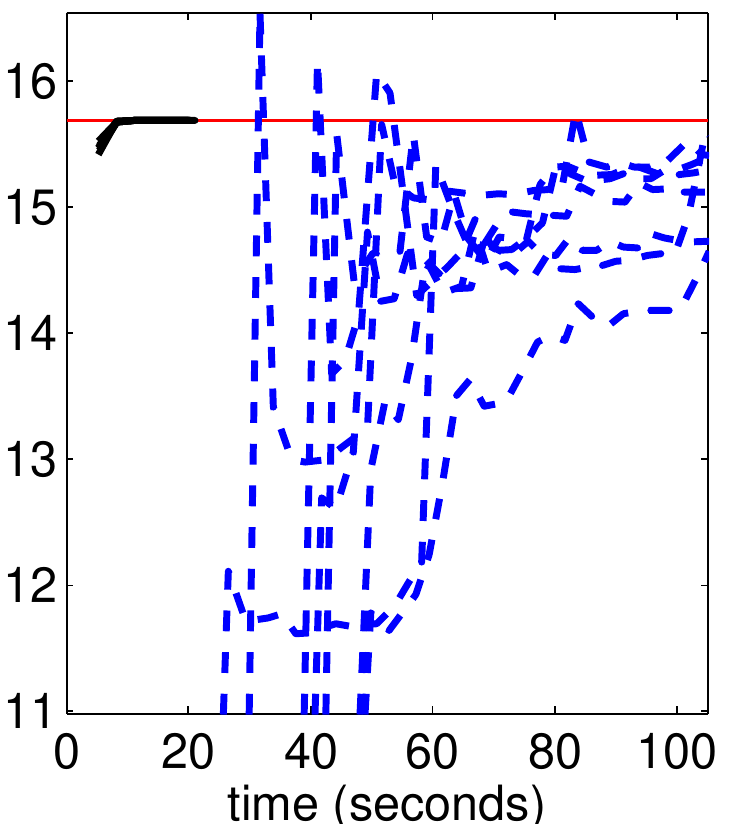}}
  \caption{Number of elements equal to zero in $\Mn{A}{1}$ found by PQN-R (solid lines,
           forming a short segment in the upper left)
           and PQN (dashed lines) as a function of compute time,
           for a subproblem with different values of $R$.
           The same subproblem was solved from six different start points, corresponding to the different colors.
           The red horizontal line shows the true value.}
      \label{fig:subprob-compare-QN-zeros}
\end{figure}

The two algorithms also differ in how they discover the number of elements
in $\Mn{A}{1}$ equal to zero.  Both eventually agree on the number of zero
elements, but PQN-R is much faster.
Figure~\ref{fig:subprob-compare-QN-zeros} shows the progress made by the two
algorithms; the behavior of PQN for this quantity is erratic and slow
to converge.

Algorithm PQN might be relatively more competitive
for tensor subproblems with a small number of rows.  Nevertheless,
it is apparent that applying L-BFGS to the block subproblem does not work
as well as applying separate instances of L-BFGS to the row subproblems.
This is not surprising since the first method ignores the block diagonal
structure of the true Hessian.
We see no advantages to using PQN and do not consider it further.

%%%%%%%%%%%%%%%%%%%%%%%%%%%
\subsubsection{PDN-R, PQN-R, and MU on the Convex Subproblem}
\label{subsec:3algs-subprob}
%%%%%%%%%%%%%%%%%%%%%%%%%%%

Next we compare our new row-based algorithms, PDN-R and PQN-R, with the
multiplicative update method MU \cite{ChiKolda:cpapr}.
Again we use the mode-1 subproblem of
Section~\ref{subsec:exp-subprob}, from the same ten random initial guesses.

As described in Section~\ref{sec:exp}, MU is a state-of-the-art representative
of the most common algorithm for nonnegative tensor factorization.  It is
a form of scaled steepest-descent with bound constraints
\cite{seung2001algorithms},
and therefore is expected to converge more slowly than Newton or quasi-Newton
methods.
We see this clearly in Table~\ref{tbl:subproblem-3-methods} for two different
stop tolerances. The MU algorithm was executed with a time limit of 1800 seconds
per problem, and failed to reach $\kktviol < 10^{-3}$ before this limit when
$R$ was 60 or larger.

\begin{table}[ht]
  \centering
  \small
  \caption{Time to reach stop tolerance for three algorithms
           (averaged over ten runs)}
      \label{tbl:subproblem-3-methods}
  \begin{tabular}{r|r|r|r|r|r|r}
        & \multicolumn{3}{c|}{$\kktviol = 10^{-2}$}
        & \multicolumn{3}{c}{$\kktviol = 10^{-3}$}   \\
    $R$ & \textsc{PDN-R} & \textsc{PQN-R} & \textsc{MU}
        & \textsc{PDN-R} & \textsc{PQN-R} & \textsc{MU}  \\
    \hline
    20 &   8.1 secs & 14.5 secs &  97.7 secs &   8.1 secs & 15.6 secs &  161.3 secs  \\
    40 &  25.1 secs & 13.1 secs & 239.2 secs &  25.2 secs & 14.6 secs &  485.9 secs  \\
    60 &  53.6 secs & 13.8 secs & 469.2 secs &  53.7 secs & 15.6 secs & $>$1800 secs  \\
    80 &  92.8 secs & 16.3 secs & 455.4 secs &  92.9 secs & 18.1 secs & $>$1800 secs  \\
   100 & 139.8 secs & 16.0 secs & 730.7 secs & 140.0 secs & 18.3 secs & $>$1800 secs
  \end{tabular}
\end{table}

\begin{figure}[htbp]
\centering
\subfloat[$R=20$]{\includegraphics[scale=0.5]{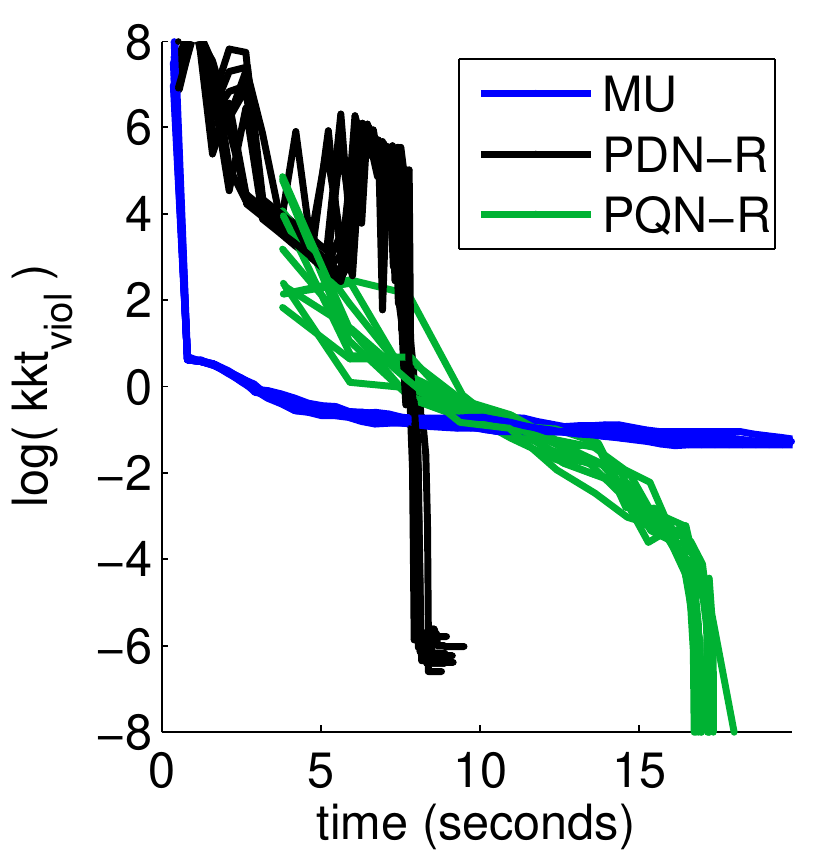}}
~~
\subfloat[$R=60$]{\includegraphics[scale=0.5,trim=0 0 0  0,clip]{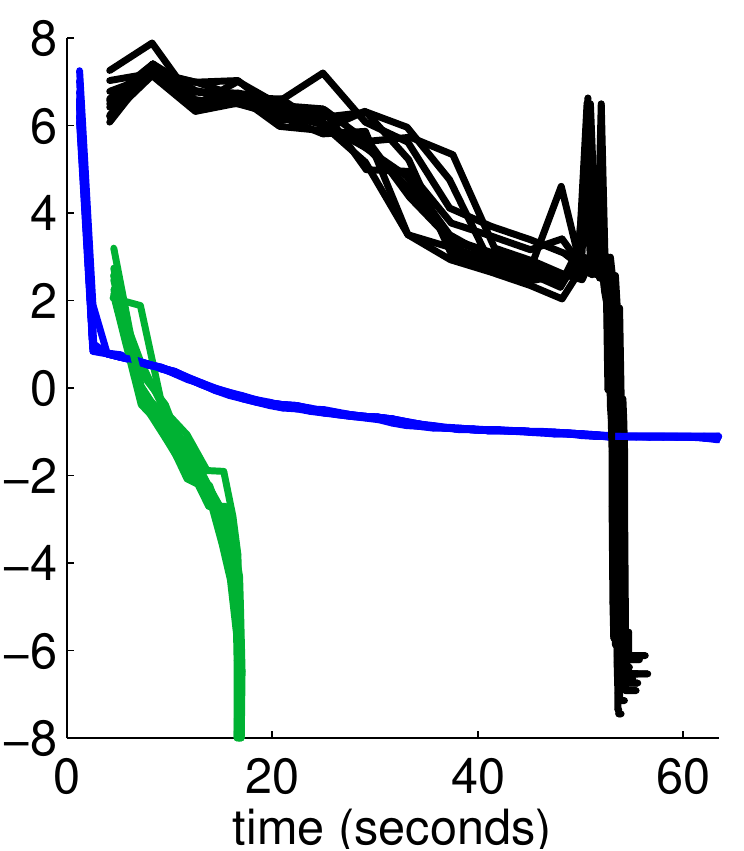}}
~~
\subfloat[$R=100$]{\includegraphics[scale=0.5,trim=0 0 0  0,clip]{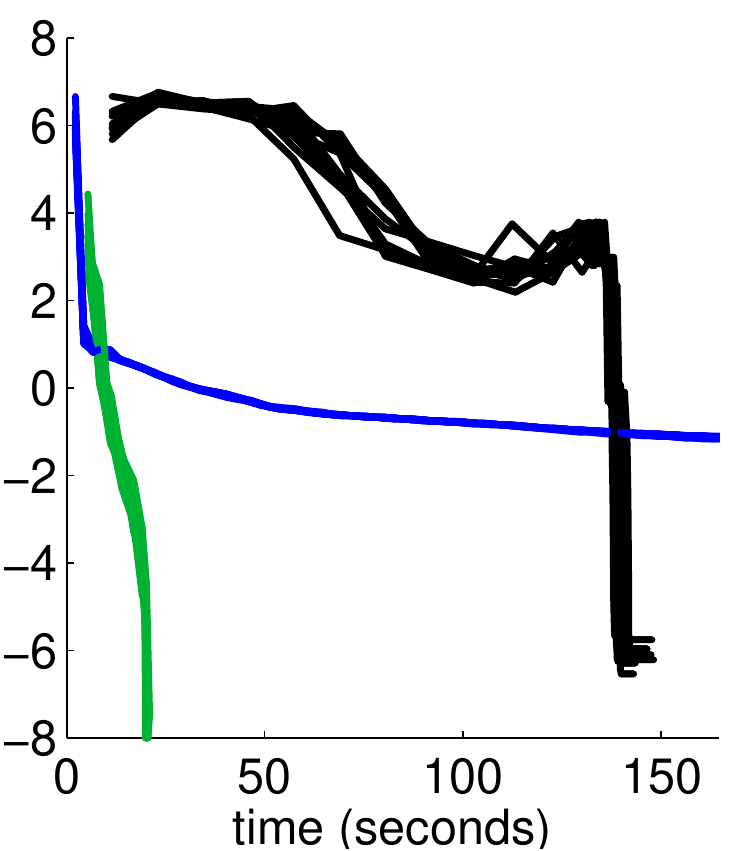}}
  \caption{Convergence behavior comparison on subproblem
           for different values of $R$ (ten runs each). Algorithm MU (blue) makes
           fast initial progress in reducing the violation, but slows
           dramatically after reaching a violation of about 0.1.
           PDN-R (black) and PQN-R (green) reduce the
           violation much further, with PQN-R being faster for higher values of $R$.}
      \label{fig:subprob-compare-KKT}
\end{figure}

Of course, the disparity in convergence time is more pronounced when a smaller
KKT error is demanded.  Figure~\ref{fig:subprob-compare-KKT} shows the
decrease in KKT violation as a function of compute time.  Here we see that
MU makes a faster initial reduction in KKT violation than PDN-R or PQN-R,
but then it slows to a linear rate of convergence.
Notice the gap from time zero for PDN-R and PQN-R, which reflects setup cost
before the first iteration result is computed.  For PQN-R the setup time
is fairly constant with $R$ (about 3.8 seconds), while PDN-R has a setup time
that increases with $R$ (11.5 seconds for $R=100$).  Unlike MU, both algorithms
must construct software structures for all row subproblems before a first
iteration result appears.
Figure~\ref{fig:subprob-compare-KKT} also reveals that PDN-R is slower
relative to PQN-R as the number of components $R$ increases.
This is because the cost of
solving a Newton-based Hessian is $O(R^3)$,
while the limited-memory BFGS Hessian cost is $O(R)$.

\begin{figure}[htbp]
\centering
\subfloat[$R=20$]{\includegraphics[scale=0.5]{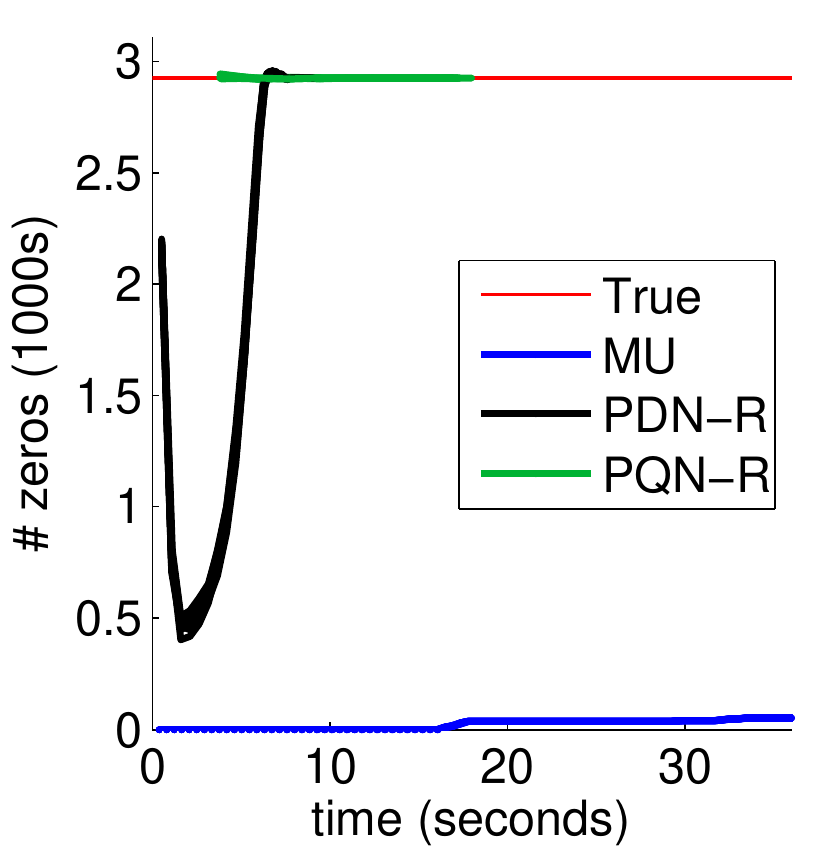}}
~~
\subfloat[$R=60$]{\includegraphics[scale=0.5,trim=0 0 0  0,clip]{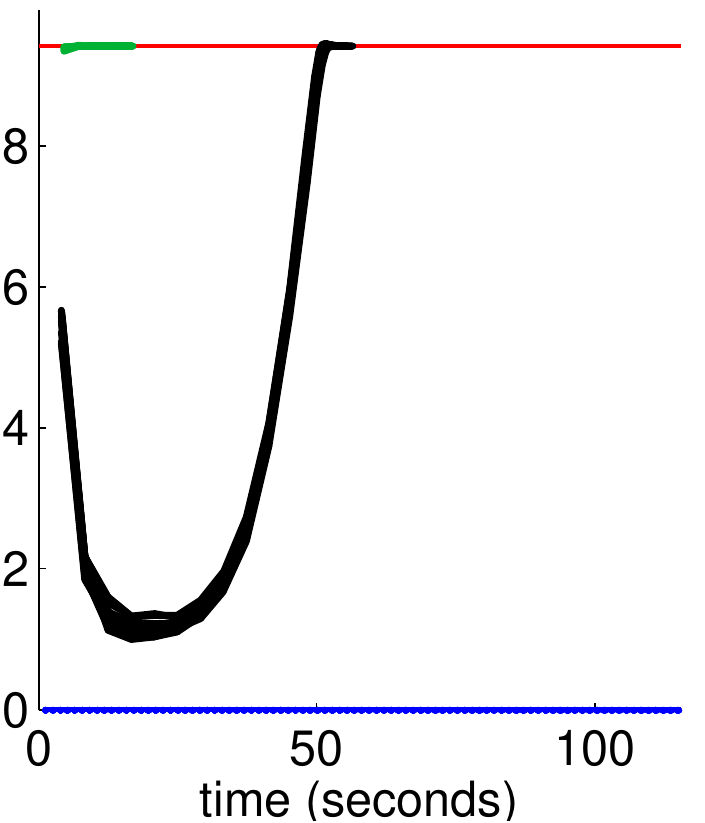}}
~~
\subfloat[$R=100$]{\includegraphics[scale=0.5,trim=0 0 0  0,clip]{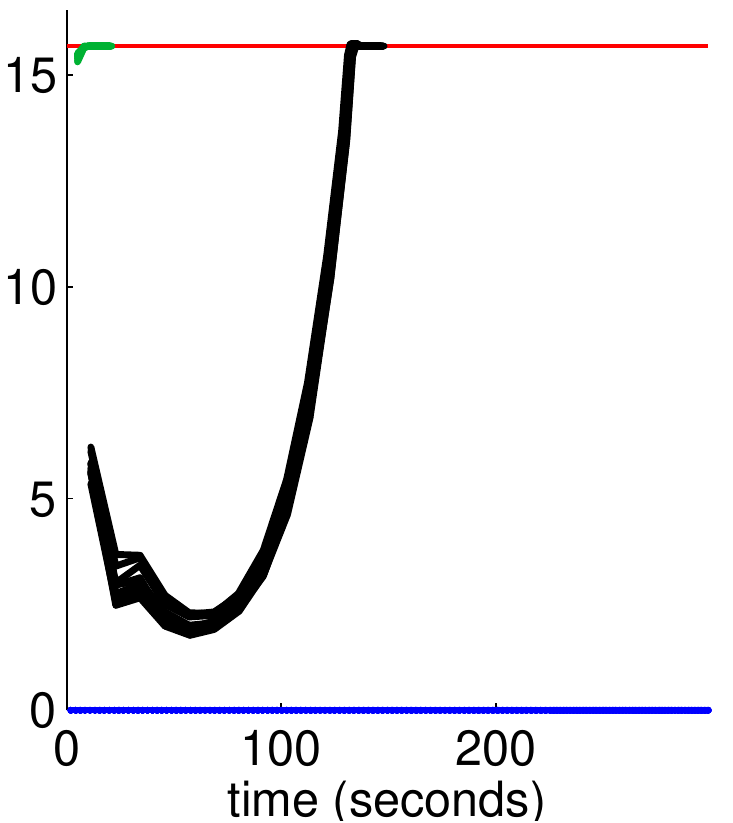}}
  \caption{Effectiveness of the three algorithms in finding a sparse solution for
           different values of $R$.  In each case the number of elements in
           $\Mn{A}{1}$ equal to zero is plotted against execution time.
           The PDN-R (black) and PQN-R (green) algorithms are
           much faster than MU (blue) at finding the zeros.
           The horizontal red line shows the true value.}
      \label{fig:subprob-compare-nzero}
\end{figure}

\begin{figure}[htbp]
\centering
\subfloat[$R=20$]{\includegraphics[scale=0.5]{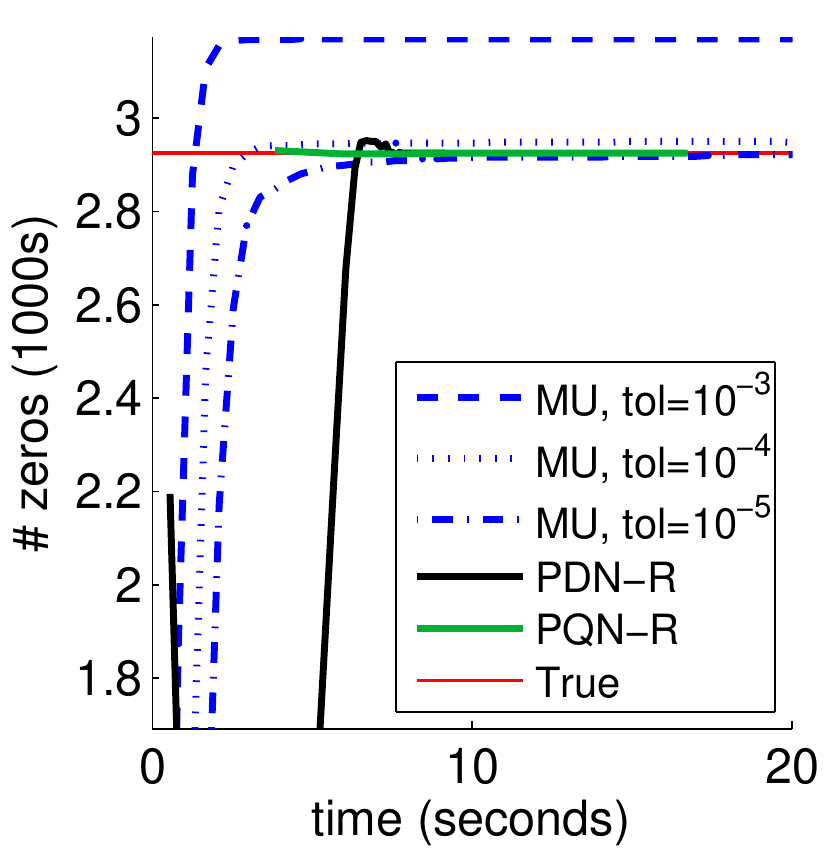}}
~~
\subfloat[$R=60$]{\includegraphics[scale=0.5]{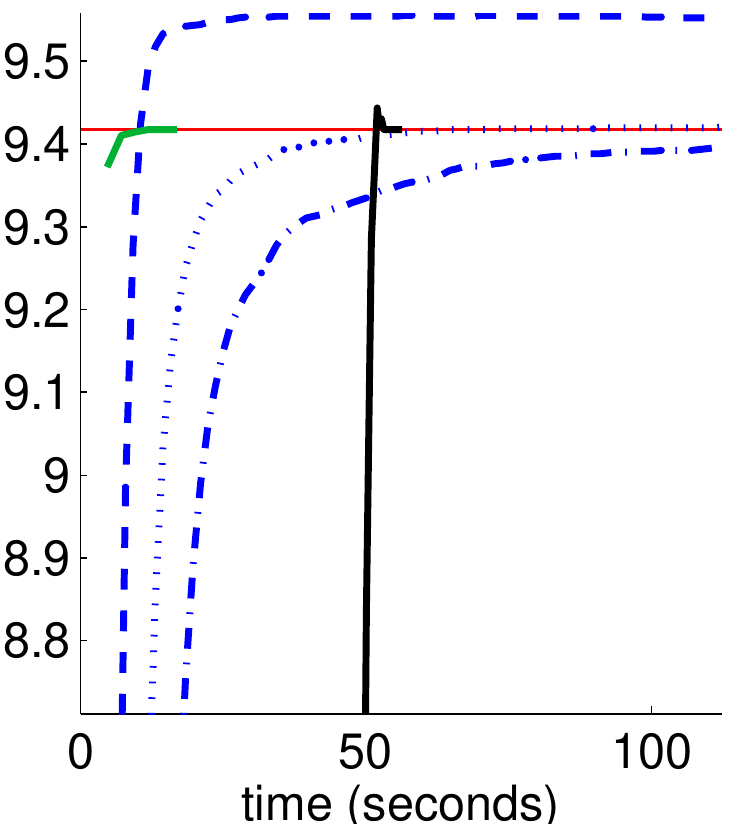}}
~~
\subfloat[$R=100$]{\includegraphics[scale=0.5]{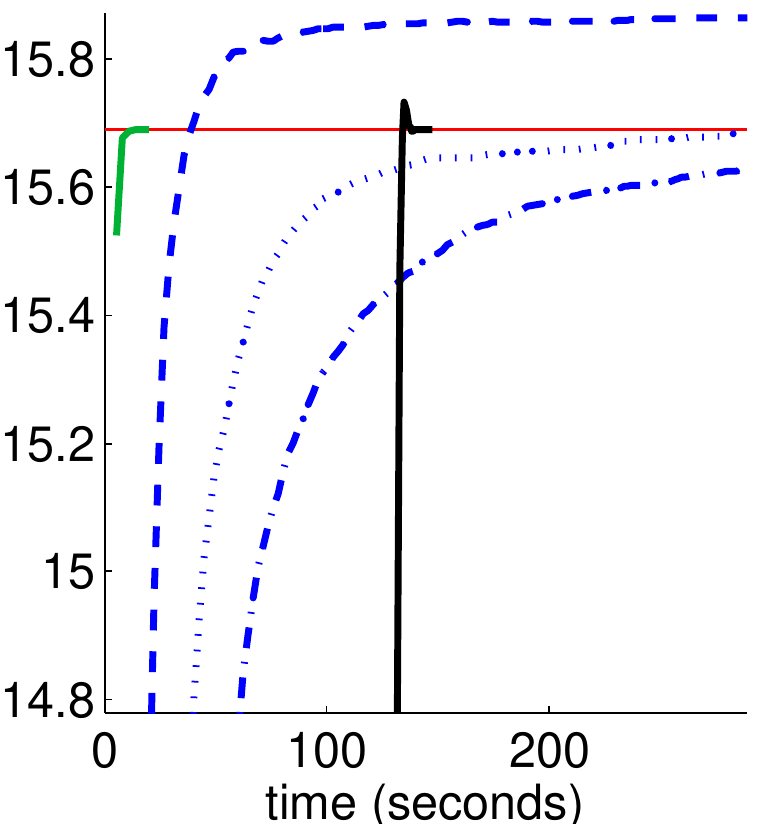}}
  \caption{Comparison of PDN-R, PQN-R and MU in finding elements of $\Mn{A}{1}$ equal to zero for a sample run.
    MU rarely finds exact zeros; therefore, we show results of applying various thresholds. Some experimentation may be needed to find the best threshold; regardless, it is slower than the methods proposed here.}
      \label{fig:subprob-MM-zeros}
\end{figure}

Figure~\ref{fig:subprob-compare-KKT} indicates that algorithm MU is preferred
if a relatively large $\kktviol$ is acceptable.  We contend that this is
not a good choice if the goal is to find a sparse solution.
Figure~\ref{fig:subprob-compare-nzero} plots the number of elements
that equal zero as a function of CPU time.  It shows that PDN-R and PQN-R
both converge to the correct number of zeros much faster than MU.

On closer inspection we see that MU is actually making factor elements small,
and is just very slow at making
them exactly zero.  If we choose a small positive threshold instead of zero,
then MU might arguably do well at finding a sparse solution.
Figure~\ref{fig:subprob-MM-zeros} summarizes an investigation of this idea.
Three different thresholds are shown: $10^{-3}$, $10^{-4}$, and $10^{-5}$.
The first threshold is clearly too large, declaring elements to be ``zero'' when
they never converge to such a value.  A threshold of $10^{-4}$ is also too
large for $R=20$, though possibly acceptable for $R=40$ and $R=60$.
The choice of $10^{-5}$
correctly identifies elements converging to zero, but PDN-R and PQN-R identifies
them much faster.
We conclude that PDN-R and PQN-R are significantly better at finding a true sparse
solution than MU, in terms of robustness (no need to choose an ad-hoc
threshold) and computation time (assuming a suitable threshold for MU is known).

%%%%%%%%%%%%%%%%%%%%%%%%%%%
\subsection{Solving the Full Problem}  \label{subsec:exp-fullprob}
%%%%%%%%%%%%%%%%%%%%%%%%%%%

In this section we move from a convex subproblem to solving the full
factorization (\ref{fullprob}).
We generate the same $200 \times 300 \times 400$ tensor data as in
Section~\ref{subsec:exp-subprob} and now treat all modes as optimization
variables.
An initial guess is constructed for all three modes in the same manner that
$\Mn{A}{1}$ was initialized in Section~\ref{subsec:exp-subprob}.
We generate ten different tensors by changing the random seed used in
Algorithm~\ref{alg:app-gen} and solve each from a single initial guess.
All tensors are factorized from the same initial guess;
however, since the full factorization
is a nonconvex optimization problem, algorithms may converge to different
local solutions.

We expect our local solutions to be reasonably close to the
multilinear model $\TM = \KTsmall{\Vl;
  \Mn{A}{1},\dots,\Mn{A}{N}}$
that generated the synthetic tensor data.
We compared computed solutions with the original model using the
Tensor Toolbox function {\tt score} with option {\tt greedy}.
This function implements the congruence test described in~\cite{TomasiBro2006}
and \cite{EvrimAcar2011a}.
The comparison considers angular differences between corresponding vectors
of the factor matrices, producing a number
between zero (poor match) and one (exact match).
Solutions computed with any algorithm to a tolerance of
$\tau = 10^{-4}$ scored above 0.84
(see the supplementary material for a detailed breakdown).
Perfect scores cannot be seen because the tensor data is generated from
the model as a noisy Poisson process.
Scores of less than 0.01 resulted when comparing the solution
to other models generated with a different random seed.
These results show that an accurate factorization can yield good approximations
to the original factors for our test problems; however, our focus is on
behavior of the algorithms in computing a solution.

%%%%%%%%%%%%%%%%%%%%%%%%%%%
\subsubsection{Comparing PDN-R with PDN}  \label{subsec:exp-PDNRfull}
%%%%%%%%%%%%%%%%%%%%%%%%%%%

We first compare the two Newton-based methods: PDN-R, which solves row
subproblems for each tensor mode, and PDN, which instead solves the block
subproblem as a single matrix.
In Section~\ref{subsec:exp-PDNRsub} we saw that the two methods behaved
similarly for the convex subproblem of a single tensor mode (except that PDN-R
was faster).  However, on the full factorization PDN is often unable to make
progress from a start point where the KKT violation is large.
Sometimes the search direction does not satisfy the sufficient decrease
condition of the Armijo line search, even after ten backtracking iterations.
More frequently, the line search puts too many variables at the bound of zero,
causing the objective function to become undefined in equation~(\ref{modeNprob})
because $\Mn{B}{n} \Mn{\Pi}{n}$ is zero for elements where $\Mz{X}{n}$ is nonzero.

If the line search fails in a subproblem, then we compute a multiplicative update
step for that iteration to make progress.  This allows PDN to reach points where
the KKT error is smaller, and we find that subsequent damped Newton steps are
successful
until convergence.  Table~\ref{tbl:linesearch-failures} quantifies the number of
line search failures over the first 20 iterations, beginning from a random
start point where $\kktviol$ is typically larger than $10^{3}$.
Columns in the table correspond to different values for the initial damping
parameter $\mu_0$.  We expect larger values of $\mu_k$ to improve robustness by
effectively shortening the step length and hopefully avoiding the mistake of
setting too many variables to zero.  However, a serious drawback to increasing
$\mu_k$ is that it damps out Hessian information, which can hinder the
convergence rate.
The table shows that improvement in robustness
is made; however, PDN still suffers from some line search failures.
In contrast, PDN-R does not have any line search failures for the same test
problems and start points, using the default $\mu_0 = 10^{-5}$.

\begin{table}[ht]
  \centering
  \small
  \caption{Line search failures by PDN in the first 20 iterations,
           averaged over five runs.  There were up to 900 possible
           line searches in each case (a maximum of 15 inner iterations per mode,
           over 20 outer iterations).}
      \label{tbl:linesearch-failures}
  \begin{tabular}{r|c|c|c}
    $R$ & $\mu_0 = 10^{1}$ & $\mu_0 = 10^{-2}$ & $\mu_0 = 10^{-5}$  \\
    \hline
     20 &  57.8 &  88.4 & 142.2  \\
     40 &  76.2 &  87.4 & 164.6  \\
     60 &  59.0 &  90.8 & 201.0  \\
     80 &  41.4 &  82.8 & 184.6  \\
    100 &  28.8 &  62.2 & 168.0
  \end{tabular}
\end{table}

Table~\ref{tbl:fullsynth-PDNR} shows that PDN-R is significantly faster than
PDN even in the region where PDN operates robustly.  These runs begin at a
start point where $\kktviol < 0.1$
and use $\mu_0 = 10^{-2}$
(PDN-R uses its default of $\mu_0 = 10^{-5}$)
so that PDN does not suffer any line search failures.
Five runs are made for each of the five values of $R$, and the method
stops when the algorithm reduces $\kktviol$ below a given
threshold (rows of Table~\ref{tbl:fullsynth-PDNR}).  PDN does not always reach
a threshold value in the three-hour-computation-time limit, but PDN-R always
succeeds.  The third column
shows that the number of outer iterations needed to reach a threshold is
very similar between PDN-R and PDN.  The fourth column shows that PDN-R
executes much faster.

\begin{table}[ht]
  \centering
  \small
  \caption{Comparison of PDN-R and PDN execution times for various stop
           tolerances.  25 experiments were run, the algorithms compared for
           each experiment that PDN completed, and the average value reported.  
           The third column is computed as
           $| \text{its}_{\text{PDN}} - \text{its}_{\text{PDN-R}} |
            / \max \{ \text{its}_{\text{PDN}}, \text{its}_{\text{PDN-R}} \}$.
           The fourth column shows that PDN-R executes from 8 to 9 times faster
           than PDN (the column reports average and standard deviation).
           The last column shows average execution time of PDN-R. }
      \label{tbl:fullsynth-PDNR}
  \begin{tabular}{c|c|c|c|c}
    $\kktviol$ & PDN failures & avg diff in outer its & PDN-R speedup
               & PDN-R time  \\
    \hline
    $10^{-2}$ &  2 &  2.48 \% &  $9.1 \pm 1.4$ & 463.3 secs  \\
    $10^{-3}$ &  5 &  3.02 \% &  $8.7 \pm 1.4$ & 595.0 secs  \\
    $10^{-4}$ &  9 &  2.68 \% &  $8.5 \pm 1.7$ & 609.7 secs  \\
    $10^{-5}$ & 13 &  7.93 \% &  $9.5 \pm 2.4$ & 573.1 secs
  \end{tabular}
\end{table}

Iterations of PDN-R run faster because each row subproblem has an individualized
step size and damping parameter (this was discussed previously in
Section~\ref{subsec:exp-PDNRsub}).
Given the large disparity in execution time and the lack of
robustness when far from a solution, we find no advantages to using PDN and
do not consider it further.

%%%%%%%%%%%%%%%%%%%%%%%%%%%
\subsubsection{Comparing PDN-R, PQN-R, and MU}
%%%%%%%%%%%%%%%%%%%%%%%%%%%

\begin{figure}[htbp]
\centering
\subfloat[$R=20$]{\includegraphics[scale=0.5]{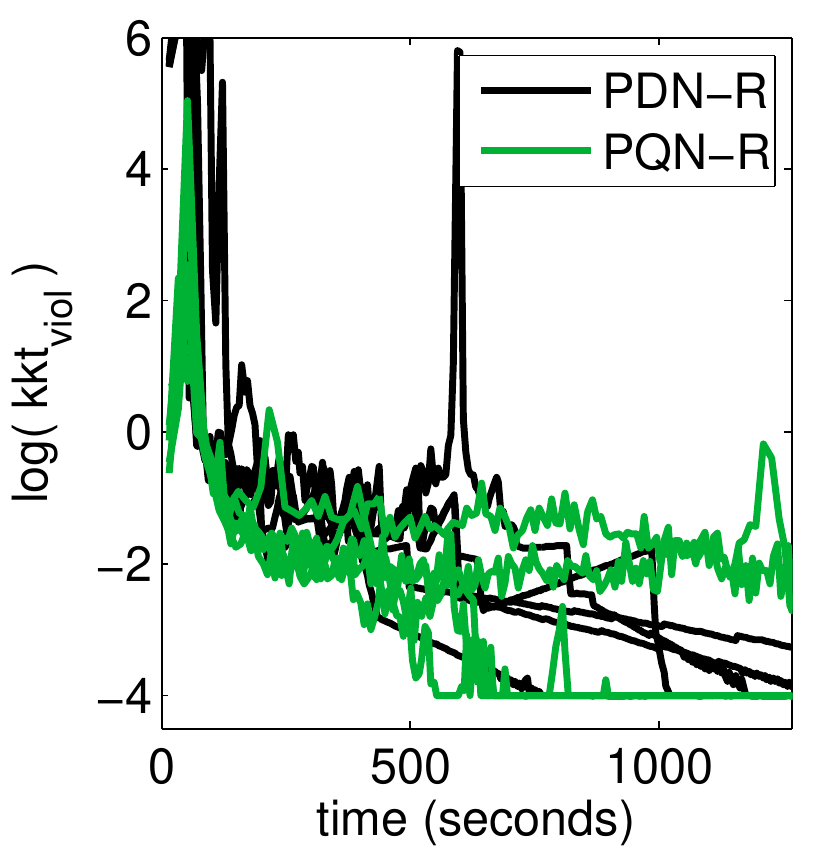}}
~~
\subfloat[$R=60$]{\includegraphics[scale=0.5]{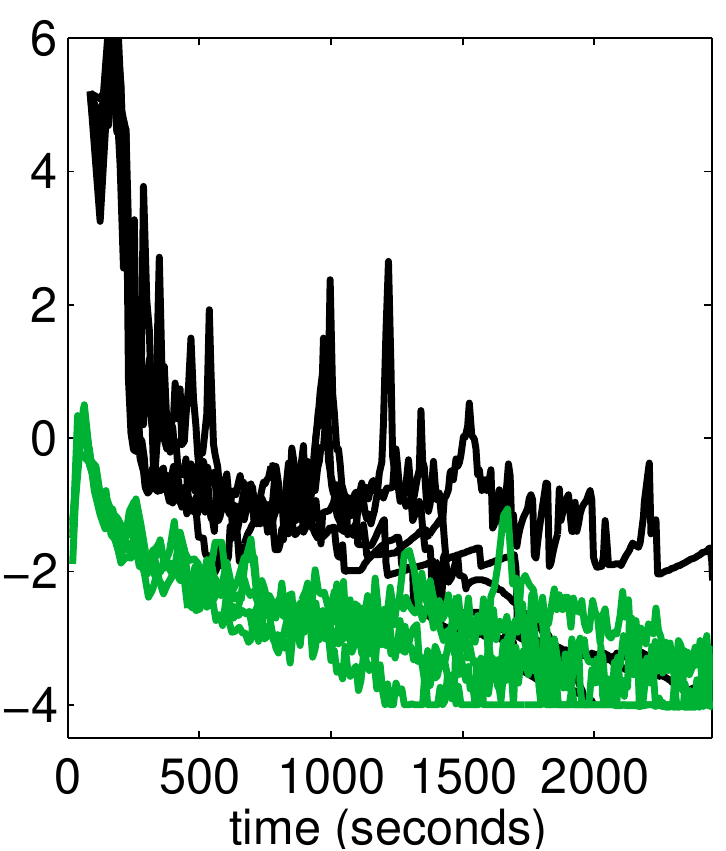}}
~~
\subfloat[$R=100$]{\includegraphics[scale=0.5]{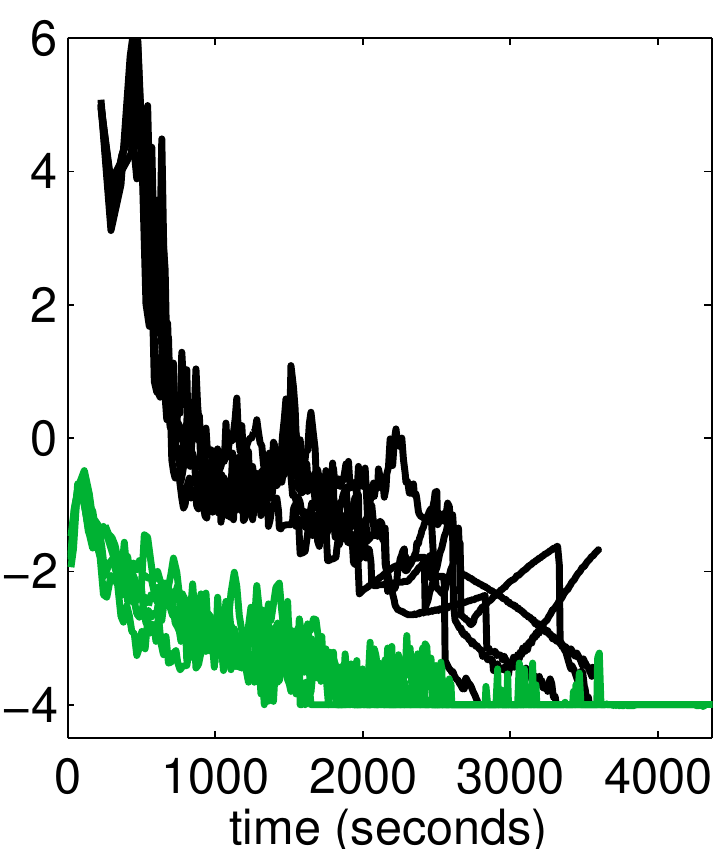}}

  \caption{Convergence behavior of the PDN-R (black lines) and PQN-R (green lines)
           algorithms in computing a full three-way solution.
           Each algorithm was run on ten different tensors.}
      \label{fig:fullprob-compare-kkt}
\end{figure}

Table~\ref{tbl:fullsynth-kkt} summarizes the time to reach a KKT
threshold of $10^{-3}$ for each algorithm over the synthetic data tensors.
Like the convex subproblem tested in Section~\ref{subsec:3algs-subprob},
the PDN-R and PQN-R methods converge to this relatively high accuracy
much faster than MU, again showing the value of second-order information.
As in the subproblem, we see that PQN-R is faster relative to PDN-R as
the number of components, $R$, increases.
Figure~\ref{fig:fullprob-compare-kkt} shows convergence behavior of the full
factorization problem in the same way that Figure~\ref{fig:subprob-compare-KKT}
showed behavior of the convex subproblem.
The KKT error of the full problem does not reach the quadratic rate of decrease
seen in the subproblem.  This is due to nonconvexity of the full factorization
problem, and the alternation between solutions of each mode.

\begin{table}[ht]
  \centering
  \small
  \caption{
           Time to reach stop tolerance $10^{-3}$ on full problem (over ten runs).
           Mean and standard deviation are reported.  Some runs of MU
           failed to reach the tolerance in three hours of execution.
           Results for different stop tolerances are in the
           supplementary material.
          }
      \label{tbl:fullsynth-kkt}
  \begin{tabular}{r|r|r|rr}
    $R$ & \multicolumn{1}{c|}{\textsc{PDN-R}}
        & \multicolumn{1}{c|}{\textsc{PQN-R}}
        & \multicolumn{2}{c}{\textsc{MU}}  \\
    \hline
    20 & \, 229 $\pm$ \, 57 secs
       & \, 397 $\pm$   123 secs
       &   3355 $\pm$   1933 secs  & (0 failures)  \\
    40 & \, 493 $\pm$   151 secs
       & \, 818 $\pm$   185 secs
       &   8101 $\pm$   2045 secs  & (2 failures)  \\
    60 &   1003 $\pm$   349 secs
       & \, 966 $\pm$   286 secs
       &   9628 $\pm$ \, 978 secs  & (5 failures)  \\
    80 &   1682 $\pm$   642 secs
       &   1639 $\pm$   390 secs
       &  no successes             & (10 failures)  \\
   100 &   2707 $\pm$   773 secs
       &   1995 $\pm$   743 secs
       &  no successes             & (10 failures)
  \end{tabular}
\end{table}

\begin{figure}[htbp]
\centering
\subfloat[$R=20$]{\includegraphics[scale=0.5]{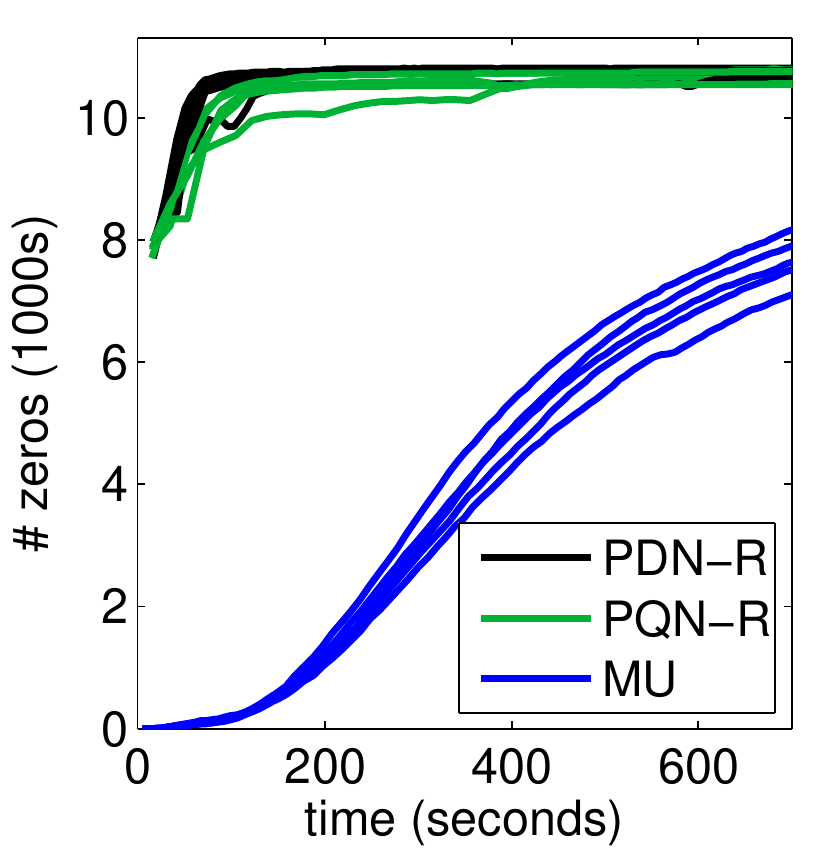}}
~~
\subfloat[$R=60$]{\includegraphics[scale=0.5]{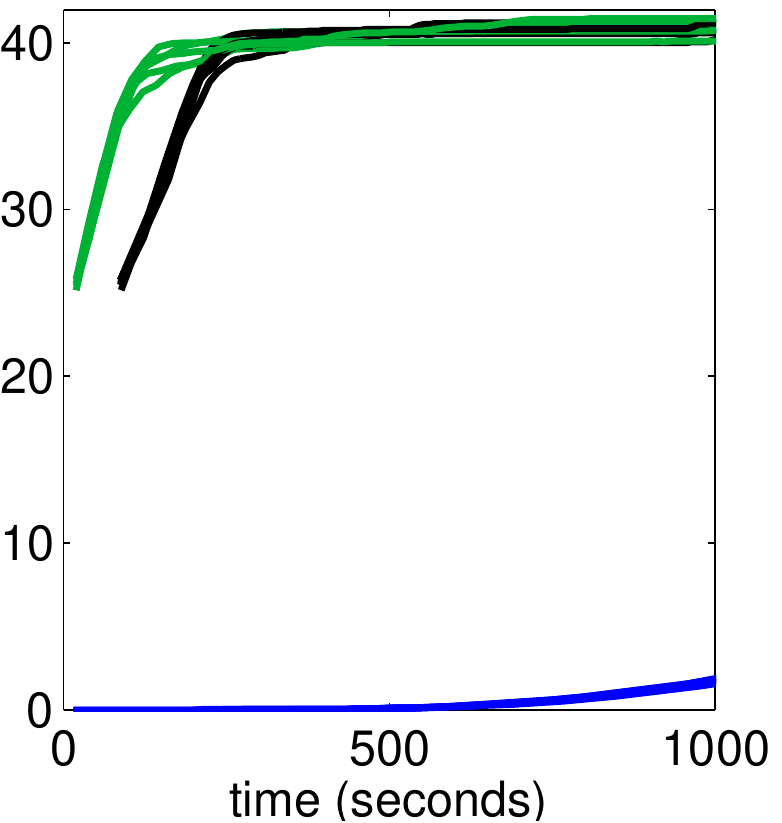}}
~~
\subfloat[$R=100$]{\includegraphics[scale=0.5]{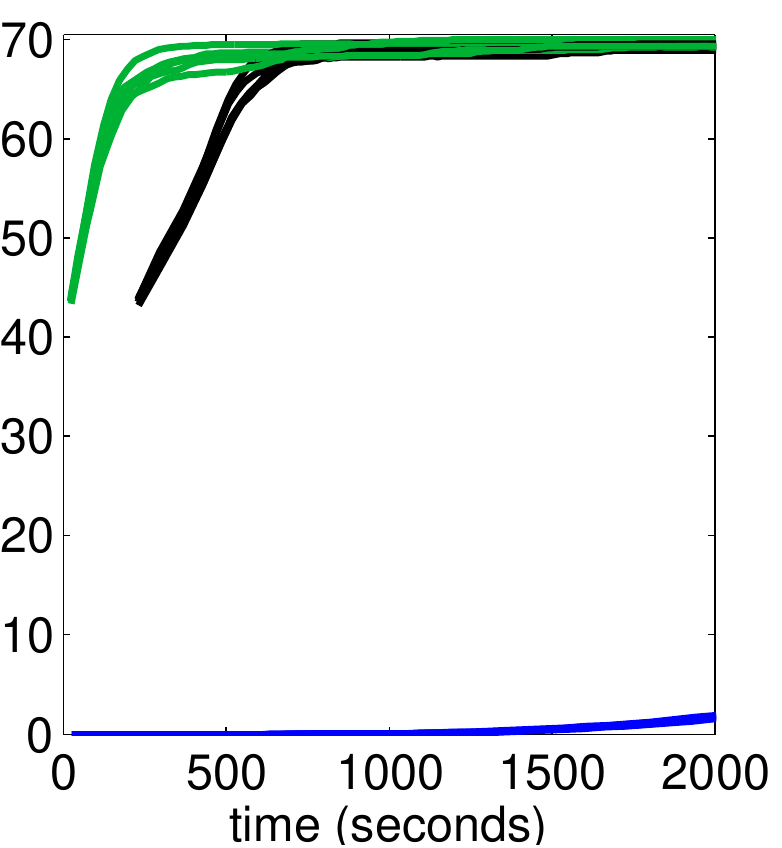}}
  \caption{Effectiveness of the algorithms in finding a sparse solution for
           a full three-way solution.  In each case the total number of elements
           in $\Mn{A}{1}$, $\Mn{A}{2}$, and $\Mn{A}{3}$
           equal to zero is plotted against execution time.
           The PDN-R (black lines) and PQN-R (green) algorithms are
           much faster than MU (blue).
           Each algorithm was run on ten different tensors,
           so the final number of zero elements has ten different values.} 
      \label{fig:fullprob-compare-zeros}
\end{figure}

As with the subproblem, we also observe better convergence by our methods
to a sparse solution.
Figure~\ref{fig:fullprob-compare-zeros} shows PDN-R and PQN-R reaching
the final count of zero elements much faster than MU.
As in Section~\ref{subsec:3algs-subprob}, we argue that PDN-R and PQN-R
are superior when the task is to find a solution with correct sparsity.

We performed similar experiments on tensors of the same size but different
sparsity.  Results are in the supplementary material.
They lead to the same conclusions as the data in Table~\ref{tbl:fullsynth-kkt};
namely, that PDN-R and PQN-R are faster than MU, and PQN-R becomes faster than
PDN-R as the number of components increases.

The supplementary material also describes a simple experiment with
sparse tensors whose factor
matrices have a high degree of collinearity between column vectors.
Such problems sometimes lead to poor algorithm performance (e.g., the ``swamps''
in \cite{Paatero1997}).  Performance of PDN-R and PQN-R was much better
than algorithm MU in this experiment as well.

%%%%%%%%%%%%%%%%%%%%%%%%%%%
\subsubsection{Comparing with DBLP Data}
%%%%%%%%%%%%%%%%%%%%%%%%%%%

\begin{figure}[htbp]
\centering
\subfloat[$R=20$]{\includegraphics[scale=0.5]{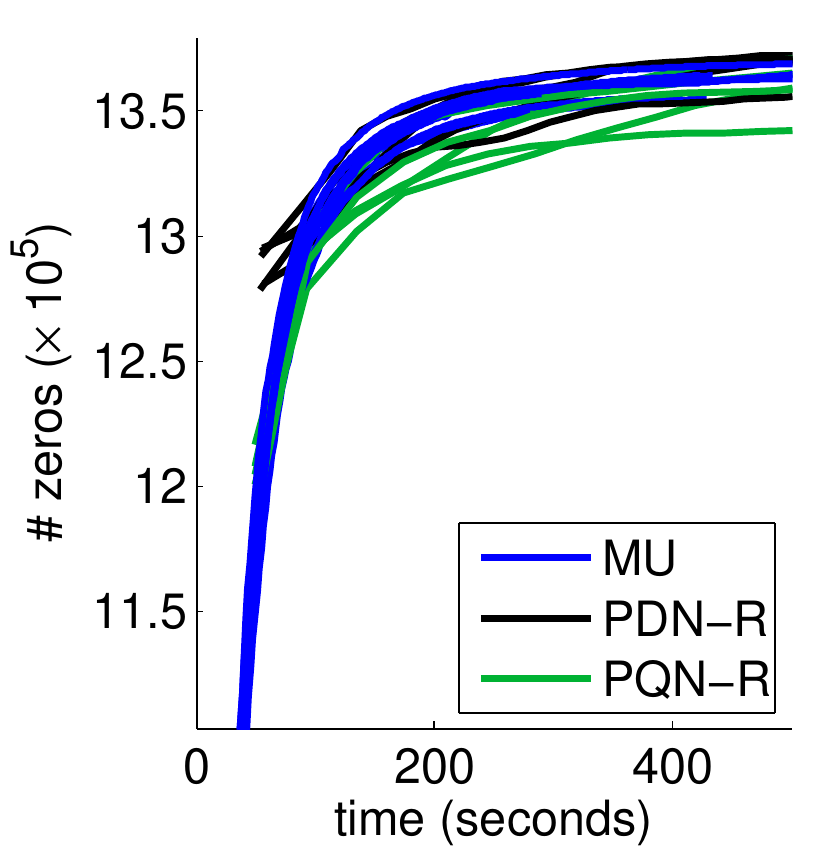}}
~~
\subfloat[$R=60$]{\includegraphics[scale=0.5]{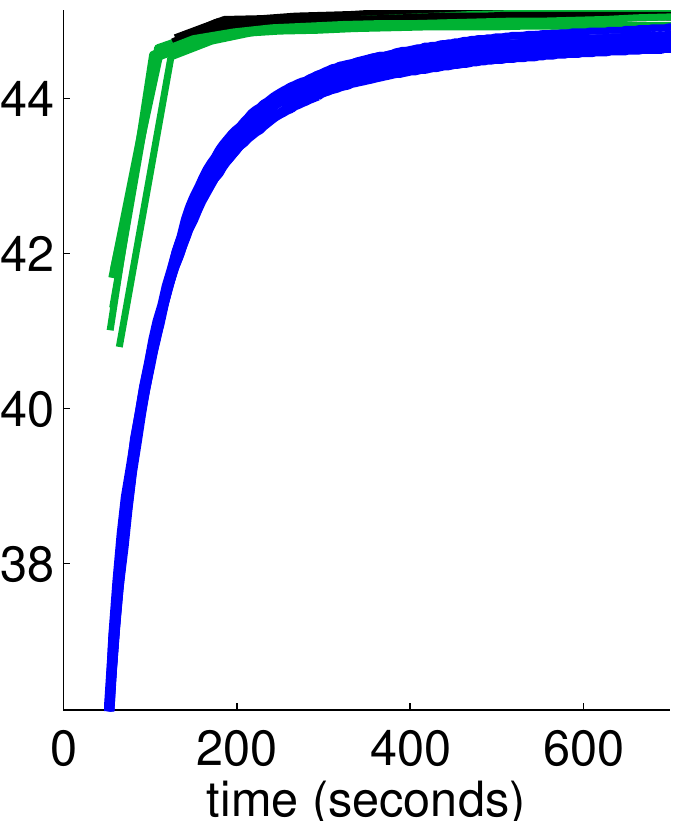}}
~~
\subfloat[$R=100$]{\includegraphics[scale=0.5]{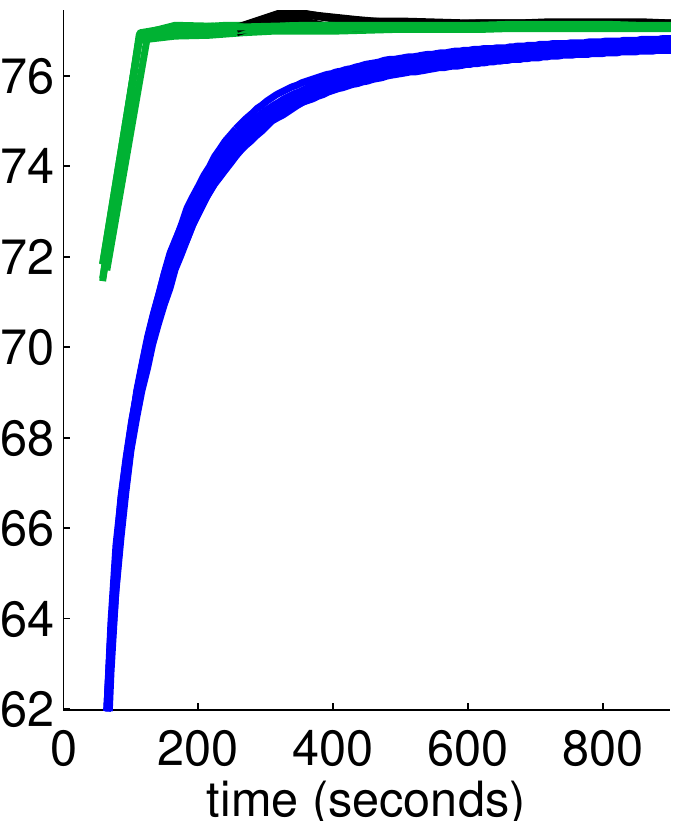}}

  \caption{Effectiveness of the algorithms in finding a sparse solution for
           the DBLP tensor.  In each case the total number of elements
           in $\Mn{A}{1}$, $\Mn{A}{2}$, and $\Mn{A}{3}$
           equal to zero is plotted against execution time.
           The PDN-R (black lines) and PQN-R (green) algorithms are
           much faster than MU (blue).}
      \label{fig:dblp-compare-zeros}
\end{figure}

We also compare the same three algorithms on the sparse three-way tensor of
DBLP data \cite{DBLP:data} examined in \cite{Kolda:DBLP}.
The data counts the number of papers published by author $i_1$ at conference
$i_2$ in year $i_3$, with dimensions $7108 \times 1103 \times 10$.
The tensor contains 112,730 nonzero elements, a density of 0.14\%.
The data was factorized for $R$ between 10 and 100 in \cite{Kolda:DBLP}
(using a least squares objective function), so we use $R \in \{ 20, 60, 100 \}$
in our experiments.
Behavior of the algorithms on the DBLP data was similar to behavior on our
synthetic data.  Figure~\ref{fig:dblp-compare-zeros} shows how the
count of elements equal to zero changes as algorithms progress, for ten runs
that start from different random initial guesses.
Again we see that PDN-R and PQN-R reach a sparse solution faster than MU.

\begin{figure}[htbp]
\centering
\subfloat{\includegraphics[scale=0.70]{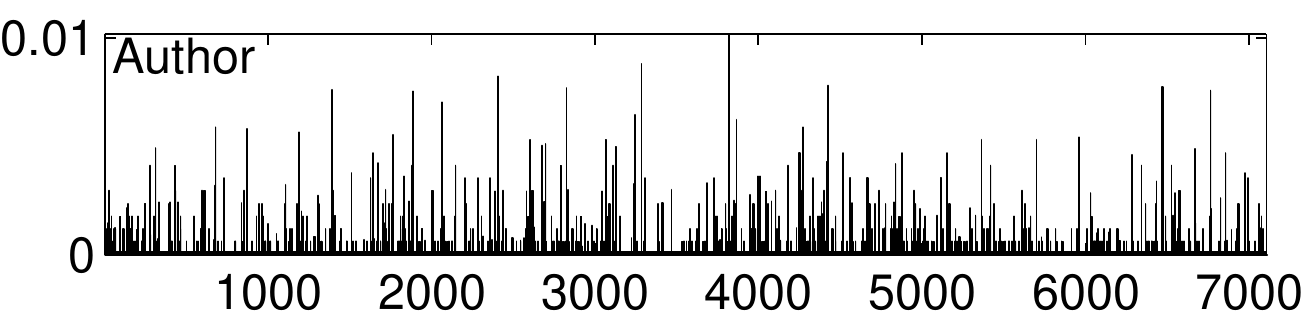}}
\\
\subfloat{\includegraphics[scale=0.70]{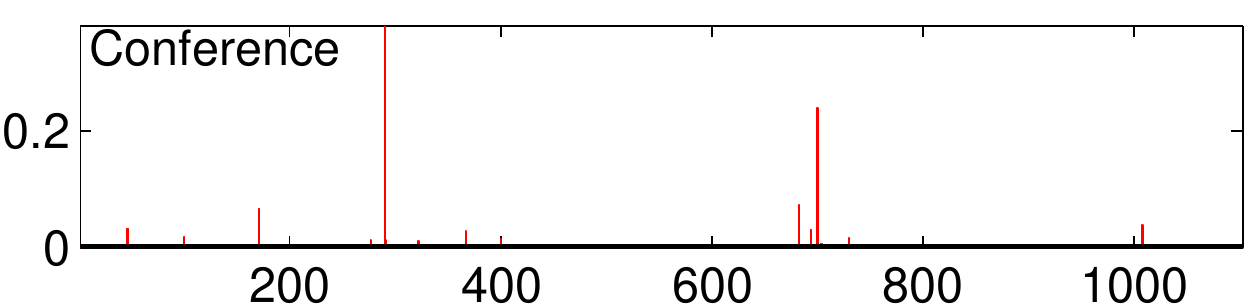}}
\\
\subfloat{\includegraphics[scale=0.70]{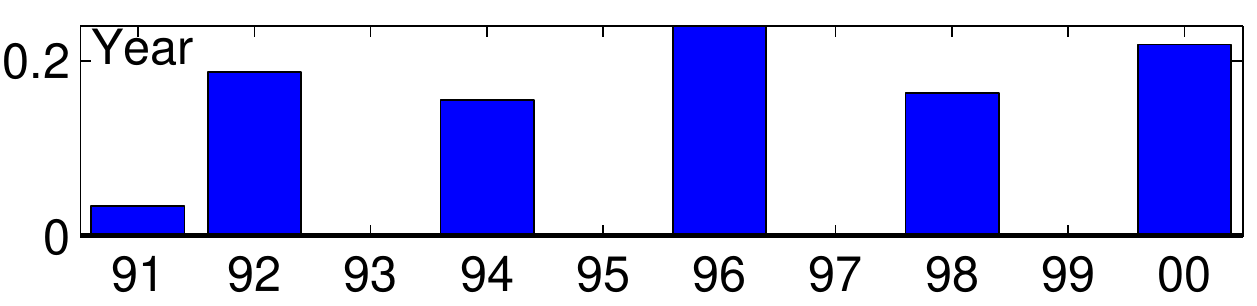}}
  \caption{Computed factors from DBLP data for component 26 (i.e., the 26th
           largest component by weight).  The two dominant conferences,
           ECAI and KR, occurred only in even years, except for KR in 1991.
           Factors are extremely sparse: 91\% (6456) of elements in the author
           factor are exactly zero,
           as are 98\% (1084) of the conference elements.}
      \label{fig:dblp-factors-comp26}
\end{figure}

\begin{figure}[htbp]
\centering
\subfloat{\includegraphics[scale=0.70]{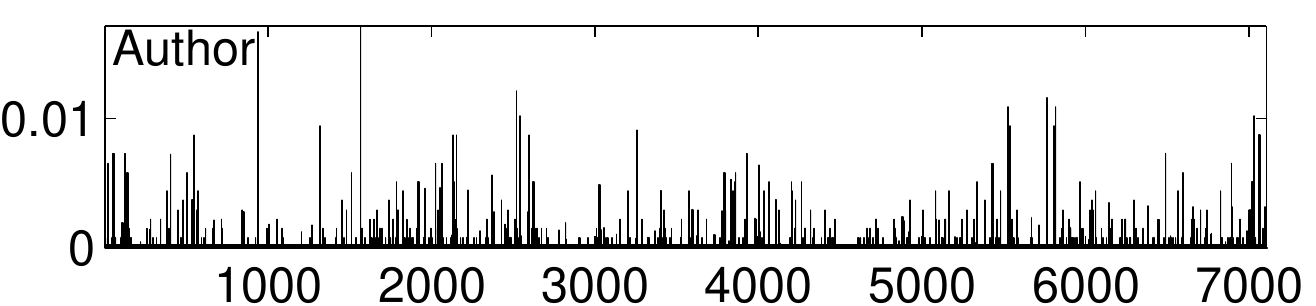}}
\\
\subfloat{\includegraphics[scale=0.70]{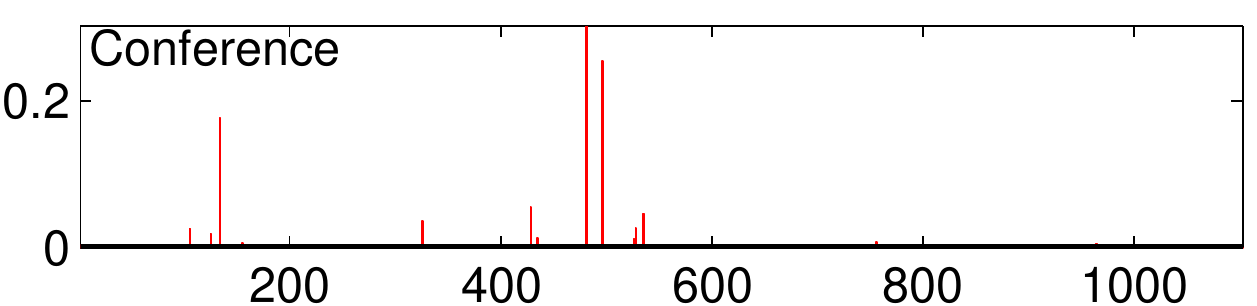}}
\\
\subfloat{\includegraphics[scale=0.70]{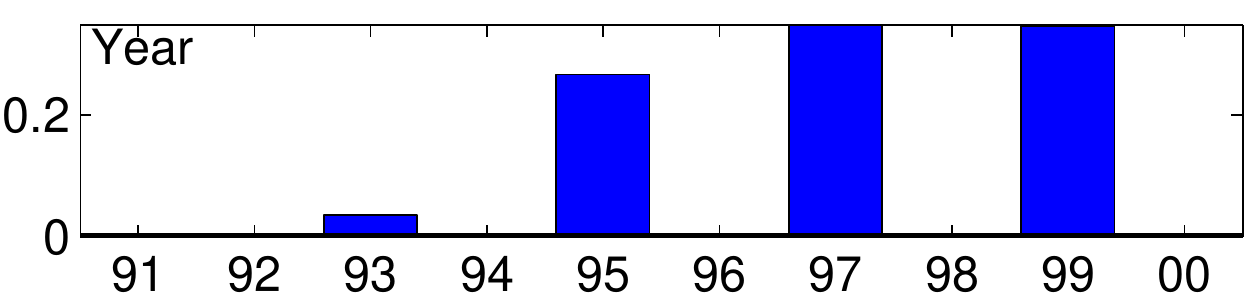}}
  \caption{Computed factors from DBLP data for component 41.
           The three dominant conferences (ICDAR, ICIAP, and CAIP)
           occurred only in odd years.
           In this component,
           93\% (6626) of elements in the author factor are exactly zero,
           as are 98\% (1083) conference elements.}
      \label{fig:dblp-factors-comp41}
\end{figure}

Factorizations of the DBLP data computed with PDN-R and PQN-R were
quite sparse, making them easier to interpret.
The fraction of elements exactly equal to zero in the computed conference
factor matrix was 98.1\%.  The author factor matrix was also very sparse,
with 95.4\% of the elements exactly zero.  These results were for a factorization
with $R=100$, stopped after 800 seconds of execution with the KKT violation
reduced to around $5 \times 10^{-4}$.
Figure~\ref{fig:dblp-factors-comp26} shows a component that detects related
conferences that took place only in even years.  The two dominant conferences
are the same as those reported in Figure 7 of \cite{Kolda:DBLP}.
Figure~\ref{fig:dblp-factors-comp41} shows another component that groups
conferences that took place only in odd years.
In both components, the sparsity is striking, especially for the conference
factor.

\section{Summary}
In this paper we consider the problem of nonnegative tensor factorization
using a \text{K-L} objective function, and we derive a \emph{row subproblem}
formulation that allows efficient use of second order information.
We present two new algorithms
that exploit the row subproblem reformulation:  PDN-R uses second derivatives
in the optimization, while PQN-R uses a quasi-Newton approximation.
We show that using the same second order information in a block subproblem
formulation is less robust and more expensive computationally than a
row subproblem formulation.
We show that both PDN-R and PQN-R are much faster than the best multiplicative
update method,especially when high accuracy solutions are desired.
We further show
that high accuracy is needed to identify zeros and compute sparse factors
without resorting to the use of ad-hoc thresholds.  This is important because
sparse count data is likely
to have sparsity in the factors, and sparse factors are always easier to
interpret.

Our Matlab algorithms will appear in Version 2.6 of the
Tensor Toolbox \cite{BaderKolda:MatlabTTB}.
We mention in section~\ref{subsec:ptf-reformulation} that row subproblems
can be solved in parallel, and we anticipate developing other versions of
the algorithms for shared and distributed memory machines.

\section*{Acknowledgments}

We thank the authors of \cite{KimSuvritDhillon:BoxConstraints} for sharing
Matlab code that we used in the experiments.
We thank the anonymous referees for clarifying some of our arguments,
and for suggesting additional experiments.

This work was partially funded by the Laboratory Directed Research
\& Development (LDRD) program at Sandia National Laboratories.
Sandia National Laboratories is a multiprogram laboratory operated by Sandia
Corporation, a wholly owned subsidiary of Lockheed Martin Corporation, for
the United States Department of Energy's National Nuclear Security
Administration under contract DE-AC04-94AL85000.

\bibliographystyle{siam}

\begin{thebibliography}{10}

\bibitem{EvrimAcar2011a}
{\sc E.~Acar, D.~M. Dunlavy, and T.~G. Kolda}, {\em A scalable optimization
  approach for fitting canonical tensor decompositions}, Journal of
  Chemometrics, 25 (2011), pp.~67--86.

\bibitem{andrew2007scalable}
{\sc G.~Andrew and J.~Gao}, {\em Scalable training of $l_1$-regularized
  log-linear models}, in Proceedings of the 24th International Conference on
  Machine Learning, ACM, 2007, pp.~33--40.

\bibitem{Bader:Enron2007}
{\sc B.~W. Bader, M.~W. Berry, and M.~Browne}, {\em Discussion tracking in
  {E}nron email using {PARAFAC}}, in Survey of Text Mining: Clustering,
  Classification, and Retrieval, M.~W. Berry and M.~Castellanos, eds.,
  Springer, 2nd~ed., 2007, ch.~8, pp.~147--162.

\bibitem{BaderKolda:MatlabTTB}
{\sc B.~W. Bader, T.~G. Kolda, et~al.}, {\em {MATLAB} {T}ensor {T}oolbox
  version 2.5}.
\newblock Available online, January 2012.
\newblock \url{http://www.sandia.gov/~tgkolda/TensorToolbox/}.

\bibitem{bertsekas1976goldstein}
{\sc D.~Bertsekas}, {\em On the {G}oldstein-{L}evitin-{P}olyak gradient
  projection method}, IEEE Transactions on Automatic Control, 21 (1976),
  pp.~174--184.

\bibitem{bertsekas1982projected}
{\sc D.~Bertsekas}, {\em Projected newton methods for optimization problems
  with simple constraints}, SIAM Journal on Control and Optimization, 20
  (1982), pp.~221--246.

\bibitem{BroPhD}
{\sc R.~Bro}, {\em Multi-way Analysis in the Food Industry: Models, Algorithms,
  and Applications}, PhD thesis, Universiteit van Amsterdam, 1998.

\bibitem{BroDeJong1997}
{\sc R.~Bro and S.~D. Jong}, {\em A fast non-negativity-constrained least
  squares algorithm}, Journal of Chemometrics, 11 (1997), pp.~393--401.

\bibitem{byrd1995limited}
{\sc R.~H. Byrd, P.~Lu, J.~Nocedal, and C.~Zhu}, {\em A limited memory
  algorithm for bound constrained optimization}, SIAM Journal on Scientific
  Computing, 16 (1995), pp.~1190--1208.

\bibitem{CarrollChang:1970}
{\sc J.~D. Carroll and J.~J. Chang}, {\em Analysis of individual differences in
  multidimensional scaling via an {N}-way generalization of `{E}ckart-{Y}oung'
  decomposition}, Psychometrika, 35 (1970), pp.~283--319.

\bibitem{ChiKolda:cpapr}
{\sc E.~C. Chi and T.~G. Kolda}, {\em On tensors, sparsity, and nonnegative
  factorizations}, {SIAM} Journal on Matrix Analysis and Applications, 33
  (2012), pp.~1272--1299.

\bibitem{cichocki2009fast}
{\sc A.~Cichocki and A.-H. Phan}, {\em Fast local algorithms for large scale
  nonnegative matrix and tensor factorizations}, IEICE Transactions on
  Fundamentals of Electronics, Communications and Computer Sciences, 92 (2009),
  pp.~708--721.

\bibitem{conn1988global}
{\sc A.~Conn, N.~Gould, and P.~Toint}, {\em Global convergence of a class of
  trust region algorithms for optimization with simple bounds}, SIAM Journal on
  Numerical Analysis, 25 (1988), pp.~433--460.

\bibitem{DBLP:data}
{\em {DBLP} data}, 2011.
\newblock \url{http://www.informatik.uni-trier.de/~ley/db/}.

\bibitem{Kolda:DBLP}
{\sc D.~M. Dunlavy, T.~G. Kolda, and E.~Acar}, {\em Temporal link prediction
  using tensor and matrix factorizations}, {ACM} Transactions on Knowledge
  Discovery from Data, 5 (2011).

\bibitem{FriedlanderHatz2008}
{\sc M.~P. Friedlander and K.~Hatz}, {\em Computing nonnegative tensor
  factorizations}, Computational Optimization and Applications, 23 (2008),
  pp.~631--647.

\bibitem{gonzalez2005accelerating}
{\sc E.~Gonzalez and Y.~Zhang}, {\em Accelerating the {L}ee-{S}eung algorithm
  for non-negative matrix factorization}, Tech. Report TR-05-02, Department of
  Computational and Applied Mathematics, Rice University, Houston, TX, 2005.

\bibitem{Harshman:1970}
{\sc R.~A. Harshman}, {\em Foundations of the {PARAFAC} procedure: Models and
  conditions for an ``explanatory'' multi-modal factor analysis}, {UCLA}
  Working Papers in Phonetics, 16 (1970).
\newblock Available online, \url{http://publish.uwo.ca/~harshman/wpppfac0.pdf}.

\bibitem{DhillonOneVar2011}
{\sc C.-J. Hsieh and I.~S. Dhillon}, {\em Fast coordinate descent methods with
  variable selection for non-negative matrix factorization}, in Proceedings of
  the 17th ACM SIGKDD International Conference on Knowledge Discovery and Data
  Mining, ACM, 2011, pp.~1064--1072.

\bibitem{KimSuvritDhillon:BoxConstraints}
{\sc D.~Kim, S.~Sra, and I.~Dhillon}, {\em Tackling box-constrained
  optimization via a new projected quasi-{N}ewton approach}, SIAM Journal on
  Scientific Computing, 32 (2010), pp.~3548--3563.

\bibitem{KimSraDhillon2008}
{\sc D.~Kim, S.~Sra, and I.~S. Dhillon}, {\em Fast projection-based methods for
  the least squares non-negative matrix approximation problem}, Statistical
  Analysis and Data Mining, 1 (2008), pp.~38--51.

\bibitem{KimPark2008}
{\sc H.~Kim and H.~Park}, {\em Nonnegative matrix factorization based on
  alternating nonnegativity constrained least squares and active set method},
  SIAM Journal on Matrix Analysis and Applications, 30 (2008), pp.~713--730.

\bibitem{kim2012fast}
{\sc J.~Kim and H.~Park}, {\em Fast nonnegative tensor factorization with an
  active-set-like method}, in High-Performance Scientific Computing, Algorithms
  and Applications, M.~W. Berry, K.~A. Gallivan, E.~Gallopoulos, A.~Grama,
  B.~Philippe, Y.~Saad, and F.~Saied, eds., Springer, 2012, pp.~311--326.

\bibitem{KoldaBader:SiamReview}
{\sc T.~G. Kolda and B.~W. Bader}, {\em Tensor decompositions and
  applications}, SIAM Review, 51 (2009), pp.~455--500.

\bibitem{seungLeeNature1999}
{\sc D.~D. Lee and H.~S. Seung}, {\em Learning the parts of objects by
  non-negative matrix factorization}, Nature, 401 (1999), pp.~788--791.

\bibitem{seung2001algorithms}
\leavevmode\vrule height 2pt depth -1.6pt width 23pt, {\em Algorithms for
  non-negative matrix factorization}, Advances in Neural Information Processing
  Systems, 13 (2001), pp.~556--562.

\bibitem{lin2007projected}
{\sc C.~Lin}, {\em Projected gradient methods for nonnegative matrix
  factorization}, Neural computation, 19 (2007), pp.~2756--2779.

\bibitem{liu2012sparse}
{\sc J.~Liu, J.~Liu, P.~Wonka, and J.~Ye}, {\em Sparse non-negative tensor
  factorization using column{\-}wise coordinate descent}, Pattern Recognition,
  45 (2012), pp.~649--656.

\bibitem{nocedalLBFGS}
{\sc J.~Nocedal}, {\em Updating quasi-{N}ewton matrices with limited storage},
  Mathematics of Computation, 35 (1980), pp.~773--782.

\bibitem{nocedal1999numerical}
{\sc J.~Nocedal and S.~J. Wright}, {\em Numerical Optimization}, Springer,
  2nd~ed., 2006.

\bibitem{Paatero1997}
{\sc P.~Paatero}, {\em A weighted non-negative least squares algorithm for
  three-way ``{PARAFAC}'' factor analysis}, Chemometrics and Intelligent
  Laboratory Systems, 38 (1997), pp.~223--242.

\bibitem{PaateroTapper1994}
{\sc P.~Paatero and U.~Tapper}, {\em Positive matrix factorization: A
  non-negative factor model with optimal utilization of error estimates of data
  values}, Environmetrics, 5 (1994), pp.~111--126.

\bibitem{PhanZdunek2010}
{\sc A.~Phan, A.~Cichocki, R.~Zdunek, and T.~Dinh}, {\em Novel alternating
  least squares algorithm for nonnegative matrix and tensor factorizations}, in
  Neural Information Processing. Theory and Algorithms, K.~W. Wong, B.~S.~U.
  Mendis, and A.~Bouzerdoum, eds., vol.~6443 of Lecture Notes in Computer
  Science, Springer Berlin Heidelberg, 2010, pp.~262--269.

\bibitem{PhanCichocki2011}
{\sc A.~H. Phan, P.~Tichavsky, and A.~Cichocki}, {\em Fast damped
  {G}auss-{N}ewton algorithm for sparse and nonnegative tensor factorization},
  in Acoustics, Speech and Signal Processing (ICASSP), 2011 IEEE International
  Conference on, IEEE, 2011, pp.~1988--1991.

\bibitem{schmidt2012}
{\sc M.~Schmidt, D.~Kim, and S.~Sra}, {\em Projected {N}ewton-type methods in
  machine learning}, in Optimization for Machine Learning, S.~Sra, S.~Nowozin,
  and S.~J. Wright, eds., {MIT} Press, 2011, pp.~305--330.

\bibitem{Sun:DynamicTensorCyber}
{\sc J.~Sun, D.~Tao, and C.~Faloutsos}, {\em Beyond streams and graphs: Dynamic
  tensor analysis}, in KDD '06, Proceedings of the 12th ACM SIGKDD
  International Conference on Knowledge Discovery and Data Mining, ACM, 2006,
  pp.~374--383.

\bibitem{TomasiBro2006}
{\sc G.~Tomasi and R.~Bro}, {\em A comparison of algorithms for fitting the
  {PARAFAC} model}, Computational Statistics and Data Analysis, 50 (2006),
  pp.~1700--1734.

\bibitem{vavasis2009complexity}
{\sc S.~Vavasis}, {\em On the complexity of nonnegative matrix factorization},
  SIAM Journal on Optimization, 20 (2009), pp.~1364--1377.

\bibitem{WellingWeber2001}
{\sc M.~Welling and M.~Weber}, {\em Positive tensor factorization}, Pattern
  Recognition Letters, 22 (2001), pp.~1255--1261.

\bibitem{xu2012block}
{\sc Y.~Xu and W.~Yin}, {\em A block coordinate descent method for multi-convex
  optimization with applications to nonnegative tensor factorization and
  completion}, Tech. Report 12-15, CAAM Rice University, Houston, TX, 2012.

\bibitem{zafeiriou2011nonnegative}
{\sc S.~Zafeiriou and M.~Petrou}, {\em Nonnegative tensor factorization as an
  alternative {C}siszar-{T}usnady procedure: algorithms, convergence,
  probabilistic interpretations and novel probabilistic tensor latent variable
  analysis algorithms}, Data Mining and Knowledge Discovery, 22 (2011),
  pp.~419--460.

\bibitem{zdunek2006non}
{\sc R.~Zdunek and A.~Cichocki}, {\em Non-negative matrix factorization with
  quasi-{N}ewton optimization}, in Eighth International Conference on
  Artificial Intelligence and Soft Computing, ICAISC, Springer, 2006,
  pp.~870--879.

\bibitem{zdunek2007nonnegative}
\leavevmode\vrule height 2pt depth -1.6pt width 23pt, {\em Nonnegative matrix
  factorization with constrained second-order optimization}, Signal Processing,
  87 (2007), pp.~1904--1916.

\bibitem{ZhengZhang}
{\sc Y.~Zheng and Q.~Zhang}, {\em Damped {N}ewton based iterative non-negative
  matrix factorization for intelligent wood defects detection}, Journal of
  Software, 5 (2010), pp.~899--906.

\end{thebibliography}

\appendix
\section{Generating Synthetic Test Data}  \label{app-generate}

The goal is to create artificial nonnegative factor matrices and use these to
compute a data tensor whose elements follow a Poisson distribution with
multilinear parameters.
Factorizing the data tensor should yield quantities that are close to the
original factor matrices.
The procedure is based on the work of \cite{ChiKolda:cpapr}.

The data tensor should be sparse, reflecting Poisson distributions whose
probability of zero is not negligible.  Each entry is a count of the number
of samples assigned to this cell, out of a given total number of samples $S$.
We generate factor matrices in each tensor mode and treat them as stochastic
quantities to draw the $S$ samples that provide data tensor counts.

Our generation procedure utilizes the function {\tt create\_problem} from
Tensor Toolbox for Matlab \cite{BaderKolda:MatlabTTB},
supplying a custom function for the {\tt Factor\_Generator} parameter
(available as Matlab code from the authors).
We create a multilinear model, $\TM = [\V{\lambda}; \Mn{A}{1} \ldots \Mn{A}{N}]$,
where $\Mn{A}{n} \in {\cal R}^{I_n \times R}$
and $\V{\lambda} \in {\cal R}^R$, and sizes $I_n$ and $R$ are given.
The model is generated by the following procedure:

\begin{algorithm}
  \caption{Generation of Synthetic Sparse Poisson Tensor Data}
  \label{alg:app-gen}
  Given tensor size, $I_1 \times \cdots \times I_N$, number of components, $R$,
  and number of samples, $S$.
  \\
  Return a model $\TM = [\V{\lambda}; \Mn{A}{1} \ldots \Mn{A}{N}]$
  and corresponding sparse data tensor $\TX$.

  \begin{algorithmic}[1]
    \State In each column of $\Mn{A}{n}$, choose $20\%$ of the elements at random
           and set their value to $1 + 10 R x$, where $x$ is a random value
           from a uniform distribution on $(0,1)$.
           Set the other elements equal to the small constant 0.1.
           \label{appgen-alg:boost}
    \State Choose random values for elements in $\V{\lambda}$ from a uniform
           distribution on $(0,1)$.
    \State Rescale each column of $\Mn{A}{n}$ so entries sum to 1, absorbing the
           scale factor into the corresponding element of $\V{\lambda}$.
    \State Rescale the vector $\V{\lambda}$ so entries sum to 1.
           \label{appgen-alg:normlambda}

    \For{$s=1, \dots, S$}
      \State Treating $\V{\lambda}$ as a distribution, choose a component $r$
             at random.
      \State Treating the $r$th column of $\Mn{A}{1}$ as a distribution,
             choose an index $i_1$ with probability proportional to
             $\V{a}^{(1)}_r$.  Do the same for indices $i_2, \ldots, i_N$,
             resulting in the index $\ib$ chosen with probability
      \begin{equation*}
        P(\ib) = a^{(1)}_{i_1r} a^{(2)}_{i_2r} \ldots a^{(N)}_{i_Nr} .
      \end{equation*}
      \State Increment the $\ib$th entry of $\TX$ by one.
    \EndFor

    \State Rescale $\V{\lambda} \gets S \V{\lambda}$ so that
           $\| \V{\lambda} \|_1 = S$.
           \Comment{Recall Step~\ref{appgen-alg:normlambda}
                    sets $\| \V{\lambda} \|_1 = 1$.}
           \label{appgen-alg:rescale-lambda}

  \end{algorithmic}
\end{algorithm}

Step~\ref{appgen-alg:boost}
defines a strong preference for certain values of each index.
As $R$ increases, the relative probability of these indices is also increased
so that they continue to stand out as strong preferences.

Step~\ref{appgen-alg:rescale-lambda} rescales $\V{\lambda}$ so that the
$\ell_1$ norms of the generated data tensor $\TX$ and K-tensor are the
same in any mode-$n$ unfolding.  For example:

\begin{align*}
  \left\| \M{X}_{(1)} \right\|_1
    & = \left\| \Mn{A}{1} \M{\Lambda}
                (\Mn{A}{N} \odot \ldots \odot \Mn{A}{2})^T \right\|_1  \\
    & = \sum_{i=1}^{I_1} \sum_{r=1}^R \lambda_r a^{(1)}_{ri}  \\
    & = \sum_{r=1}^R \lambda_r
\end{align*}
The second equality uses the fact that rows of the Khatri-Rao product
sum to one when columns of  $\Mn{A}{n}$ sum to one (see the comments after
equation~(\ref{pimatrix})).

\bigskip
\section{Supplementary Material}
\bigskip

This supplementary section provides detailed results for some experiments
mentioned in the paper.  Refer to \RefSecExp of the paper for a description
of default parameters used for the algorithms and characteristics of the
workstation used for testing.

\subsection{Performance on $200 \times 1000$ Data Set}

\RefSecIntro of the paper states that a 2-way tensor of size $200 \times 1000$
was too large to factorize by the method of \RefBibZdunekCichocki but not
difficult for our algorithms.
Here we support our claim with experimental evidence.
Reference \RefBibZdunekCichockiNumber does not specify tensor details except
that the desired factorization has $R=10$ components.
We generated a synthetic data tensor according to the procedure
in \RefSecAppendix of the paper, modifying Step 1 to boost $9\%$ of the elements
to $1 + 4 R x$, where $x$ is a random number chosen from a uniform distribution
on $(0,1)$.
We tested $R=10$ and $R=50$ components.
We compare algorithms PDN-R and PQN-R using default parameters
(see \RefSecExp for values) over ten
instances of synthetic data.
Table~\ref{tbl:sup-2way-nnz} shows characteristics of the data,
and results in Table~\ref{tbl:sup-2way-results}
demonstrate that the problems are solved to high accuracy in under ten minutes.

\begin{table}[ht]
  \centering
  \small
  \caption{Sparsity of synthetic tensors versus the number of components $R$,
           for tensors of size $200 \times 1000$, generated by
           boosting $9\%$ of the elements to $1 + 4 R x$.
           Number of nonzeros is the average over ten tensors.}
      \label{tbl:sup-2way-nnz}

  \begin{tabular}{r|cc}
    $R$ & \textsc{Number Nonzeros} & \textsc{Density}  \\
    \hline
     10  &  50,144  & 25.1\%  \\
     50  &  61,826  & 30.9\%
  \end{tabular}
\end{table}

\begin{table}[htb]
  \centering
  \small
  \caption{Time to reach various stop tolerances for tensors of size
           $200 \times 1000$.
           Mean and standard deviation are reported over ten runs.
           Both methods reach the tolerance of $\tau = 10^{-4}$ in well under
           ten minutes for $R=50$ components.
          }
      \label{tbl:sup-2way-results}
 \[ \begin{array}{rl}
       \tau = 10^{-2} \;\; &
  \begin{tabular}{r|r|r}
     $R$ & \multicolumn{1}{c|}{\textsc{PDN-R}}
         & \multicolumn{1}{c}{\textsc{PQN-R}}  \\
     \hline
     10 & \, \,    63 $\pm$ \,    26 secs & \,    155 $\pm$ \,  84 secs  \\
     50 & \,      170 $\pm$ \,    31 secs & \,    311 $\pm$ \,  59 secs  \\
  \end{tabular}
 \end{array} \]
 \[ \begin{array}{rl}
       \tau = 10^{-3} \;\; &
  \begin{tabular}{r|r|r}
     $R$ & \multicolumn{1}{c|}{\textsc{PDN-R}}
         & \multicolumn{1}{c}{\textsc{PQN-R}}  \\
     \hline
     10 & \, \,    80 $\pm$ \,    25 secs & \,    184 $\pm$    100 secs  \\
     50 & \,      191 $\pm$ \,    26 secs & \,    357 $\pm$ \,  86 secs  \\
  \end{tabular}
 \end{array} \]
 \[ \begin{array}{rl}
       \tau = 10^{-4} \;\; &
  \begin{tabular}{r|r|r}
     $R$ & \multicolumn{1}{c|}{\textsc{PDN-R}}
         & \multicolumn{1}{c}{\textsc{PQN-R}}  \\
     \hline
     10 & \, \,    91 $\pm$ \,    25 secs & \,    208 $\pm$    116 secs  \\
     50 & \,      231 $\pm$ \,    59 secs & \,    385 $\pm$ \,  84 secs  \\
  \end{tabular}
 \end{array} \]
\end{table}

\vspace*{1.0in}

%\clearpage   % To separate this text from the floating table above.
%%%%%%%%%%%%%%%%%%%%%%%%%%
%\newpage
\subsection{Detail for Results in \RefSecExpFullprobPdnrPqnrMu  }
\label{sec:sup-1}
%%%%%%%%%%%%%%%%%%%%%%%%%%%

\RefSecExpFullprobPdnrPqnrMu of the paper contains a table showing the
performance of algorithms PDN-R, PQN-R, and MU on a tensor of size
$200 \times 300 \times 400$ for different values of $R$.
The table reports execution time for the algorithms to reach a stop
tolerance of $\tau = 10^{-3}$.
Supplementary Tables~\ref{tbl:sup-perf-full20} -
                     \ref{tbl:sup-score-full20-beststart}
contain additional results from the same tests.
Table~\ref{tbl:sup-perf-full20} shows that the relative advantage of
the algorithms holds for all tolerances; in particular, PDN-R and PQN-R are
faster than MU at all tolerances, and PQN-R becomes faster than PDN-R for a
sufficiently large number of components.
Table~\ref{tbl:sup-obj-full20} shows that the final values of the
objective function attained by the three algorithms are nearly identical.

\begin{table}[htb]
  \centering
  \small
  \caption{Time to reach various stop tolerances
           ($\tau = 10^{-3}$ repeats \RefTblFullprob
            in \RefSecExpFullprobPdnrPqnrMu of the paper).
           Mean and standard deviation are reported over ten runs.
           The last column for MU is the number of runs that failed to reach the
           stop tolerance after three hours (10,800 seconds) of execution.
          }
      \label{tbl:sup-perf-full20}
 \[ \begin{array}{rl}
       \tau = 10^{-2} \;\; &
  \begin{tabular}{r|r|r|rr}
    $R$ & \multicolumn{1}{c|}{\textsc{PDN-R}}
        & \multicolumn{1}{c|}{\textsc{PQN-R}}
        & \multicolumn{2}{c}{\textsc{MU}}  \\
    \hline
     20 & \,    170 $\pm$  \,  37 secs & \,    317 $\pm$    127 secs & \,   997 $\pm$ \,   347 secs & \,  (0 failures) \\
     40 & \,    339 $\pm$  \,  71 secs & \,    573 $\pm$    268 secs &     3275 $\pm$     1357 secs & \,  (0 failures) \\
     60 & \,    611 $\pm$    126 secs & \,     740 $\pm$    235 secs &     5758 $\pm$     1275 secs & \,  (0 failures) \\
     80 &      1145 $\pm$    358 secs &       1340 $\pm$    358 secs &     9134 $\pm$     1030 secs & \,  (1 failure) \, \\
    100 &      1739 $\pm$    359 secs &       1623 $\pm$    583 secs &     9687 $\pm$ \,   372 secs & \,  (6 failures) \\
  \end{tabular}
 \end{array} \]
 \[ \begin{array}{rl}
       \tau = 10^{-3} \;\; &
  \begin{tabular}{r|r|r|rr}
    $R$ & \multicolumn{1}{c|}{\textsc{PDN-R}}
        & \multicolumn{1}{c|}{\textsc{PQN-R}}
        & \multicolumn{2}{c}{\textsc{MU}}  \\
    \hline
     20 & \,    229 $\pm$ \, 57 secs & \,    397 $\pm$    123 secs &      3355 $\pm$     1933 secs & \,  (0 failures) \\
     40 & \,    493 $\pm$   151 secs & \,    818 $\pm$    185 secs &      8101 $\pm$     2045 secs & \,  (2 failures) \\
     60 &      1003 $\pm$   349 secs & \,    966 $\pm$    286 secs &      9628 $\pm$ \,   978 secs & \,  (5 failures) \\
     80 &      1682 $\pm$   642 secs &      1639 $\pm$    390 secs &  - \qquad \qquad \qquad &   (10 failures) \\
    100 &      2707 $\pm$   773 secs &      1995 $\pm$    743 secs &  - \qquad \qquad \qquad &   (10 failures) \\
  \end{tabular}
 \end{array} \]
 \[ \begin{array}{rl}
       \tau = 10^{-4} \;\; &
  \begin{tabular}{r|r|r|rr}
    $R$ & \multicolumn{1}{c|}{\textsc{PDN-R}}
        & \multicolumn{1}{c|}{\textsc{PQN-R}}
        & \multicolumn{2}{c}{\textsc{MU}}  \\
    \hline
     20 & \,    270 $\pm$ \,  67 secs & \,    451 $\pm$   108 secs &      5207 $\pm$     1975 secs & \,  (1 failure) \, \\
     40 & \,    756 $\pm$    383 secs & \,    905 $\pm$   171 secs &  - \qquad \qquad \qquad &   (10 failures) \\
     60 &      1166 $\pm$    390 secs &      1102 $\pm$   305 secs &  - \qquad \qquad \qquad &   (10 failures) \\
     80 &      1890 $\pm$    617 secs &      1859 $\pm$   413 secs &  - \qquad \qquad \qquad &   (10 failures) \\
    100 &      2834 $\pm$    741 secs &      2280 $\pm$   808 secs &  - \qquad \qquad \qquad &   (10 failures) \\
  \end{tabular}
 \end{array} \]
\end{table}

\begin{table}[htb]
  \centering
  \small
  \caption{Final objective function values reached by the algorithms.
           Algorithms minimize K-L divergence, but the negative of K-L is shown
           because this equals the (unnormalized) maximum likelihood.
           Mean and maximum (most optimal) values are reported over ten runs.
           Values are in units of $10^6$; for example, the average unnormalized
           log likelihood reached by algorithm PDN-R for $R=20$
           was $1.146 \times 10^6$.
           The likelihood increases with the number of components because this
           allows for a better fit to the data.
          }
      \label{tbl:sup-obj-full20}
  \begin{tabular}{r|rr|rr|rr}
        & \multicolumn{2}{c|}{\textsc{PDN-R}}
        & \multicolumn{2}{c|}{\textsc{PQN-R}}
        & \multicolumn{2}{c}{\textsc{MU}}  \\
    $R$ & \textsc{mean} & \textsc{max}
        & \textsc{mean} & \textsc{max}
        & \textsc{mean} & \textsc{max}  \\
     \hline
     20 &  1.146 &  1.205 &  1.144 &  1.193 &  1.145 &  1.193   \\
     40 &  1.448 &  1.475 &  1.448 &  1.475 &  1.447 &  1.473   \\
     60 &  1.610 &  1.643 &  1.611 &  1.647 &  1.610 &  1.640   \\
     80 &  1.716 &  1.735 &  1.715 &  1.734 &  1.715 &  1.738   \\
    100 &  1.794 &  1.811 &  1.794 &  1.810 &  1.794 &  1.811   \\
  \end{tabular}
\end{table}

Tables~\ref{tbl:sup-score-full20} and \ref{tbl:sup-score-full20-beststart}
report the agreement between the models computed by each algorithm and the
factor matrices from which tensor data was generated (the ``true model'').
Agreement is measured by the Tensor Toolbox {\tt score} function,
which implements the method described in~\RefBibScoreK.
The score considers the angles between all factor matrix column vectors.
Let the two $R$-component models for a three-way tensor have
weight vectors and factor matrices
$\{ \V{\lambda}^A, \Mn{A}{1}, \Mn{A}{2}, \Mn{A}{3} \}$ and
$\{ \V{\lambda}^B, \Mn{B}{1}, \Mn{B}{2}, \Mn{B}{3} \}$.  Assume columns of
factor matrices are normalized to one in the $\ell_2$ norm, adjusting
the weights as needed.  The comparison is computed from
\begin{equation}
  s = \frac{1}{R} \sum_{r=1}^R
                      [ (\V{a}^{(1)}_r)^T \V{b}^{(1)}_r ]
                      [ (\V{a}^{(2)}_r)^T \V{b}^{(2)}_r ]
                      [ (\V{a}^{(3)}_r)^T \V{b}^{(3)}_r ] .
\end{equation}
The score is the largest $s$ over column permutations of the models.
Ideally, $s$ is computed for all possible permutations, but this is not
feasible when $R$ is large.  Instead, we use the {\tt greedy} option to
limit the number of permutations considered.
The score is a number between zero and one, with one indicating
a perfect match between the two models.

Table~\ref{tbl:sup-score-full20} shows that all three algorithms produce
solutions with scores in the range of $[0.84, 0.92]$, with no clear
superiority for any algorithm.  To give a sense of score magnitudes,
the table also reports the score between computed solutions and factor matrices
used to generate other models.  The generation procedure in \RefSecAppendix
is based on random values, so we expect poor scores between a solution
and random factor matrices.  This is indeed the case, with no scores greater
than 0.01.
Table~\ref{tbl:sup-score-full20-beststart} reports the score for computed
solutions when the algorithm is initialized with the factor matrices of
the true model.  The solution should be a local minimum close to the true
model, representing an easily found solution that gives one of the best
possible fits to the data.  Since tensor data is generated with statistical
noise, a perfect fit is not seen.  These ``ideal'' scores are slightly higher
than the scores from algorithm solutions in Table~\ref{tbl:sup-score-full20}.

We notice that scores in Table~\ref{tbl:sup-score-full20-beststart}
decrease as $R$ increases.  This is because the tensor data has about the
same number of samples in each test case (see Table~\ref{tbl:sup-nnz-full20}),
but a model with more components has more variables.
As the ratio of data samples to model variables decreases, the variables are
affected more strongly by noise in the data, leading to solutions with
a lower score.

\begin{table}[htb]
  \centering
  \small
  \caption{Comparison of final computed solutions with the true model
           from which data was generated.
           The comparison uses the Tensor Toolbox {\tt score} function,
           executed with option {\tt greedy}.
           The score is a number between zero (poor match) and one (exact match).
           The \textsc{true model} column shows the lowest (worst) score obtained
           over ten runs when the computed solution is scored against the true
           model.  For comparison, the \textsc{others} column
           shows the highest (best) score between the computed solution and
           any of the nine other models, again averaged over ten runs.
          }
      \label{tbl:sup-score-full20}
  \begin{tabular}{r|cc|cc|cc}
        & \multicolumn{2}{c|}{\textsc{PDN-R}}
        & \multicolumn{2}{c|}{\textsc{PQN-R}}
        & \multicolumn{2}{c}{\textsc{MU}}  \\
    $R$ & \textsc{true model} & \textsc{others}
        & \textsc{true model} & \textsc{others}
        & \textsc{true model} & \textsc{others}  \\
     \hline
     20 &  0.919 & 0.009  &  0.845 & 0.009 &  0.854 & 0.009  \\
     40 &  0.880 & 0.010  &  0.888 & 0.010 &  0.892 & 0.010  \\
     60 &  0.846 & 0.010  &  0.873 & 0.010 &  0.844 & 0.010  \\
     80 &  0.876 & 0.010  &  0.873 & 0.010 &  0.889 & 0.010  \\
    100 &  0.865 & 0.010  &  0.864 & 0.010 &  0.847 & 0.010  \\
  \end{tabular}
\end{table}

\begin{table}[htb]
  \centering
  \small
  \caption{Comparison of final computed solutions with the true model
           from which data was generated using the Tensor Toolbox {\tt score}
           function.
           The PDN-R algorithm was initialized with
           the factor matrices of \RefSecAppendix that generated the test data.
           Columns in the table show the highest (best) and lowest (worst)
           scores obtained over ten runs.
          }
      \label{tbl:sup-score-full20-beststart}
  \begin{tabular}{r|cc}
        & \multicolumn{2}{c}{\textsc{PDN-R}}  \\
    $R$ & \textsc{best score} & \textsc{worst score}  \\
     \hline
     20 &  0.996 & 0.961   \\
     40 &  0.989 & 0.964   \\
     60 &  0.982 & 0.913   \\
     80 &  0.965 & 0.937   \\
    100 &  0.952 & 0.907   \\
  \end{tabular}
\end{table}

\clearpage   % To separate this text from the floating table above.
%%%%%%%%%%%%%%%%%%%%%%%%%%
%\newpage
\subsection{\RefSecExpFullprobPdnrPqnrMu Repeated with Different Sparsity}
%%%%%%%%%%%%%%%%%%%%%%%%%%

The data tensor of \RefSecExpFullprobPdnrPqnrMu was generated by the method
described in Step 1 of \RefSecAppendix with $20\%$ of element values boosted
to $1 + 10 R x$.
In this section, we report on convergence of the algorithms for tensors of
the same size but different sparsity.  Sparsity was modified by boosting
different numbers of elements.  The boosting procedure provides tensor data
with appropriate Poisson distributions but does not give exact control
over the number of nonzeros.  We select boost parameters that yield approximately
the same number of nonzeros for $R \in \{ 20, 40, 60, 80, 100 \}$
so we can more easily compare algorithm performance as a function of $R$.
Tables~\ref{tbl:sup-nnz-full20} - \ref{tbl:sup-nnz-full03}
describe the sparsity of the different test tensors.

\begin{table}[ht]
  \centering
  \small
  \caption{Sparsity of synthetic tensors versus the number of components $R$,
           for tensors of size $200 \times 300 \times 400$, generated by
           {\bf boosting $20\%$ of the elements} to $1 + 10 R x$.
           Number of nonzeros is the average over ten tensors.
           (This table repeats \RefTblSparsity from the paper.)}
      \label{tbl:sup-nnz-full20}

  \begin{tabular}{r|cc}
    $R$ & \textsc{Number Nonzeros} & \textsc{Density}  \\
    \hline
     20  & 413,458  &  1.72\%  \\
     40  & 450,756  &  1.88\%  \\
     60  & 464,443  &  1.94\%  \\
     80  & 470,950  &  1.96\%  \\
    100  & 475,455  &  1.98\%
  \end{tabular}
\end{table}

\begin{table}[ht]
  \centering
  \small
  \caption{Sparsity of synthetic tensors versus the number of components $R$,
           for tensors of size $200 \times 300 \times 400$, generated by
           {\bf boosting $5\%$ of the elements} to $1 + 2 R x$.
           Number of nonzeros is the average over ten tensors.}
      \label{tbl:sup-nnz-full05}

  \begin{tabular}{r|cc}
    $R$ & \textsc{Number Nonzeros} & \textsc{Density}  \\
    \hline
     20  & 158,616  &  0.66\%  \\
     40  & 141,778  &  0.59\%  \\
     60  & 148,273  &  0.62\%  \\
     80  & 161,212  &  0.67\%  \\
    100  & 177,060  &  0.74\%
  \end{tabular}
\end{table}

\begin{table}[ht]
  \centering
  \small
  \caption{Sparsity of synthetic tensors versus the number of components $R$,
           for tensors of size $200 \times 300 \times 400$, generated by
           {\bf boosting $3\%$ of the elements} to $1 + 10 R x$.
           Number of nonzeros is the average over ten tensors.}
      \label{tbl:sup-nnz-full03}

  \begin{tabular}{r|cc}
    $R$ & \textsc{Number Nonzeros} & \textsc{Density}  \\
    \hline
     20  & 55,471  &  0.23\%  \\
     40  & 44,862  &  0.19\%  \\
     60  & 47,171  &  0.20\%  \\
     80  & 51,827  &  0.22\%  \\
    100  & 57,700  &  0.24\%
  \end{tabular}
\end{table}

The experiments in \RefSecExpFullprobPdnrPqnrMu were performed for each
of these tensors.  Supplementary Table~\ref{tbl:sup-perf-full20} in the
previous section shows algorithm performance for the sparse tensor
of Table~\ref{tbl:sup-nnz-full20} ($20\%$ boost).  The tables below show
performance for the $5\%$ and $3\%$ boosted tensors.
We see that the conclusions from \RefSecExpFullprobPdnrPqnrMu hold for these
cases as well: PDN-R and PQN-R are faster than MU in nearly every case,
and PQN-R becomes faster than PDN-R for sufficiently large $R$.
Comparing the three tables with each other, we see that performance improves
when the number of nonzeros in the data decreases.

\begin{table}[htb]
  \centering
  \small
  \caption{Time to reach various stop tolerances $\tau$ on the tensor
           with the sparsity in supplementary Table~\ref{tbl:sup-nnz-full05}
           resulting {\bf from $5\%$ boost}.
           Mean and standard deviation are reported over ten runs.
           The last column for MU is the number of runs that failed to reach the
           stop tolerance after three hours of execution.
          }
      \label{tbl:sup-perf-full05}
 \[ \begin{array}{rl}
       \tau = 10^{-2} \;\; &
  \begin{tabular}{r|r|r|r}
    $R$ & \multicolumn{1}{c|}{\textsc{PDN-R}}
        & \multicolumn{1}{c|}{\textsc{PQN-R}}
        & \multicolumn{1}{c}{\textsc{MU}}  \\
    \hline
     20 & \,   165 $\pm$ \,  40 secs & \,   353 $\pm$   137 secs & \,   345 $\pm$   214 secs \\
     40 & \,   393 $\pm$    125 secs & \,   454 $\pm$   148 secs & \,   541 $\pm$   129 secs \\
     60 & \,   516 $\pm$    124 secs & \,   394 $\pm$   103 secs & \,   865 $\pm$   186 secs \\
     80 & \,   558 $\pm$ \,  48 secs & \,   386 $\pm$ \, 53 secs &     1235 $\pm$   210 secs \\
    100 & \,   658 $\pm$ \,  59 secs & \,   343 $\pm$ \, 44 secs &     1633 $\pm$   228 secs \\
  \end{tabular}
 \end{array} \]
 \[ \begin{array}{rl}
       \tau = 10^{-3} \;\; &
  \begin{tabular}{r|r|r|r}
    $R$ & \multicolumn{1}{c|}{\textsc{PDN-R}}
        & \multicolumn{1}{c|}{\textsc{PQN-R}}
        & \multicolumn{1}{c}{\textsc{MU}}  \\
    \hline
     20 & \,   270 $\pm$   139 secs & \,   580 $\pm$   323 secs & \,\,\,  no successes \,\; \\
     40 & \,   444 $\pm$   145 secs & \,   509 $\pm$   147 secs & \,\,\,  no successes \,\; \\
     60 & \,   623 $\pm$   128 secs & \,   475 $\pm$   135 secs & \,\,\,  no successes \,\; \\
     80 & \,   726 $\pm$   130 secs & \,   487 $\pm$ \, 46 secs & \,\,\,  no successes \,\; \\
    100 & \,   893 $\pm$   150 secs & \,   442 $\pm$ \, 86 secs & \,\,\,  no successes \,\; \\
  \end{tabular}
 \end{array} \]
 \[ \begin{array}{rl}
       \tau = 10^{-4} \;\; &
  \begin{tabular}{r|r|r|r}
    $R$ & \multicolumn{1}{c|}{\textsc{PDN-R}}
        & \multicolumn{1}{c|}{\textsc{PQN-R}}
        & \multicolumn{1}{c}{\textsc{MU}}  \\
    \hline
     20 & \,   287 $\pm$   139 secs & \,   585 $\pm$   320 secs & \,\,\,  no successes \,\; \\
     40 & \,   456 $\pm$   142 secs & \,   526 $\pm$   146 secs & \,\,\,  no successes \,\; \\
     60 & \,   700 $\pm$   199 secs & \,   566 $\pm$   193 secs & \,\,\,  no successes \,\; \\
     80 & \,   810 $\pm$   214 secs & \,   662 $\pm$   119 secs & \,\,\,  no successes \,\; \\
    100 & \,   964 $\pm$   173 secs & \,   631 $\pm$   187 secs & \,\,\,  no successes \,\; \\
  \end{tabular}
 \end{array} \]
\end{table}

\begin{table}[h!]
  \centering
  \small
  \caption{Time to reach various stop tolerances $\tau$ on the tensor with
           with the sparsity in supplementary Table~\ref{tbl:sup-nnz-full03}
           resulting {\bf from $3\%$ boost}.
           Mean and standard deviation are reported over ten runs.
           The last column for MU is the number of runs that failed to reach the
           stop tolerance after three hours of execution.
          } 
      \label{tbl:sup-perf-full03}
 \[ \begin{array}{rl}
       \tau = 10^{-2} \;\; &
  \begin{tabular}{r|r|r|rr}
    $R$ & \multicolumn{1}{c|}{\textsc{PDN-R}}
        & \multicolumn{1}{c|}{\textsc{PQN-R}}
        & \multicolumn{2}{c}{\textsc{MU}}  \\ 
    \hline
     20 & \,   143 $\pm$ \, 75 secs & \,   182 $\pm$ \, 73 secs & \,   108 $\pm$ \, 34 secs & (0 failures) \\
     40 & \,   141 $\pm$ \, 19 secs & \,   134 $\pm$ \, 53 secs & \,   143 $\pm$ \, 49 secs & (0 failures) \\
     60 & \,   174 $\pm$ \, 29 secs & \,   114 $\pm$ \, 23 secs & \,   218 $\pm$ \, 46 secs & (0 failures) \\
     80 & \,   219 $\pm$ \, 30 secs & \,   114 $\pm$ \, 20 secs & \,   347 $\pm$ \, 69 secs & (0 failures) \\
    100 & \,   230 $\pm$ \, 35 secs & \,   166 $\pm$ \, 37 secs & \,   481 $\pm$ \, 86 secs & (0 failures) \\
  \end{tabular}
 \end{array} \]
 \[ \begin{array}{rl}
       \tau = 10^{-3} \;\; &
  \begin{tabular}{r|r|r|rr}
    $R$ & \multicolumn{1}{c|}{\textsc{PDN-R}}
        & \multicolumn{1}{c|}{\textsc{PQN-R}}
        & \multicolumn{2}{c}{\textsc{MU}}  \\
    \hline
     20 & \,   172 $\pm$ \, 80 secs & \,   208 $\pm$ \, 77 secs & \,   278 $\pm$   135 secs & (0 failures) \\
     40 & \,   169 $\pm$ \, 40 secs & \,   151 $\pm$ \, 41 secs & \,   363 $\pm$   226 secs & (0 failures) \\
     60 & \,   328 $\pm$   122 secs & \,   139 $\pm$ \, 37 secs & \,   552 $\pm$   229 secs & (0 failures) \\
     80 & \,   290 $\pm$ \, 70 secs & \,   135 $\pm$ \, 21 secs & \,   897 $\pm$   302 secs & (0 failures) \\
    100 & \,   302 $\pm$ \, 84 secs & \,   190 $\pm$ \, 55 secs &     1320 $\pm$   397 secs & (0 failures) \\
  \end{tabular}
 \end{array} \]
 \[ \begin{array}{rl}
       \tau = 10^{-4} \;\; &
  \begin{tabular}{r|r|r|rr}
    $R$ & \multicolumn{1}{c|}{\textsc{PDN-R}}
        & \multicolumn{1}{c|}{\textsc{PQN-R}}
        & \multicolumn{2}{c}{\textsc{MU}}  \\
    \hline
     20 & \,   178 $\pm$ \, 80 secs & \,   213 $\pm$ \, 76 secs & \qquad \,\, - \qquad \,\,\,\, \, &   (10 failures) \\
     40 & \,   237 $\pm$   188 secs & \,   159 $\pm$ \, 39 secs & \qquad \,\, - \qquad \,\,\,\, \, &   (10 failures) \\
     60 & \,   342 $\pm$   127 secs & \,   176 $\pm$ \, 65 secs & \qquad \,\, - \qquad \,\,\,\, \, &   (10 failures) \\
     80 & \,   290 $\pm$ \, 70 secs & \,   189 $\pm$ \, 94 secs & \qquad \,\, - \qquad \,\,\,\, \, &   (10 failures) \\
    100 & \,   321 $\pm$ \, 98 secs & \,   259 $\pm$   129 secs & \qquad \,\, - \qquad \,\,\,\, \, &   (10 failures) \\
  \end{tabular}
 \end{array} \]
\end{table}

\clearpage   % To separate this text from the floating table above.
%%%%%%%%%%%%%%%%%%%%%%%%%%
\newpage
\subsection{Collinear Factor Matrices}
%%%%%%%%%%%%%%%%%%%%%%%%%%

Tensor data in these experiments was designed to reflect underlying factor
matrices that have nearly collinear columns.  Algorithms PDN-R and PQN-R
do much better than MU in these experiments.
The basic idea, proposed in \RefBibPhanZdunek2010,
is to generate random
factor matrices and then modify column vectors according to
\begin{equation}  \label{eq:sup-collinear}
  \V{a}^{(n)}_r = \V{a}^{(n)}_1 + \alpha \V{a}^{(n)}_r
                    \quad \mbox{for} \;\; r \in 2, \ldots, R ,
\end{equation}
with $\alpha = 0.5$.
We generated synthetic data according to \RefSecAppendix but added the
collinearity modification after Step 1.

The modification significantly impacts the sparsity of the generated tensor,
adding many nonzero elements because the boosted elements in column $r=1$
become boosted elements in all other columns.
This makes it difficult to compare performance ``before'' and ``after'' the
collinearity modification.  Instead, we made two experiments with different
values of $\alpha$ and view them as two different results.

Collinearity is measured as the cosine between pairs of column vectors.
For real-valued $N$-vectors $\V{x}$ and $\V{y}$,
\begin{equation*}
  \mbox{cos}(\V{x},\V{y}) = \frac{ \V{x}^T \V{y} }
                                 { \| \V{x} \| \| \V{y} \| } \; .
\end{equation*}
If elements of $\V{x}$ and $\V{y}$ are independent random variables chosen from
a uniform distribution on $(0,1)$, then the expected value of their cosine
is 0.75.  However, the sparse columns generated by
\RefSecAppendix are not uniform.  For instance, average collinearity among
column pairs in the experiments of supplementary Section~\ref{sec:sup-1}
is close to $0.10$.

Tables~\ref{tbl:sup-nnz-collinear-pt5} - \ref{tbl:sup-perf-coll-pt1}
show tensor characteristics and computational performance
for two experiments.  The tensors in the experiments have the same size but
use a different value of $\alpha$ in equation (\ref{eq:sup-collinear})
to modify collinearity, which also changes their sparsity; hence, the two
experiments should not be compared with each other.  What they each show
is that algorithms PDN-R and PQN-R are much faster than MU, especially
when high accuracy is desired.

\begin{table}[p!]
  \centering
  \small
  \caption{Sparsity and collinearity of synthetic tensors versus the number 
           of components $R$, for tensors of size $50 \times 50 \times 50$,
           generated by boosting $10\%$ of the elements to $1 + 10 R x$,
           and modifying collinearity using $\alpha = 0.5$ in
           equation (\ref{eq:sup-collinear}).
           \textsc{All Pairs} shows the average collinearity
           between all pairs of column vectors within each factor matrix.
           Data under $\V{a}_1$ \textsc{Pairs} shows the average
           collinearity between the $\V{a}^{(n)}_1$ column vector
           and all others within each factor matrix.
           All measurements are an average over ten tensors.
          }
      \label{tbl:sup-nnz-collinear-pt5}

  \begin{tabular}{r|cc|cc}
        & \multicolumn{2}{c}{\textsc{Sparsity}}
        & \multicolumn{2}{c}{\textsc{Collinearity}}  \\
    $R$ & \textsc{Number Nonzeros} & \textsc{Density}
        & \textsc{All Pairs} & $\V{a}_1$ \textsc{Pairs}  \\
    \hline
     10  & 16,442  &  13.2\%   &  0.839  & 0.900  \\
     20  & 14,908  &  11.9\%   &  0.828  & 0.897  \\
     30  & 15,572  &  12.5\%   &  0.830  & 0.900  \\
     40  & 16,892  &  13.5\%   &  0.800  & 0.882
  \end{tabular}
\end{table}

\begin{table}[hb]
  \centering
  \small
  \caption{Time to reach various stop tolerances $\tau$ on tensors of size
           $50 \times 50 \times 50$, generated by boosting $10\%$ of the elements
           to $1 + 10 R x$, and modifying collinearity
           using $\alpha = 0.5$ in equation (\ref{eq:sup-collinear}).
           Mean and standard deviation are reported over ten runs.
           The last column for MU is the number of runs that failed to reach
           the stop tolerance after 1000 seconds of execution.
          }
      \label{tbl:sup-perf-coll-pt5}
 \[ \begin{array}{rl}
       \tau = 10^{-2} \;\; &
  \begin{tabular}{r|r|r|rr}
    $R$ & \multicolumn{1}{c|}{\textsc{PDN-R}}
        & \multicolumn{1}{c|}{\textsc{PQN-R}}
        & \multicolumn{2}{c}{\textsc{MU}}  \\
    \hline
     10 & \, \, \,    9 $\pm$ \, \,   3 secs & \,     14 $\pm$ \, \,   3 secs & \, 26 $\pm$ \, \, 7 secs & \, (0 failures)  \\
     20 & \, \,     18 $\pm$ \, \,    7 secs & \,     23 $\pm$ \, \,   6 secs & \, 75 $\pm$ \, 35 secs & \, (0 failures)  \\
     30 & \, \,     20 $\pm$ \, \,    5 secs & \,     35 $\pm$ \,     14 secs & \, 97 $\pm$ \, 36 secs & \, (0 failures)  \\
     40 & \, \,     28 $\pm$ \, \,    6 secs & \,     46 $\pm$ \,     17 secs &   191 $\pm$ \, 84 secs & \, (0 failures)  \\
  \end{tabular}
 \end{array} \]
 \[ \begin{array}{rl}
       \tau = 10^{-3} \;\; &
  \begin{tabular}{r|r|r|rr}
    $R$ & \multicolumn{1}{c|}{\textsc{PDN-R}}
        & \multicolumn{1}{c|}{\textsc{PQN-R}}
        & \multicolumn{2}{c}{\textsc{MU}}  \\
    \hline
     10 & \, \,     20 $\pm$ \,    11 secs & \,     16 $\pm$ \, \,   4 secs & 223 $\pm$ 126 secs & \, (0 failures)  \\
     20 & \, \,     25 $\pm$ \,    11 secs & \,     25 $\pm$ \, \,   6 secs & 511 $\pm$ 191 secs & \, (1 failure) \,  \\
     30 & \, \,     25 $\pm$ \, \,   6 secs & \,     38 $\pm$ \,    16 secs & 719 $\pm$ 117 secs & \, (2 failures)  \\ 
     40 & \, \,     35 $\pm$ \, \,   7 secs & \,     59 $\pm$ \,    26 secs & 734 $\pm$ 172 secs & \, (4 failures)  \\
  \end{tabular}
 \end{array} \]
 \[ \begin{array}{rl}
       \tau = 10^{-4} \;\; &
  \begin{tabular}{r|r|r|rr}
    $R$ & \multicolumn{1}{c|}{\textsc{PDN-R}}
        & \multicolumn{1}{c|}{\textsc{PQN-R}}
        & \multicolumn{2}{c}{\textsc{MU}}  \\
    \hline
     10 & \, \,     32 $\pm$ \,    17 secs & \,     17 $\pm$ \, \,   5 secs & 295 $\pm$ 113 secs & \, (1 failure) \,  \\
     20 & \, \,     28 $\pm$ \,    12 secs & \,     26 $\pm$ \, \,   6 secs & 784 $\pm$ 221 secs & \, (8 failures)  \\
     30 & \, \,     28 $\pm$ \, \,  6 secs & \,     40 $\pm$ \,     17 secs &  - \qquad \qquad \qquad &   (10 failures) \\
     40 & \, \,     46 $\pm$ \,    21 secs & \,     63 $\pm$ \,     26 secs &  - \qquad \qquad \qquad &   (10 failures) \\
  \end{tabular}
 \end{array} \]
\end{table}

%\clearpage   % To separate this table from the previous one.
\begin{table}[pt!]
  \centering
  \small
  \caption{Sparsity and collinearity of synthetic tensors versus the number 
           of components $R$, for tensors of size $50 \times 50 \times 50$,
           generated by boosting $10\%$ of the elements to $1 + 10 R x$,
           and modifying collinearity using $\alpha = 0.1$ in
           equation (\ref{eq:sup-collinear}).
           \textsc{All Pairs} shows the average collinearity
           between all pairs of column vectors within each factor matrix.
           Data under $\V{a}_1$ \textsc{Pairs} shows the average
           collinearity between the $\V{a}^{(n)}_1$ column vector
           and all others within each factor matrix.
           All measurements are an average over ten tensors.
          }
      \label{tbl:sup-nnz-collinear-pt1}

  \begin{tabular}{r|cc|cc}
        & \multicolumn{2}{c}{\textsc{Sparsity}}
        & \multicolumn{2}{c}{\textsc{Collinearity}}  \\
    $R$ & \textsc{Number Nonzeros} & \textsc{Density}
        & \textsc{All Pairs} & $\V{a}_1$ \textsc{Pairs}  \\
    \hline
     10  & 9,432  &  7.55\%   &  0.991  & 0.995  \\
     20  & 8,572  &  6.86\%   &  0.990  & 0.995  \\
     30  & 8,828  &  7.06\%   &  0.991  & 0.995  \\
     40  & 9,419  &  7.54\%   &  0.988  & 0.993
  \end{tabular}
\end{table}

\begin{table}[hb]
  \centering
  \small
  \caption{Time to reach various stop tolerances $\tau$ on tensors of size
           $50 \times 50 \times 50$, generated by boosting $10\%$ of the elements
           to $1 + 10 R x$, and modifying collinearity
           using $\alpha = 0.1$ in equation (\ref{eq:sup-collinear}).
           Mean and standard deviation are reported over ten runs.
           The last column for MU is the number of runs that failed to reach
           the stop tolerance after 1000 seconds of execution.
          }
      \label{tbl:sup-perf-coll-pt1}
 \[ \begin{array}{rl}
       \tau = 10^{-2} \;\; &
  \begin{tabular}{r|r|r|rr}
    $R$ & \multicolumn{1}{c|}{\textsc{PDN-R}}
        & \multicolumn{1}{c|}{\textsc{PQN-R}}
        & \multicolumn{2}{c}{\textsc{MU}}  \\
    \hline
     10 & \, \,     11 $\pm$ \, \,   4 secs & \,     16 $\pm$ \, \,   4 secs & \, \,     33 $\pm$ \, \, \,   8 secs & \,  (0 failures) \\
     20 & \, \,     12 $\pm$ \, \,   3 secs & \,     23 $\pm$ \, \,   8 secs & \, \,     99 $\pm$ \, \,    31 secs & \,  (0 failures) \\
     30 & \, \,     14 $\pm$ \, \,   3 secs & \,     30 $\pm$ \, \,   6 secs & \,      226 $\pm$ \, \,    66 secs & \,  (0 failures) \\
     40 & \, \,     17 $\pm$ \, \,   4 secs & \,     30 $\pm$ \, \,   7 secs & \,      340 $\pm$ \,     114 secs & \,  (0 failures) \\
  \end{tabular}
 \end{array} \]
 \[ \begin{array}{rl}
       \tau = 10^{-3} \;\; &
  \begin{tabular}{r|r|r|rr}
    $R$ & \multicolumn{1}{c|}{\textsc{PDN-R}}
        & \multicolumn{1}{c|}{\textsc{PQN-R}}
        & \multicolumn{2}{c}{\textsc{MU}}  \\
    \hline
     10 & \, \,     19 $\pm$ \, \,   9 secs & \,     26 $\pm$ \, \,   6 secs & \,      181 $\pm$ \, \,    62 secs & \,  (0 failures) \\
     20 & \, \,     22 $\pm$ \, \,   9 secs & \,     32 $\pm$ \, \,   7 secs & \,      512 $\pm$ \,     200 secs & \,  (2 failures) \\
     30 & \, \,     26 $\pm$ \, \,   4 secs & \,     41 $\pm$ \, \,   6 secs & \,      787 $\pm$ \,     174 secs & \,  (7 failures) \\
     40 & \, \,     30 $\pm$ \, \,   7 secs & \,     43 $\pm$ \, \,   8 secs & \,      827 $\pm$ \, \, \,   0 secs & \,  (9 failures) \\
  \end{tabular}
 \end{array} \]
 \[ \begin{array}{rl}
       \tau = 10^{-4} \;\; &
  \begin{tabular}{r|r|r|rr}
    $R$ & \multicolumn{1}{c|}{\textsc{PDN-R}}
        & \multicolumn{1}{c|}{\textsc{PQN-R}}
        & \multicolumn{2}{c}{\textsc{MU}}  \\
    \hline
     10 & \, \,     24 $\pm$ \,    14 secs & \,     28 $\pm$ \, \,   6 secs & \,      545 $\pm$ \,     228 secs & \,  (2 failures) \\
     20 & \, \,     25 $\pm$ \, \,   8 secs & \,     36 $\pm$ \,    11 secs &  - \qquad \qquad \qquad &   (10 failures) \\
     30 & \, \,     30 $\pm$ \, \,   5 secs & \,     43 $\pm$ \, \,   7 secs &  - \qquad \qquad \qquad &   (10 failures) \\
     40 & \, \,     33 $\pm$ \, \,   6 secs & \,     45 $\pm$ \, \,   8 secs &  - \qquad \qquad \qquad &   (10 failures) \\
  \end{tabular}
 \end{array} \]
\end{table}

\end{document}